\documentclass[a4paper, 12pt, final]{article}

\usepackage[latin1]{inputenc}
\usepackage[T1]{fontenc}
\usepackage[english]{babel}
\usepackage{fullpage}
\usepackage{amsmath}
\usepackage{dsfont}\let\mathbb\mathds
\usepackage{setspace}
\usepackage{amsfonts}
\usepackage{amssymb}
\usepackage{graphicx}
\usepackage{amscd}
\usepackage{mathrsfs}
\usepackage{fancybox}
\usepackage{mathtools}

	\makeatletter
	\renewenvironment{thebibliography}[1]
	{\section*{\refname}%
		\@mkboth{\MakeUppercase\refname}{\MakeUppercase\refname}%
		\list{\@biblabel{\@arabic\c@enumiv}}%
		{\settowidth\labelwidth{\@biblabel{#1}}%
			\leftmargin\labelwidth
			\advance\leftmargin\labelsep
			\@openbib@code
			\usecounter{enumiv}%
			\let\p@enumiv\@empty
			\itemsep=0pt
			\parsep=0pt
			\leftmargin=\parindent
			\itemindent=-\parindent
			\renewcommand\theenumiv{\@arabic\c@enumiv}}%
		\sloppy
		\clubpenalty4000
		\@clubpenalty \clubpenalty
		\widowpenalty4000%
		\sfcode`\.\@m}
	{\def\@noitemerr
		{\@latex@warning{Empty `thebibliography' environment}}%
		\endlist}
	\makeatother

\usepackage{lmodern}
\usepackage{graphicx}
\usepackage{fancyheadings}
\usepackage{chngcntr}
\counterwithout{figure}{section}
\usepackage{tikz}
\usetikzlibrary{calc}
\usepackage{pgfplots}
\usepgfplotslibrary{groupplots}
\usepgfplotslibrary{units}

\usetikzlibrary{arrows}
\usepackage{tkz-berge} 

\usetikzlibrary{graphs,shapes.geometric,arrows,fit,matrix,positioning}

\usepackage[extra]{tipa}

\usepackage[OT2,T1]{fontenc}

\usepackage{amsmath, amsthm, amsfonts}

\DeclareMathOperator{\conv}{conv}

\DeclareMathOperator{\Lie}{Lie}
\DeclareMathOperator{\weight}{weight}

\DeclareMathOperator{\id}{id}

\DeclareMathOperator{\pr}{pr}

\DeclareMathOperator{\Map}{Map}

\DeclareMathOperator{\depth}{depth}
\DeclareMathOperator{\gr}{gr}

\DeclareMathOperator{\DR}{DR}
\DeclareMathOperator{\rig}{rig}
\DeclareMathOperator{\Spec}{Spec}
\DeclareMathOperator{\B}{B}
\DeclareMathOperator{\an}{an}

\DeclareMathOperator{\Ad}{Ad}
\DeclareMathOperator{\Isom}{Isom}

\DeclareMathOperator{\dch}{dch}

\DeclareMathOperator{\Gal}{Gal}

\DeclareMathOperator{\shft}{shft}

\DeclareMathOperator{\crys}{crys}
\DeclareMathOperator{\Col}{Col}

\DeclareMathOperator{\har}{har}
\DeclareMathOperator{\loc}{loc}

\DeclareMathOperator{\mot}{mot}
\DeclareMathOperator{\KZ}{KZ}

\DeclareMathOperator{\Vect}{Vect}

\DeclareMathOperator{\comp}{comp}
\DeclareMathOperator{\un}{un}

\DeclareMathOperator{\iter}{iter}

\DeclareMathOperator{\DMR}{DMR}
\DeclareMathOperator{\MT}{MT}
\DeclareMathOperator{\dec}{dec}
\DeclareMathOperator{\Li}{Li}
\DeclareMathOperator{\Aut}{Aut}

\cornersize{2}
\theoremstyle{definition}

\newtheorem{Théorème}{Théorème}[section]
\newtheorem{Proposition}[Théorème]{Proposition}
\newtheorem{Theorem}[Théorème]{Theorem}

\newtheorem{Notation}[Théorème]{Notation}
\newtheorem{Conjecture}[Théorème]{Conjecture}

\newtheorem{Question}[Théorème]{Question}

\newtheorem{Nota Bene}[Théorème]{Nota Bene}

\newtheorem{Principle}[Théorème]{Principle}
\newtheorem{Key Lemma}[Théorème]{Key Lemma}

\newtheorem{Fact}[Théorème]{Fact}

\newtheorem{Example}[Théorème]{Example}

\newtheorem{Convention-Notation}[Théorème]{Convention-Notation}
\newtheorem{Definition}[Théorème]{Definition}

\newtheorem{Proposition-Definition}[Théorème]{Proposition-Definition}

\newtheorem{N.B.}[Théorème]{N.B.}
\newtheorem{Convention}[Théorème]{Convention}

\usepackage{lipsum}
	
\input cyracc.def
\DeclareFontFamily{U}{russian}{}
\DeclareFontShape{U}{russian}{m}{n}
        { <5><6> wncyr5
        <7><8><9> wncyr7
        <10><10.95><12><14.4><17.28><20.74><24.88> wncyr10 }{}
\DeclareSymbolFont{Russian}{U}{russian}{m}{n}
\DeclareSymbolFontAlphabet{\mathcyr}{Russian}
\makeatletter
\let\@math@cyr\mathcyr
\renewcommand{\mathcyr}[1]{\@math@cyr{\cyracc #1}}
\makeatother
\newcommand{\sh}{\mathcyr{sh}} 

\date{}

\newcommand{\simlra}{\buildrel \sim \over \longrightarrow}
\makeatletter
\newcounter {subsubsubsection}[subsubsection]
\renewcommand\thesubsubsubsection{\thesubsubsection .\@alph\c@subsubsubsection}
\newcommand\subsubsubsection{\@startsection{subsubsubsection}{4}{\z@}%
                                     {-3.25ex\@plus -1ex \@minus -.2ex}%
                                     {1.5ex \@plus .2ex}%
                                     {\normalfont\normalsize\bfseries}}
\newcommand*\l@subsubsubsection{\@dottedtocline{3}{10.0em}{4.1em}}
\newcommand*{\subsubsubsectionmark}[1]{}
\makeatother

\usepackage[top=3cm, bottom=3cm, left=2.5cm, right=2.5cm]{geometry}

\title{$p$-adic multiple zeta values and $p$-adic pro-unipotent harmonic actions : summary of parts I and II}

\author{David Jarossay}

\begin{document}

\maketitle

\begin{abstract}
This is a review on the two first parts of our work on $p$-adic multiple zeta values at $N$-th roots of unity ($p$MZV$\mu_{N}$'s), the $p$-adic periods of the crystalline pro-unipotent fundamental groupoid of $\mathbb{P}^{1} - \{0,\mu_{N},\infty\}$ (where $N$ and $p$ are coprime). We restrict for simplicity the review to the case of $N=1$, i.e. the case of $p$-adic multiple zeta values ($p$MZV's).
\newline The main tools are new objects which we call $p$-adic pro-unipotent harmonic actions. These are continuous group actions on a space containing the non-commutative generating series of weighted multiple harmonic sums, they are related to the motivic Galois action on $\pi_{1}^{\un}(\mathbb{P}^{1} - \{0,1,\infty\})$ and to the Poisson-Ihara bracket, and interrelated by some maps. They are defined in \cite{J2} and \cite{J3} ; the definition relies on a simplification of the differential equation of the Frobenius, proved as a preliminary technical fact by \cite{J1}.
\newline Part I (\cite{J1},\cite{J2},\cite{J3}) is an explicit computation of the Frobenius of $\pi_{1}^{\un,\crys}(\mathbb{P}^{1} - \{0,1,\infty\})$, and in particular of $p$MZV's.
We give formulas which keep a track of the motivic Galois action.
\newline Part II (\cite{J4},\cite{J5},\cite{J6}) is a study of the algebraic properties of $p$MZV's brought together with the formulas of part I. We state an explicit elementary version of the Galois theory of $p$MZV's.
\newline In this text, we emphasize the ideas and the intuition of this work. We review the general context and the motivations (\S1), a technical description of $\pi_{1}^{\un}(\mathbb{P}^{1} - \{0,1,\infty\})$ (\S2), and we explain our general strategy (\S3). We state most of the main results of part I (\S4) and of part II (\S5), and we also summarize and motivate the methods of the proofs. We conclude on the main messages of this work (\S6).
\end{abstract}


\normalsize

\newpage

\section{Context and motivation}

We introduce the two general notions that we will study : periods, and the pro-unipotent fundamental groupoid (\S1.1, \S1.2), then the particular example of $\pi_{1}^{\un}(\mathbb{P}^{1} - \{0,1,\infty\})$ and multiple zeta values, and its $p$-adic aspects (\S1.3, \S1.4).

\subsection{Periods}

\subsubsection{Definition}

There are several possible ways to define periods.

\begin{Definition} (See \cite{Gr1}, \cite{Hu-MS})
\footnote{This definition is not the most general one.} \label{first definition of periods}Let $X$ be a smooth algebraic variety over $\mathbb{Q}$, and $D$ a normal crossings divisor. One can associate to them :
\begin{itemize} \item The algebraic De Rham cohomology $H^{\DR}(X,D)$ of $X$ relative to $D$, which is a graded $\mathbb{Q}$-vector space of finite dimension
\item The rational Betti cohomology $H^{B}(X,D)$ of $X(\mathbb{C})$ relative to $D(\mathbb{C})$, which is a graded $\mathbb{Q}$-vector space of finite dimension
\item The isomorphism of comparison $\comp : H^{\DR}(X,D) \otimes_{k} \mathbb{C} \simlra H^{B}(X,D)^{\vee} \otimes_{\mathbb{Q}} \mathbb{C}$, defined by $(\omega,\Delta) \mapsto \int_{\Delta}\omega$.
\end{itemize}
\noindent The coefficients of any matrix of the isomorphism of comparison associated with a couple of bases of $H^{\DR}(X,D)$ and $H^{B}(X,D)^{\vee}$ over $\mathbb{Q}$ are called the periods of $(X,D)$.
\end{Definition}

\begin{Definition} \label{second definition of periods} (Kontsevich-Zagier, \cite{KoZ}, \S1.1) A period is a complex number whose real and imaginary parts are values of absolutely convergent integrals of rational functions with rational coefficients, over domains in $\mathbb{R}^{n}$ given by polynomial inequalities with rational coefficients.
\end{Definition}

\begin{Example} All algebraic numbers are periods, and one can replace above "rational" by "algebraic". $2i\pi = \int_{\gamma} \frac{dz}{z}$, where $\gamma$ is a simple counterclockwise loop around $0$ in $\mathbb{C}$, is a period of $\mathbb{A}^{1} - \{0\}$. For $r \in \mathbb{Q}^{+\ast}$, $\log(r)= \int_{1}^{r} dt/t$ is a period. Conjecturally, $e=\exp(1)$ is not a period.
\end{Example}

\noindent The set of periods is countable. In particular, there are complex numbers that are not periods.
\newline 
\newline There exists similarly a notion of $p$-adic periods : it arises from the comparison between De Rham cohomology and étale cohomology, resp. crystalline cohomology and étale cohomology \cite{Fa1}, \cite{Fa2} ; $p$-adic periods are elements of Fontaine rings \cite{Fo}. In this text, we will consider only the image in $\mathbb{Q}_{p}$ of certain elements of Fontaine rings by the reduction maps (see \S1.4.1) and study them through the concept of pro-unipotent fundamental groupoid (reviewed in \S1.2) thus the notion of $p$-adic periods is not technically needed for our purposes.

\subsubsection{The conjecture of periods and motivic Galois actions}

The theory of motives, initiated by Grothendieck, aims for a unification of the cohomology theories of algebraic geometry. Several different categories of motives have been constructed so far (see \cite{An1}, \cite{An2} for reviews). 
\newline Let us assume that we have such a category of motives $M$, with the following properties. First, for a $(X,D)$ as above, the triple $(H^{\DR}(X,D),H^{B}(X,D),\comp)$ lifts to an object $H^{M}(X,D)$ in $M$, and one has functors of "Betti realization" and "De Rham realization" $\omega_{B}$, $\omega_{\DR}$ from $M$ to the category of graded finite dimensional vector spaces over $\mathbb{Q}$, which send $H^{M}(X,D) \mapsto H^{\DR}(X,D)$ and $H^{M}(X,D) \mapsto H^{\B}(X,D)$. Let us assume moreover that $\omega_{\DR}$ and $\omega_{B}$ make $M$ into a neutral Tannakian over $\mathbb{Q}$ (see \cite{De3}, or \cite{An1} \S2, for generalities on Tannakian categories \footnote{A $\otimes$-category over a commutative ring $F$ is a category $\mathcal{T}$ equipped with a functor $\otimes : \mathcal{T} \times \mathcal{T} \rightarrow \mathcal{T}$, and a unit object $1$, which is "associative, commutative, unitary" ; a $\otimes$-cateory is said to be rigid when it has a functor $\mathcal{T} \rightarrow \mathcal{T}^{op}$ called autoduality which satisfies certain axioms (\cite{An1}, \S2.2.2).
	\newline If $\mathcal{T}$ be an abelian rigid $\otimes$-category over a field $F$, a fiber functor over $\mathcal{T}$ is a faithful and exact $\otimes$-functor $\mathcal{T} \rightarrow \Vect(K)$, where $K$ is a field extension of $F$. If a fiber functor exists, $\mathcal{T}$ is said to be Tannakian and the $K$-affine group scheme $G = \underline{\Aut}^{\otimes}(\omega)$ is called the Tannakian group of $\mathcal{T}$ attached to $\omega$. It is compatible with base-change by field extensions of $K$.
	\newline If there exists a fiber functor with $K=F$, $\mathcal{T}$ is said to be neutral Tannakian, and $\omega$ defines an equivalence of $\otimes$-rigid categories $\mathcal{T} \rightarrow Rep_{F}(G)$, where $Rep_{F}(G)$ is the category of finite dimensional representations of $G$. (\cite{An1}, \S2.3.2).}) 
: $M$ is then equivalent to the category of representations of an affine algebraic group over $\mathbb{Q}$ : $G_{B} = \underline{\Aut}^{\otimes}(\omega_{B})$, resp. $G_{DR} = \underline{\Aut}^{\otimes}(\omega_{\DR})$. The scheme $P = \underline{\Isom}^{\otimes}(\omega_{DR},\omega_{B})$ is a torsor over $G_{B}$ and $G_{\DR}$. 
$G_{B}$ and $G_{\DR}$ are called \emph{motivic Galois groups}, and their action on $P$ is called a \emph{motivic Galois action} ; the isomorphism $\comp$ is a point of $P$ over $\mathbb{C}$.
\newline The important \emph{conjecture of periods} of Grothendieck (unpublished), is that every polynomial equation in a $\mathbb{Q}$-algebra of periods "is of geometric origin" i.e. is reflected the framework above. More precisely, 

\begin{Conjecture} (Grothendieck, see \cite{An1}, \cite{An2}, \cite{Hu-MS} for precise statements) In a context as above, $\comp$ is a generic point of the torsor of periods.
\end{Conjecture}

\noindent Variant : the previous facts imply that the transcendence degree of the extension of $\mathbb{Q}$ generated by the periods of $(X,D)$ is bounded by the dimension of the motivic Galois group and :

\begin{Conjecture} (Grothendieck, see \cite{An1}, \cite{An2}, \cite{Hu-MS} for precise statements) In a context as above, the transcendence degree of the extension of $\mathbb{Q}$ generated by the periods of $(X,D)$ is equal to the dimension of the motivic Galois group.
\end{Conjecture}

\noindent A different and more elementary formulation of the conjecture of periods is as follows (Kontsevich-Zagier, \cite{KoZ}, \S4.1). Let $P_{formal}$, the algebra of "formal periods", be the algebra generated the equivalence classes of quadruples $(X,D,\omega,\gamma)$ where $(X,D)$ is as in Definition \ref{first definition of periods}, $\omega \in \Omega^{d}(X)$, $\gamma \in H_{d}(X(\mathbb{C}),D(\mathbb{C}),\mathbb{Q})$, $d= \dim(X)$, where the equivalence is defined by the operations of linearity of integration, changes of variables and Stokes formula. Let us denote by $P_{num} \subset \mathbb{C}$ the algebra of ("numerical") periods of Definition \ref{second definition of periods}.

\begin{Conjecture} \label{conjecture of periods Kontsevich Zagier}
	(Kontsevich-Zagier, \cite{KoZ}, \S4.1) The map from $P_{formal} \rightarrow P_{num}$ sending the equivalence class of $(X,D,\omega,\gamma)$ to $\int_{\gamma} \omega$ is injective.
\end{Conjecture}

\subsection{The pro-unipotent fundamental groupoid}

The pro-unipotent fundamental groupoid $\pi_{1}^{\un}$ is the collection of several functors, which send certain algebraic varieties to certain groupoids in pro-affine schemes over them : $\pi_{1}^{\un,\B}$ (Betti), $\pi_{1}^{\un,\DR}$ (De Rham), $\pi_{1}^{\un,\crys}$ (crystalline), $\pi_{1}^{\un,l}$ ($l$-adic) ; some morphisms comparisons between them  ; and, finally, a functor $\pi_{1}^{\un,\mot}$ (motivic) whose realizations (in the sense of the realizations of a motive) are the previous ones. One also has the Hodge realization which englobes the Betti and De Rham realizations and the comparison between them. The notion of $\pi_{1}^{\un}$ has been defined by Deligne \cite{De2} and has been extended by Goncharov \cite{Go2} and by Deligne and Goncharov \cite{DeGo}.

\subsubsection{\label{Betti and De Rham pi1 un}The Betti and De Rham realizations and the comparison}

The concept of $\pi_{1}^{\un}$ relies on the Tannakian point of view on the topological fundamental groupoid, and the Riemann-Hilbert correspondence (see for example \cite{De1}, and also \cite{De2}, \S10.10-\S10.13).

\begin{Definition} \label{De Rham definition} (Deligne, \cite{De2}, \S13.5, \S10.30 ii)) Let $X$ be a smooth algebraic variety over a field $K$ of characteristic zero ; let us view $X$ as $\overline{X}-D$, where $\overline{X}$ is a projective smooth variety over $K$ and $D$ is a normal crossings divisor ; then $\pi_{1}^{\un,\DR}(X)$ is the fundamental groupoid associated with the Tannakian category $C^{\un,\DR}(X)$ of vector bundles with integrable connection, with logarithmic singularities at $D$, and which are unipotent.
\end{Definition}

\begin{Definition} (Deligne, \cite{De2}, \S13.6)\label{Betti definition} Same hypothesis plus the existence of an embedding $k\hookrightarrow \mathbb{C}$. Then $\pi_{1}^{\un,\B}(X)$ (relative to this embedding) is the algebraic unipotent envelope of the topological fundamental groupoid of $X(\mathbb{C})$, in the sense of \cite{De2}, \S10.24.
\end{Definition}

\begin{Theorem} \label{Betti De Rham comparison}(Deligne, \cite{De2}, equation (10.33.4) and \S10.43
	\footnote{$\pi_{1}^{\un}$ is compatible with base-change (\cite{De2}, \S10.36-\S10.43) ; note that this is not true for its variant defined without the condition of unipotence of the bundles (\cite{De2}, \S10.35)}) Under the previous hypothesis, the Riemann-Hilbert correspondence implies an isomorphism 
	$\pi_{1}^{\un,\DR}(X) \otimes \mathbb{C} \simeq \pi_{1}^{\un,\B}(X) \otimes \mathbb{C}$.
\end{Theorem}

\subsubsection{The $p$-adic aspects \label{paragraph crystalline Frobenius}}

In the $p$-adic world, one has Tannakian definitions similar to Definition \ref{De Rham definition}, where "bundle with connection" is replaced by "$F$-isocristal" in the sense of \cite{BerO}, \cite{O1}, or \cite{Kat}. Let $k$ be a perfect field of characteristic $p>0$. By Shiho \cite{Sh1}, \cite{Sh2}, if $(X,M)$ is a fine log scheme over $k$ satisfying certain conditions (\cite{Sh1}, \S4) one has a notion of (log-)crystalline pro-unipotent fundamental groupoid $\pi_{1}^{\un,\crys}(X_{0})$. Alternatively, by Chiarellotto and Le Stum \cite{ChiL}, if $X_{0}$ is separated scheme of finite type over $k$, one has a notion of rigid pro-unipotent fundamental groupoid $\pi_{1}^{\un,\rig}(X_{0})$.
\newline Finally, by Deligne (\cite{De2}, \S11) one has a notion of crystalline Frobenius of the De Rham $\pi_{1}^{\un}$. Let us review briefly \cite{De2}, \S11.9 - \S11.12. Let $k$ be a perfect field of characteristic $p>0$, $W(k)$ its ring of Witt vectors, $K$ the field of fractions of $W(k)$ ; let $S=\Spec(W(k))$. Let $X$ be a smooth algebraic variety over $S$, $X_{K}= X \times_{S} \Spec(K)$ ; assume that $X = \overline{X} - D$, with $\overline{X}$ smooth over $S$ and $D$ a relative normal crossings divisor, sum of smooth divisors. Let $\sigma$ be the Frobenius automorphism of $W(k)$, $F = \sigma^{\ast} : S \rightarrow S$, $X^{(p)}$ the pull-back of $X$ by $F$, and $X_{K}^{(p)} = X^{(p)} \times_{S} \Spec(K)$. One has a notion of pull-back by Frobenius $F^{\ast}\mathcal{V}$ of certain bundles with integrable connection $\mathcal{V}$ over $X$ (\cite{De2} Definition 11.10). It gives rise to a $K$-linear functor $F_{X/K}^{\ast} : C^{\un,\DR}(X_{K}^{(p)}) \rightarrow C^{\un,\DR}(X_{K})$ (\cite{De2}, equation (11.11.1)), where $C^{\un,\DR}$ is in the sense of Definition \ref{De Rham definition}. There exists an isomorphism of groupoids, which is horizontal with respect to the connexions ($\pi_{1}^{\un,\DR}$ is an initial object of the category of groupoids with integrable connection : \cite{De2}, \S10.49, in the sense of \cite{De2}, \S10.28) : \cite{De2}, equation (11.11.2)
$$ F_{X/K\ast} : \pi_{1}^{\un,\DR}(X_{K}) \simlra F_{X/K}^{\ast}(\pi_{1}^{\un,\DR}(X_{K}^{(p)})) $$

\begin{Definition} \label{def of Frobenius}(Deligne, \cite{De2}, \S13.6) Under the assumptions above, the crystalline Frobenius of $\pi_{1}^{\un,\DR}(X_{K})$ is $\phi=(F_{X/K})_{\ast}^{-1}$.
\end{Definition}

\noindent Thus, the crystalline $\pi_{1}^{\un}$ of  $X_{0}= X \times_{\Spec(W(k))} \Spec(k)$ is $\pi_{1}^{\un,\DR}(X_{K})$ equipped with $\phi$. We will follow this point of view in this text. By contrast, in the two other points of view, there is a theorem of comparison involving $\pi_{1}^{\un,\DR}(X_{K})$ when $X_{0}$ is liftable.

\subsubsection{Relation with cohomology and periods}

A theorem of Beilinson (unpublished, reviewed in \cite{Go1}, \S4, Theorem 4.1 and \cite{DeGo}, \S3, Proposition 3.4) expresses $\pi_{1}^{\un,\B}(X)$ in terms of Betti cohomology groups of $X(\mathbb{C})^{n}$ relative to certain divisors for all $n\in \mathbb{N}^{\ast}$. This theorem has variants for the other realizations of the $\pi_{1}^{\un}$, which are compatible with usual structures on cohomology (\cite{Go1}, \S4). This is the starting point of the definition of the motivic $\pi_{1}^{\un}$ (\cite{Go1}, \S4 and \cite{DeGo}, \S3-\S4). Thus, the $\pi_{1}^{\un}$ is actually a very specific part of the cohomology of algebraic varieties. When $X$ is defined over $k \subset \overline{\mathbb{Q}} \subset \mathbb{C}$, it thus makes sense, by \S\ref{Betti and De Rham pi1 un}, to talk about the periods of the Betti-De Rham comparison of $\pi_{1}^{\un}(X)$.

\subsection{$\pi_{1}^{\un}(\mathbb{P}^{1} - \{0,1,\infty\})$ and multiple zeta values}

\subsubsection{$\mathbb{P}^{1} - \{0,1,\infty\}$ and the motivation for the $\pi_{1}^{\un}$}

Despite the generality of the definition of the $\pi_{1}^{\un}$, the paper \cite{De2} is entitled \emph{Le groupe fondamental de la droite projective moins trois points}, which can be explained as follows (see the introductions of \cite{De2} and of \cite{I}). The origin of the notion of $\pi_{1}^{\un}$ comes from the desire to understand $\Gal(\overline{\mathbb{Q}}/\mathbb{Q})$ via its action on geometric objects ; a theorem of Belyi \cite{Bel} says that the natural map from $\Gal(\overline{\mathbb{Q}}/\mathbb{Q})$ to the group $\text{Out} (\widehat{\pi_{1}}(\mathbb{P}^{1} - \{0,1,\infty\})(\mathbb{C}))$ of outer automorphisms of the pro-finite completion of $\pi_{1}(\mathbb{P}^{1} - \{0,1,\infty\})(\mathbb{C}))$ is injective. Grothendieck has suggested in \cite{Gr2} a way to understand the image of this map by studying more generally the outer action of $\Gal(\overline{\mathbb{Q}}/\mathbb{Q})$ on the $\widehat{\pi_{1}}$ of moduli spaces $\mathcal{M}_{g,n}$ of curves of genus $g$ with $n$ marked points ; and, in particular, $\mathcal{M}_{0,4} \simeq \mathbb{P}^{1} - \{0,1,\infty\}$, $\mathcal{M}_{0,5}$, $\mathcal{M}_{1,1}$ and $\mathcal{M}_{1,2}$. In \cite{De2}, Deligne has proposed to replace in this program the pro-finite completion of the $\pi_{1}$ by its pro-unipotent completion, which is still equipped with an action of $\Gal(\overline{\mathbb{Q}}/\mathbb{Q})$. The advantage of $\pi_{1}^{\un}$ comes from its proximity with cohomology, which follows from anterior works ; whence also a second type of motivation for studying the $\pi_{1}^{\un}$ : to test the theory of motives and periods.

\subsubsection{Multiple zeta values and their expressions as periods, including as Betti-De Rham periods of 
$\pi_{1}^{\un}(\mathbb{P}^{1} - \{0,1,\infty\},-\vec{1}_{1},\vec{1_{0}})$}

Multiple zeta values were discovered, in one and two variables, by Euler, and were forgotten during more than two centuries. Their reapparition in the 1990's is partly due to the role that they play in quantum field theory.

\begin{Definition} \label{Euler Zagier}(Euler-Zagier) Multiple zeta values (MZV's) are the following numbers, for $n_{d},\ldots,n_{1} \in \mathbb{N}^{\ast}$, such that $n_{d} \geq 2$ :
\begin{equation}
\label{eq:1} \zeta(n_{d},\ldots,n_{1}) 
= \sum_{0<m_{1}<\ldots<m_{d}} \frac{1}{m_{1}^{n_{1}} \ldots m_{d}^{n_{d}}} \in \mathbb{R}
\end{equation}
\end{Definition}

\noindent One says that $n=n_{d}+\ldots+n_{1}$ is the weight of $(n_{d},\ldots,n_{1})$, and that $d$ is the depth of $(n_{d},\ldots,n_{1})$. One recognizes the values at positive integers of the Riemann zeta function as the $d=1$ case. It has been noticed by Kontsevich that we have, for all $(n_{d},\ldots,n_{1})$,
\begin{equation} \label{eq:integral} \zeta(n_{d},\ldots,n_{1}) 
= (-1)^{d} \int_{0}^{1} \frac{dt_{n}}{t_{n}-\epsilon_{n}} \int_{0}^{t_{n}} \ldots  \frac{dt_{2}}{t_{2}-\epsilon_{2}} \int_{0}^{t_{2}}\frac{dt_{1}}{t_{1}-\epsilon_{1}} 
\end{equation}
\noindent where $(\epsilon_{n},\ldots,\epsilon_{1}) =(\underbrace{0,\ldots,0}_{n_{d}-1},1,\ldots,\underbrace{0,\ldots,0}_{n_{1}-1},1)$. This shows that multiple zeta values as periods by Definition \ref{second definition of periods}. What one wants, then, is to construct out of this observation a framework which enables to apply the tools of algebraic geometry to study multiple zeta values according to \S1.1. There are two canonical ways to do it ; only the second one will be used in this text.

\begin{itemize}
\item After a sequence of blowing-ups applied to the integral of (\ref{eq:integral}), one obtains an expression of each multiple zeta value $\zeta(n_{d},\ldots,n_{1})$ as a period in the sense of Definition \ref{first definition of periods}, of a relative cohomology group $H^{n}(\overline{\mathcal{M}_{0,n+3}}-A_{(n_{d},\ldots,n_{1})}, B - A_{(n_{d},\ldots,n_{1})} \cap B)$, where $n=n_{d}+\ldots+n_{1}$, $\overline{\mathcal{M}_{0,n+3}}$ is the Deligne-Mumford compactification of $\mathcal{M}_{0,n+3}$ and $A_{(n_{d},\ldots,n_{1})}$ and $B$ are some unions of irreducible components of the normal crossings divisor $\partial\overline{\mathcal{M}_{0,n+3}} = \overline{\mathcal{M}_{0,n+3}} - \mathcal{M}_{0,n+3}$ \cite{GoMan}.
\item An integral such as the one of (\ref{eq:integral}) is called a \emph{iterated path integral}, and is the typical expression of a period of a pro-unipotent fundamental groupoid. The work of Chen on (homotopy-invariant) iterated path integrals \cite{Che} is actually a sort of preliminary Betti - De Rham theory of the $\pi_{1}^{\un}$. Equation (\ref{eq:integral}) shows, precisely, that multiple zeta values are Betti-De Rham periods of $\pi_{1}^{\un}(\mathbb{P}^{1} - \{0,1,\infty\})$ at the couple of tangential base-points $(-\vec{1}_{1},\vec{1_{0}})$ (see \S2.3.3 for details.)
\end{itemize}

\subsubsection{The question of studying multiple zeta values as periods}

The question of understanding the polynomial equations in $\mathbb{Q}$-algebras of periods (\S1.1.2) applies to multiple zeta values. The main objects involved are the motivic Galois action, or certain substitutes to it, and families of polynomial equations which conjecturally generate the whole of the ideal of polynomial equations of multiple zeta values. We will review these objects in \S2.2 - \S2.3. They are all described by fully explicit and simple formulas ; this is a rare and advantageous situation, which enables an explicit theory of multiple zeta values as periods ; in the present work on $p$-adic multiple zeta values, we will extend the area of explicitness of this theory.
\newline The known results of the type "a given polynomial equation does not hold" are the transcendence of $\pi$ (and thus the one of $\zeta(2n) = \frac{|B_{2n}|\pi^{2s}}{2(2n)!}$, for all $n \in \mathbb{N}^{\ast}$), Apéry's result that $\zeta(3) \not\in \mathbb{Q}$, and other results such as Rivoal-Ball's theorem that the values $\zeta(n)$ with $n \in \mathbb{N}^{\ast}$ odd generate an infinite dimensional $\mathbb{Q}$-vector space. It is conjectured that all multiple zeta values are transcendental numbers ; more precisely, the theory started out with the following conjecture :
\begin{Conjecture} \label{conjecture Zagier for multiple zeta values}(Zagier) For each $n \in \mathbb{N}$, let $\mathcal{Z}_{n}$ be the $\mathbb{Q}$-vector space generated by multiple zeta values of weight $n$ (by convention $\mathcal{Z}_{0} =\mathbb{Q}$ is generated by $\zeta(\emptyset) = 1$). Then :
\newline i) The $\mathcal{Z}_{n}$ $(n \in \mathbb{N})$ are in direct sum.
\newline ii) The dimensions of the $\mathcal{Z}_{n}$ are given by the formula : $\sum_{n=0}^{\infty} \dim(\mathcal{Z}_{n})\Lambda^{n} = \frac{1}{1-\Lambda^{2}-\Lambda^{3}}$.
\end{Conjecture}

\noindent Note that we have $\mathcal{Z}_{n}\mathcal{Z}_{n'} \subset \mathcal{Z}_{n+n'}$ for all $n,n'$ (by the two shuffle equations reviewed in \S\ref{double shuffle relations definition})
\footnote{The conjectural homogeneity of algebraic relations with respect to the weight expressed by i) of Conjecture \ref{conjecture Zagier for multiple zeta values} does not remain true if we replace weight by depth : we have $\zeta(3) = \zeta(2,1)$}. The framework of algebraic geometry has enabled to prove that we have $\sum_{n=0}^{\infty} \dim(Z_{n})\Lambda^{n} \leq \frac{1}{1-\Lambda^{2} -\Lambda^{3}}$, \footnote{i.e. for all $n \in \mathbb{N}$, $\dim(\mathcal{Z}_{n}) \leq $ the coefficient of $\Lambda^{n}$ in $ \frac{1}{1-\Lambda^{2} -\Lambda^{3}}$} \cite{T}, \cite{Go1} which is a special case of the majoration reviewed in \S1.1.2. However, the inequality $\geq$ is out of reach at present.

\subsection{Motivations for the explicit theory of $\pi_{1}^{\un,\crys}(\mathbb{P}^{1} - \{0,1,\infty\})$}

The main subject of this text is $\pi_{1}^{\un,\crys}(X_{0})$ with $X_{0} = \mathbb{P}^{1} - \{0,1,\infty\}/\mathbb{F}_{p}$, and more specifically, the families (one for each $\alpha \in (\mathbb{Z} \cup \{\pm \infty\}) - \{0\}$) of numbers $\zeta_{p,\alpha}(n_{d},\ldots,n_{1}) \in \mathbb{Q}_{p}$ called $p$-adic multiple zeta values defined through it.

\subsubsection{Nature of $p$-adic multiple zeta values \label{nature}}

Unlike real multiple zeta values, which express the comparison of two realizations of the $\pi_{1}^{\un}$, $p$-adic multiple zeta values are defined (see \S2.4.2 for details) as the numbers in $\mathbb{Q}_{p}$ expressing canonically the crystalline Frobenius (Definition \ref{def of Frobenius}) iterated $\alpha \in (\mathbb{Z} \cup \{\pm \infty\}) - \{0\}$ times, applied to $X_{0}=\mathbb{P}^{1} - \{0,1,\infty\} / \mathbb{F}_{p}$, at the base-points  $(\vec{1}_{0},-\vec{1}_{1})$. The questions on $p$-adic multiple zeta values are similar to those on multiple zeta values, and we have :

\begin{Conjecture} \label{conjecture p-adic multiple zeta values}(for $\alpha=1$ : Deligne-Goncharov, \cite{DeGo}, \S5.28)
\newline For each $n \in \mathbb{N}$, $p$ prime and $\alpha \in (\mathbb{Z} \cup \{\pm \infty\}) - \{0\}$, let $\mathcal{Z}_{n}^{(p,\alpha)}$ be the $\mathbb{Q}$-vector space generated by numbers $\zeta_{p,\alpha}(w)$ with $w$ of weight $n$ (where $\mathcal{Z}_{0}=\mathbb{Q}$ generated by $\zeta_{p,\alpha}(\emptyset)=1$). 
\newline i) The $\mathcal{Z}_{n}^{(p,\alpha)}$ are in direct sum.
\newline ii) We have $\sum_{n=0}^{\infty} \dim(Z^{(p,\alpha)}_{n})\Lambda^{n} = \frac{1-\Lambda^{2}}{1-\Lambda^{2}-\Lambda^{3}}$.
\newline iii) The ideal of polynomial equations over $\mathbb{Q}$ satisfied by the numbers $\zeta_{p,\alpha}(w)$ is generated by the ideal of polynomial equations over $\mathbb{Q}$ satisfied by multiple zeta values and the relation $\zeta_{p,\alpha}(2)=0$.
\end{Conjecture}

\noindent 
We have, as in the real framework, $Z^{(p,\alpha)}_{n}Z^{(p,\alpha)}_{n'} \subset Z^{(p,\alpha)}_{n+n'}$ for all $n,n'$ (by the two shuffle equations whose $p$-adic analogues are reviewed in \S\ref{p-adic double shuffle relations}).
\newline By (unpublished) work of Yamashta, $p$-adic multiple zeta values in this sense are the images of some $p$-adic periods living in $B = \oplus_{n=0}^{\infty}(B_{crys} \cup \text{Fil}^{0}B_{dR})^{\varphi=p^{n}}$, by the reduction map $B \rightarrow \mathbb{C}_{p}$,
\footnote{Whence the vanishing of the $p$-adic analogue of $\zeta(2)$ : the analogue of $2i\pi$ in $B$ is sent to zero by the reduction map.}and $\sum_{n=0}^{\infty} \dim(Z^{(p,\alpha)}_{n})\Lambda^{n} \leq \frac{1-\Lambda^{2}}{1-\Lambda^{2}-\Lambda^{3}}$ \cite{Yam}. The inequality $\geq$ is out of reach at present like in the real case.
\newline Aside from their intrinsic interest, a motivation for studying $p$-adic multiple zeta values is that this is believed to be necessary in the long-term to understand certain arithmetic aspects of real multiple zeta values. We now explain two more precise and shorter-term motivations ; both of them are related to the fact that, unlike the case of multiple zeta values, there is a priori no straightforward way to compute explicitly $p$-adic multiple zeta values (see \S2.4.3).

\subsubsection{The quasi-shuffle relation and Deligne-Goncharov's question}

One of the standard families of algebraic equations of multiple zeta values, also proved for $p$-adic multiple zeta values, is the so-called quasi-shuffle (sometimes called "stuffle") relation (see \S\ref{double shuffle relations definition}, \S2.4.3), whose example in depth $(1,1)$ is
$$ \zeta(n)\zeta(n') = \zeta(n+n') + \zeta(n,n') + \zeta(n',n) $$
\noindent For real multiple zeta values, the quasi-shuffle relation is a direct consequence of the formula of series of (\ref{eq:1}) combined to a canonical way of writing a product of two sets $\{m_{1}<\ldots<m_{d}\} \subset (\mathbb{N}^{\ast})^{d}$,  $\{m'_{1}<\ldots<m'_{d'}\} \subset (\mathbb{N}^{\ast})^{d'}$ as a disjoint union of sets of the same type in $(\mathbb{N}^{\ast})^{r}$, $r \in \{\max(d,d'),\ldots,d+d'\}$ : in depth $(1,1)$ see equation (\ref{eq:proof of quasi shuffle}) below. The fact that the quasi-shuffle relation holds for $p$-adic multiple zeta values can be proved without explicit formulas (see \S2.4.3) but suggests the following question :

\begin{Question} (Deligne-Goncharov, \cite{DeGo}, \S5.28)
\newline "\emph{Il serait aussi intéressant de disposer pour ces coefficients d'expressions p-adiques qui rendent clair qu'ils vérifient des identités du type [...]
\begin{equation} \label{eq:proof of quasi shuffle} \sum \frac{1}{k^{n}} \sum \frac{1}{l^{m}} = \sum_{k>l} \frac{1}{k^{n}} \frac{1}{l^{m}} + \sum_{l>k} \frac{1}{k^{n}}\frac{1}{l^{m}} + \sum_{l=k} \frac{1}{k^{n+m}}" 
\end{equation}}
\end{Question}

\noindent i.e. they ask to compute explicitly $p$-adic multiple zeta values in a way which enables a $p$-adic analogue of the above proof that multiple zeta values satisfy the quasi-shuffle relation.

\subsubsection{Finite multiple zeta values and Kaneko-Zagier's conjecture \label{paragraph def of finite}}

Let $\mathcal{P}$ be the set of prime numbers, and let $\mathbb{Z}/_{p \rightarrow \infty}$ be the ring of "integers modulo infinitely large primes" \cite{Ko}, which is of characteristic zero \footnote{Kaneko and Zagier as well as people working on finite multiple zeta values denote this ring by $\mathcal{A}$ ; our notation $\mathbb{Z}/_{p \rightarrow \infty}$ is not standard} :
$$ \mathbb{Z}/_{p \rightarrow \infty} = \big( \prod_{p \in \mathcal{P}} \mathbb{Z}/p\mathbb{Z} \big) / \big( \bigoplus_{p \in \mathcal{P}} \mathbb{Z}/p\mathbb{Z} \big) = \big( \prod_{p \in \mathcal{P}} \mathbb{Z}/p\mathbb{Z} \big) / \big( \prod_{p \in \mathcal{P}} \mathbb{Z}/p\mathbb{Z} \big)_{\text{tors}} \simeq \big( \prod_{p \in \mathcal{P}} \mathbb{Z}/p\mathbb{Z} \big) \otimes_{\mathbb{Z}} \mathbb{Q} $$
\noindent where the subscript tors refers to the torsion subgroup. We have an embedding $\mathbb{Q} \hookrightarrow \mathbb{Z}/_{p \rightarrow \infty}$.

\begin{Definition} \label{definition finite multiple zeta values}(Zagier, unpublished) Finite multiple zeta values are the following numbers, for $(n_{d},\ldots,n_{1}) \in (\mathbb{N}^{\ast})^{d}$, $d \in \mathbb{N}^{\ast}$ \footnote{The usual notation is $\zeta_{\mathcal{A}}(n_{d},\ldots,n_{1})$}:
\begin{equation} \zeta_{\mathbb{Z}/_{p \rightarrow \infty}}(n_{d},\ldots,n_{1}) =  \text{ the image of } \left(\sum_{0<m_{1}<\ldots<m_{d}<p} \frac{1}{m_{1}^{n_{1}} \ldots m_{d}^{n_{d}}} \mod p \right)_{p \in \mathcal{P}} \text{ in }\mathbb{Z}/_{p \rightarrow \infty}
\end{equation}
\end{Definition}

\noindent This notion is close to other notions and ideas that appeared earlier in the work of several authors, in particular Hoffman \cite{Ho} and Zhao \cite{Zh}. The explanation for this terminology is the following :

\begin{Conjecture} \label{conjecture Kaneko-Zagier}(Kaneko-Zagier) The following correspondence defines an isomorphism of $\mathbb{Q}$-algebras from the algebra generated by finite multiple zeta values to the algebra generated by multiple zeta values modulo $(\zeta(2))$
\footnote{By Definition \ref{Euler Zagier},  multiple zeta values $\zeta(n_{d},\ldots,n_{1})$ are defined when $n_{d} \geq 2$ whereas by Definition \ref{definition finite multiple zeta values} finite multiple zeta values $\zeta_{\mathbb{Z}/_{p \rightarrow \infty}}(n_{d},\ldots,n_{1})$ are defined for all $n_{d} \geq 1$. When $n_{d}=1$, the map (\ref{eq:32})  refers to the two notions of "regularized" $\zeta(n_{d},\ldots,n_{1})$ which are defined respectively by a regularization of the right-hand side of (\ref{eq:1}) and (\ref{eq:integral}) (see \cite{Ra} for a review of the regularizations) ; the map is well-defined because the two regularizations of $\sum_{d'=0}^{d}(-1)^{n_{d'+1}+\ldots+n_{d}}\zeta(n_{d'+1},\ldots,n_{d})\zeta(n_{d'},\ldots,n_{1})$ have the same reduction modulo $\zeta(2)$.} :
\begin{equation} \label{eq:32} \zeta_{\mathbb{Z}/_{p \rightarrow \infty}}(n_{d},\ldots,n_{1}) \mapsto \sum_{d'=0}^{d}(-1)^{n_{d'+1}+\ldots+n_{d}}\zeta(n_{d'+1},\ldots,n_{d})\zeta(n_{d'},\ldots,n_{1}) \mod \zeta(2)
\end{equation}
\end{Conjecture}
\noindent This conjecture is striking, both because it involves the ring $\mathbb{Z}/_{p \rightarrow \infty}$ and because of the explicit formula. Comparing it to Conjecture \ref{conjecture p-adic multiple zeta values} and writing
$$ \mathbb{Z}/_{p \rightarrow \infty} = \bigg\{ (x_{p})_{p} \in \prod_{p\in \mathcal{P}}\mathbb{Q}_{p} \text{ }|\text{ } v_{p}(x_{p}) \geq 0 \text{ for p large} \bigg\} \text{ }\bigg/\text{ }\bigg\{ (x_{p})_{p} \in \prod_{p\in \mathcal{P}} \mathbb{Q}_{p} \text{ }|\text{ } v_{p}(x_{p}) \geq 1 \text{ for p large} \bigg\} $$
\noindent enables to imagine a relation between finite multiple zeta values and $p$-adic multiple zeta values, which could be an accessible part of the conjecture and partially explain it. This suggests, by extrapolating a little more, the existence of readable explicit formulas for $p$-adic multiple zeta values.

\subsubsection{The initial technical question \label{initial question}}

By \S1.4.1, \S1.4.2 and \S1.4.3, we are motivated to solve the problem of computing explicitly $p$-adic multiple zeta values, and we have in mind that a test of the success of the computation would be to see whether it sheds light in some way on the questions of \S1.4.2 and \S1.4.3.
\newline It is granted at least that there is a certain way to compute $p$-adic multiple zeta values, by the overconvergence of the differential equation satisfied by the Frobenius (\S2.4.3). To our knowledge, not much is granted regarding the form of the answer. We will see that there exists a primary form of answer which looks complicated, but that, by introducing certain particular principles of computation (\S3) combined to certain ideas of $p$-adic analysis (\S4.2), we can obtain different, simpler and more exploitable formulas ; this will also lead us to enlarge our set of questions.

\section{Review on $\pi_{1}^{\un}(\mathbb{P}^{1} - \{0,1,\infty\}$)}

We review some elements of explicit description of $\pi_{1}^{\un}(\mathbb{P}^{1} - \{0,1,\infty\})$, which make concrete the question of computing $p$-adic multiple zeta values. The pro-unipotent harmonic actions, which we will define in \S4.2-\S4.3, are connected to the Poisson-Ihara bracket which we review in \S2.2.

\subsection{Explicit description of $\pi_{1}^{\un,\DR}(\mathbb{P}^{1} - \{0,1,\infty\})$}

In this paragraph, $X$ is $\mathbb{P}^{1} - \{0,1,\infty\}$ over any field $K$ of characteristic zero. The canonical base-point $\omega_{\DR}$ is defined by \cite{De2}, \S12.4 ; the tangential base-points are defined by \cite{De2}, \S15 and for $u \in \{0,1,\infty\}$, $T_{u}$ means the tangent space of $\mathbb{P}^{1}$ at $u$.

\subsubsection{The groupoid $\pi_{1}^{\un,\DR}(\mathbb{P}^{1} - \{0,1,\infty\})$ and its canonical base-point $\omega_{\DR}$}

For each element $x$ of $\mathbb{P}^{1} - \{0,1,\infty\}(K) \bigcup \big( \cup_{u \in \{0,1,\infty\}} (T_{u} - \{0\})(K)\big)$, one has a tensor functor "fiber at $x$" : $C^{\un,\DR} \rightarrow \Vect(K)$, thus a base-point of the fundamental groupoid, which we will denote also by $x$.
\newline For each couple $(x,y)$ of base-points, one has a pro-affine scheme $\pi_{1}^{\un,\DR}(X,y,x)$ over $\mathbb{Z}$ : it is, by definition, the scheme of tensor automorphisms between the fiber functors $x$ and $y$. The points of $\pi_{1}^{\un,\DR}(X,y,x)$ are called De Rham paths from $x$ to $y$.
\newline For each triple $(x,y,z)$ of base-points one has a morphism of schemes 
$\pi_{1}^{\un,\DR}(X,z,y) \times  \pi_{1}^{\un,\DR}(X,y,x) \rightarrow \pi_{1}^{\un,\DR}(X,z,x)$ called the groupoid multiplication.
\newline When $x=y$, $\pi_{1}^{\un,\DR}(X,x)=\pi_{1}^{\un,\DR}(X,x,x)$ is a group scheme. The groupoid multiplication makes each $\pi_{1}^{\un,\DR}(X,y,x)$ into a bi-torsor under $(\pi_{1}^{\un,\DR}(X,x),\pi_{1}^{\un,\DR}(X,y))$.
\newline 
\newline The functor "global section" : $\omega_{\DR} : C^{\un,\DR} \rightarrow \Vect(k)$ is a tensor functor (this is because  $H^{1}(\mathbb{P}^{1},\mathcal{O}_{\mathbb{P}^{1}}) = 0$), and for each base-point $x$, we have a canonical isomorphism $x \simeq \omega_{\DR}$. As a consequence, we have a pro-affine scheme $\pi_{1}^{\un,\DR}(X,\omega_{\DR})$, with, for each $(x,y)$, a canonical isomorphism $\pi_{1}^{\un,\DR}(X,\omega_{\DR}) \simeq \pi_{1}^{\un,\DR}(X,y,x)$ and a canonical path ${}_y 1_{x}$ of $\pi_{1}^{\un,\DR}(X,y,x)$. They are compatible with the groupoid structure : $({}_z 1_{y})({}_y 1_{x}) = {}_z 1_{x}$.
\newline Describing explicitly the groupoid $\pi_{1}^{\un,\DR}(\mathbb{P}^{1} - \{0,1,\infty\})$ is thus reduced to describing explicitly the group scheme $\pi_{1}^{\un,\DR}(\mathbb{P}^{1} - \{0,1,\infty\},\omega_{\DR})$.

\subsubsection{Explicit description of $\pi_{1}^{\un}(\mathbb{P}^{1} - \{0,1,\infty\},\omega_{\DR})$}

The shuffle Hopf algebra is a functor from the category of (finite) sets to the category of graded Hopf algebras over $\mathbb{Q}$. We need only the following example.

\begin{Proposition} The $\mathbb{Q}$-vector space $\mathcal{O}^{\sh}=\mathbb{Q}\langle e_{0},e_{1}\rangle$, freely generated by words on $\{e_{0},e_{1}\}$ including the empty word, endowed with the following operations $\sh,\Delta_{\dec},\epsilon,S$, is a Hopf algebra over $\mathbb{Q}$, graded by the number of letter of words :
\begin{itemize}
\item the shuffle product $\sh:\mathcal{O}^{\sh} \otimes \mathcal{O}^{\sh} \rightarrow \mathcal{O}^{\sh}$ defined by $e_{i_{n+n'}}\ldots e_{i_{n+1}}\text{ }\sh\text{ }e_{i_{n}} \ldots e_{i_{1}} =$
\newline $\sum_{\sigma}
e_{i_{\sigma^{-1}(n+n')}} \ldots e_{i_{\sigma^{-1}(1)}}$, 
where the sum is over permutations $\sigma$ of $\{1,\ldots,n+n'\}$ such that $\sigma(l)<\ldots<\sigma(1)$  and $\sigma(n+n')<\ldots<\sigma(l+1)$
\item the deconcatenation coproduct $\Delta_{\dec}:\mathcal{O}^{\sh} \rightarrow \mathcal{O}^{\sh} \otimes \mathcal{O}^{\sh}$, defined by  $\Delta_{\dec}(e_{i_{n}}\ldots e_{i_{1}}) = \sum_{n''=0}^{n} e_{i_{n}}\ldots e_{i_{n''+1}} \otimes e_{i_{n''}} \ldots e_{i_{1}}$
\item the counit $\epsilon : \mathcal{O}^{\sh} \rightarrow \mathbb{Q}$ equal to the augmentation map of $\mathcal{O}^{\sh}$, 
\item the antipode $S : \mathcal{O}^{\sh} \rightarrow \mathcal{O}^{\sh}$, defined by $S(e_{i_{n}}\ldots e_{i_{1}}) = (-1)^{n} e_{i_{1}}\ldots e_{i_{n}}$.
\end{itemize}
\end{Proposition}

\begin{Definition} $\mathcal{O}^{\sh}$ is called the shuffle Hopf algebra over the alphabet $\{e_{0},e_{1}\}$, and the number of letters of a word is called its weight.
\end{Definition}

\begin{Proposition} \label{shuffle group scheme}(follows from \cite{De2}, \S12, Proposition 12.7, Corollaire 12.9) The group scheme $\Spec(\mathcal{O}^{\sh})$ is pro-unipotent and is canonically isomorphic to $\pi_{1}^{\un,\DR}(X,\omega_{\DR})$.
\end{Proposition}

\noindent We will denote in the same way $e_{0},e_{1} \in \mathcal{O}^{\sh}$ and their duals below. Let 
$K\langle \langle e_{0},e_{1} \rangle\rangle$ be the non-commutative $K$-algebra of formal power series with variables $e_{0},e_{1}$ with coefficients in $K$. An element $f$ of $K\langle \langle e_{0},e_{1} \rangle\rangle$ can be written in a unique way as $f = \sum_{w\text{ word over } \{e_{0},e_{1}\}} f[w] w$ with $f[w] \in K$ for all $w$. The notation $f[w]$ extends to linear combinations of words by linearity. 

\begin{Proposition} \label{points of spec o sh} 
The dual of the topological Hopf algebra $\mathcal{O}^{\sh} \otimes_{\mathbb{Q}} K$ is $K\langle\langle e_{0},e_{1} \rangle\rangle$ viewed as the topological Hopf algebra associated with the universal enveloping algebra of the complete free Lie algebra over the two variables $e_{0},e_{1}$.
\newline The group $\Spec(\mathcal{O}^{\sh})(K)$ is the group of the elements of $K\langle\langle e_{0},e_{1} \rangle\rangle$ which satisfy the following property called the shuffle equation : for all words $w,w'$, $f[w]f[w'] = f[w\text{ }\sh\text{ }w']$.
\end{Proposition}

\subsubsection{The fundamental torsor of paths and its connection}

The groupoid $\pi_{1}^{\un,\DR}(\mathbb{P}^{1} - \{0,1,\infty\})$ is an initial object of the category of groupoids with integrable connection (by \cite{De2}, \S10.49), in the sense of \cite{De2}, \S10.28. It follows from \cite{De2}, \S12 that :

\begin{Proposition} The connection on $\pi_{1}^{\un,\DR}(\mathbb{P}^{1} - \{0,1,\infty\})$ is, in the sense of \S7.30.2$_{S}$, the map
$$ \nabla_{\KZ} : f \mapsto f^{-1} \big( df - (e_{0} \frac{dz}{z} + e_{1} \frac{dz}{z-1})f \big) $$
\end{Proposition}

\begin{Definition} $\nabla_{\KZ}$ is called the Khnizhnik-Zamolodchikov connection.
\end{Definition}

\subsection{The motivic Galois action on $\pi_{1}^{\un,\DR}(\mathbb{P}^{1} - \{0,1,\infty\},-\vec{1}_{1},\vec{1}_{0})$ \label{paragraph Ihara}}

Following \cite{DeGo}, we now review how the Tannakian framework of the motivic Galois theory of periods evoked in \S1.1.2 applies to $\pi_{1}^{\un}(\mathbb{P}^{1} - \{0,1,\infty\})$.

\subsubsection{Mixed Tate motives and $\pi_{1}^{\un,\mot}(\mathbb{P}^{1} - \{0,1,\infty\})$ \label{deux deux un}}

Let $k$ be a number field. By Levine \cite{L}, one has the Tannakian category $\MT(k)$ of mixed Tate motives over $k$ and, for $S$ a set of finite places $k$, one has the Tannakian category $\MT(\mathcal{O}_{S})$ of mixed Tate motives over $\mathcal{O}_{S}$, which is a subcategory of $\MT(k)$. They are described in \cite{DeGo}, \S1,\S2.
\newline The De Rham realization functor $\DR : \MT(k) \rightarrow \Vect_{k}$ is canonically isomorphic to the extension of scalars from $\mathbb{Q}$ to $k$ of the "canonical realization functor" $\omega$ (\cite{DeGo}, Proposition \S2.10). The motivic Galois group $G^{\omega} = \Aut^{\otimes}(\omega)$ associated with $\omega$ on $\MT(\mathcal{O}_{S})$ with $S$ as above has a semi-direct product decomposition $G^{\omega} = \mathbb{G}_{m} \ltimes U^{\omega}$ where $U^{\omega}$ is a pro-unipotent group scheme, described explicitly in \cite{DeGo}, Proposition 2.2.
\newline The groupoid $\pi_{1}^{\un,\mot}(\mathbb{P}^{1} - \{0,1,\infty\})$ is defined as a groupoid over $\mathbb{P}^{1} - \{0,1,\infty\}$ in affine schemes in $\MT(\mathbb{Q})$ (\cite{DeGo}, \S3.12, \S3.13, \S4) ; the notion of affine schemes in a Tannakian category is reviewed in \cite{DeGo}, \S2.6.

\subsubsection{Description of the action of $G^{\omega}$ on $\pi_{1}^{\un,\DR}(\mathbb{P}^{1} - \{0,1,\infty\},-\vec{1}_{1},\vec{1}_{0})$ ; the Goncharov coproduct and the Ihara bracket
\label{formula for motivic Galois action}}

The motivic Galois group $G^{\omega}$ associated with $\MT(\mathbb{Z})$ acts on $\pi_{1}^{\un,\DR}(\mathbb{P}^{1} - \{0,1,\infty\},-\vec{1}_{1},\vec{1}_{0})$. We review its description.
\newline 
\newline $\bullet$ The action of $\mathbb{G}_{m} \subset G^{\omega}$ expresses the weight grading : it provides the motivic way to formulate the conjecture that all polynomial relations between multiple zeta values are homogeneous for the weight (Conjecture \ref{conjecture Zagier for multiple zeta values}, i)) :

\begin{Definition}
	Let	$\tau$ be the action of $\mathbb{G}_{m}$ on $\pi_{1}^{\un,\DR}(\mathbb{P}^{1} - \{0,1,\infty\},-\vec{1}_{1},\vec{1}_{0})$ which sends $\displaystyle(\lambda, f = \sum_{w\text{ word}} f[w]w)$ to $\displaystyle\sum_{w\text{ word}} \lambda^{\weight(w)}f[w]w$.
\end{Definition}

\noindent $\bullet$ Since the functors of Hodge realization of categories of mixed Tate motives are fully faithful (\cite{DeGo}, Proposition 2.14), in order to compute the action of $U^{\omega}$, it is sufficient to compute its Hodge realization. This has been done by Goncharov : \cite{Go1}, Theorem 6.4, Theorem 6.5. The formula is usually written in terms Hopf algebras, is often called the Goncharov coaction or Goncharov coproduct. Here, what we are interested in is its De Rham realization. 
We will denote, as in \cite{DeGo}, \S5, by 
$\Pi_{1,0} = \pi_{1}^{\un,\DR}(\mathbb{P}^{1} - \{0,1,\infty\},-\vec{1}_{1},\vec{1}_{0})$, and by $\Pi_{0,0}$, $\Pi_{0,1}$, $\Pi_{1,0}$ the other similar schemes, where we associate $1$ with $-\vec{1}_{1}$ and $0$ with $-\vec{1}_{0}$.

\begin{Definition} (\cite{DeGo}, \S5.6 - \S5.9) 
	\label{def de V omega}Let $V_{\omega}$ be the group of automorphisms of the four schemes  $\Pi_{0,0},\Pi_{1,0},\Pi_{0,1},\Pi_{1,1}$ which preserve the groupoid multiplication and fix $\exp(e_{0})$ at $(0,0)$ and $\exp(e_{1})$ at $(1,1)$.
\end{Definition}

\begin{Theorem} (\cite{DeGo}, equation (5.10.3)) There exists a morphism $\mathbb{G}_{m} \ltimes U_{\omega} \rightarrow \mathbb{G}_{m} \ltimes V_{\omega}$, compatible to the semi-direct product decompositions and the actions on $\Pi_{1,0}$.
\end{Theorem}

\begin{Theorem} \label{formula for Ihara action}
	(\cite{DeGo}, Proposition 5.11) \label{prop Ihara action}
\newline i) $\Pi_{1,0}$ is a $V_{\omega}$-torsor, and we thus have an isomorphism of schemes $V_{\omega} \simeq \Pi_{1,0}$, $v \mapsto v({}_1 1_{0})$ (where ${}_1 1_{0} = {}_{-\vec{1}_{1}} 1_{\vec{1}_{0}}$ is the canonical path as in \S2.1.1).
\newline ii) More precisely, for $v$ a point of $V_{\omega}$ and $g = v({}_1 1_{0})$, the action of $v$ on $\Pi_{0,0}$ is $g(e_{0},e_{1}) \mapsto (f(e_{0},e_{1}) \mapsto f(e_{0},g^{-1}e_{1}g))$ and the action of $v$ on $\Pi_{1,0}$ is $g(e_{0},e_{1}) \mapsto (f(e_{0},e_{1}) \mapsto g(e_{0},e_{1}).f(e_{0},g^{-1}e_{1}g))$.
\newline iii) The map $(g(e_{0},e_{1}),f(e_{0},e_{1}) \mapsto g(e_{0},e_{1}).f(e_{0},g^{-1}e_{1}g)$ defines a group law on $\Pi_{1,0}$.
\end{Theorem}

\begin{Notation} We denote the group law above by $\circ^{\smallint_{0}^{1}}$. (This notation is not standard.)
\end{Notation}

\noindent Some authors call Ihara action or Ihara product the above action of $V_{\omega}$ on $\Pi_{1,0}$, since it has been first written explicitly by Ihara \cite{I}. The Lie bracket of $\Lie V_{\omega}$ is sent by the isomorphism $V_{\omega} \simeq \Pi_{1,0}$ to a Lie bracket on $\Lie \Pi_{1,0}$ called the Ihara bracket, or the Poisson-Ihara bracket, or the Poisson bracket on the free Lie algebra in two generators.  Some other authors call ($\Pi_{1,0},\circ^{\smallint_{0}^{1}})$ the twisted Magnus group. 

\subsection{Betti-De Rham comparison of $\pi_{1}^{\un}(\mathbb{P}^{1} - \{0,1,\infty\})$ and multiple zeta values}

\subsubsection{Iterated path integrals}

\begin{Definition} (Chen, \cite{Che})
Let $\eta_{1},\ldots,\eta_{r}$ be differential $1$-forms on a manifold $M$, and let $\gamma : [0,1] \rightarrow M$ be a smooth path. Denote by $f_{i}(t)dt = \gamma^{\ast}(\eta_{i})$, $i=1,\ldots,r$. The iterated path integral of $(\eta_{1},\ldots,\eta_{r})$ along $\gamma$ is 
$$ \int_{\gamma} \eta_{r} \ldots \eta_{1} = \int_{0}^{1} f_{r}(t)dt \int_{0}^{t_{r}} \ldots \int_{0}^{t_{2}} f_{1}(t)dt $$
\end{Definition}

\noindent If the differential forms have logarithmic singularities at infinity, the definition can be extended to the case where the extremities of $\gamma$ are tangential-base points \footnote{For example in the case of $\mathbb{P}^{1} - \{0,1,\infty\}$ : let $\gamma$ be a smooth path $[0,1] \rightarrow \mathbb{P}^{1}(\mathbb{C})$ such that for $t \in ]0,1[$, we have $\gamma(t) \in (\mathbb{P}^{1} - \{0,1,\infty\})(\mathbb{C})$, and $\gamma(0) \in \{0,1,\infty\}$, then $\gamma$ is said to have extremity $\overrightarrow{\gamma'(0)}_{\gamma(0)}$ ; analogous definition at $t=1$. The notions of homotopy of paths and of topological fundamental groupoid extend to such paths, as does the Malcev completion defining $\pi_{1}^{\un,\B}(\mathbb{P}^{1} - \{0,1,\infty\})$ from Definition \ref{Betti definition}.} ; this sometimes requires to regularize a divergent iterated integral
\footnote{For example in the case of $\mathbb{P}^{1} - \{0,1,\infty\}$ : if $\gamma(0) \in \{0,1,\infty\}$, we consider the analogue iterated integral on $\gamma([\epsilon,1])$ with $0<\epsilon<1$, which has an asymptotic expansion in $\mathbb{C}[[\epsilon]][\log(\epsilon)]$ when $\epsilon \rightarrow 0$ and define the regularized iterated integral on $\gamma$ as the constant coefficient of its asymptotic expansion}.
\newline In all the rest of this text, we will choose paths in $\mathbb{P}^{1} - \{0,1,\infty\}$ starting at the tangential base-point $\vec{1}_{0}$.

\subsubsection{Multiple polylogarithms and their series expansion at $0$}

\begin{Definition} (Goncharov, \cite{Go1}) \label{def des multiple polylog}Multiple polylogarithms (MPL's) on $\mathbb{P}^{1} - \{0,1,\infty\}$ are multivalued holomorphic functions on $\mathbb{P}^{1} - \{0,1,\infty\}(\mathbb{C})$ solutions to $\nabla_{\KZ}$ of \S2.1.3 : for each $\gamma$ a path in $\mathbb{P}^{1} - \{0,1,\infty\}$ starting at $\vec{1}_{0}$, let
$$ \Li(\gamma) = 1+ \sum_{\substack{n \in \mathbb{N}^{\ast} \\  i_{n},\ldots,i_{1} \in \{0,1\}^{n}}} \bigg( \int_{\gamma} \frac{dz}{z- i_{n}} \ldots \frac{dz}{z- i_{1}} \bigg) e_{i_{n}}\ldots e_{i_{1}} \in \pi_{1}^{\un,\DR}(\mathbb{P}^{1} - \{0,1,\infty\},z,\vec{1}_{0}) $$
\noindent Indeed, $\Li(\gamma)$ depends only on the homotopy class of $\gamma$ ; it defines a function on the topological $\pi_{1}$, and then on the $\pi_{1}^{\un,\B}$ via the Malcev completion map.
\end{Definition}

\subsubsection{Series expansion of multiple polylogarithms ; multiple harmonic sums \label{deux point trois point deux}}

\begin{Proposition} (Goncharov, \cite{Go1})\label{mhs and mpl} For $d \in \mathbb{N}^{\ast}$, and $n_{1},\ldots,n_{d} \in \mathbb{N}^{\ast}$, let $w= e_{0}^{n_{d}-1}e_{1} \ldots e_{0}^{n_{1}-1}e_{1}$. For $z \in \mathbb{C}$ such that $|z|<1$, we have, denoting by $\Li[w](z)$ the value of $\Li[w]$ at the straight path from $0$ to $z$,
\footnote{Note that $\Li[e_{1}](z) = - \log(1-z)$. We take the convention that the values $\Li[w]$ for $w$ a word of the form $w'e_{0}$ are determined by  $\Li[e_{0}](z)=\log(z)$, the values of Proposition 2.14 and the shuffle equation of Proposition 2.4 : namely, using that for all words $w$, and $x \in \{0,1\}$, we have $\Li[e_{1}e_{x}w](z) = -\Li[e_{x}(e_{1}\text{ }\sh\text{ }w)](z)$ and $\Li[we_{x}e_{0}](z) = -\Li[(w\text{ }\sh\text{ }e_{0})e_{x}](z)$.}
\begin{equation}
\label{eq: series expansion polylog} \Li[w](z) = \sum_{0<m_{1}<\ldots<m_{d}} \frac{z^{m_{d}}}{m_{1}^{n_{1}} \ldots m_{d}^{n_{d}}}
\end{equation}
\end{Proposition}

\noindent We retrieve the formula (\ref{eq:1}) for multiple zeta values by taking the limit $z \rightarrow 1$ : see  \S2.3.3 below. We are led to consider intrinsically the coefficients of this series expansion :

\begin{Definition} \label{definition weighted multiple harmonic sums}i) Let weighted multiple harmonic sums be the following numbers, where $d \in \mathbb{N}^{\ast}$, and $n_{1},\ldots,n_{d} \in \mathbb{N}^{\ast}$,
$$ \har_{m}(n_{d},\ldots,n_{1}) = m^{n_{d}+\ldots+n_{1}} \sum_{0<m_{1}<\ldots<m_{d}<m} \frac{1}{m_{1}^{n_{1}} \ldots m_{d}^{n_{d}}} $$
\noindent ii) Let prime weighted multiple harmonic sums be the weighted multiple harmonic sums of the form $\har_{p^{\alpha}}(n_{d},\ldots,n_{1})$ with $p$ a prime number and $\alpha \in \mathbb{N}^{\ast}$.
\end{Definition}

\noindent The adjective weighted \footnote{the use of the word weighted for this purpose appears in \cite{Ro} in a particular case} refers to the factor $m^{n_{d}+\ldots+n_{1}}$. By (\ref{eq: series expansion polylog}) we have, for all $l \in \mathbb{N}^{\ast}$,
\begin{equation} \label{eq:M mhs and mpl} \har_{m}(w) = \tau(n)\Li[e_{0}^{l}e_{1} w][z^{m}] = \tau(m)\sum_{m'=1}^{m-1}\Li[w][z^{m'}]
\end{equation}
\noindent where $[z^{m'}]$ means the coefficient of degree $m'$ in the power series expansion.

\subsubsection{Expression of multiple zeta values as Betti-De Rham periods of $\pi_{1}^{\un}(\mathbb{P}^{1} - \{0,1,\infty\},-\vec{1}_{1},\vec{1}_{0})$}

\begin{Definition} (\cite{Dr}, \S2, \cite{DeGo}, \S5.18) 
\newline Let $\dch \in \pi_{1}^{\un,\B}(\mathbb{P}^{1} - \{0,1,\infty\},-\vec{1}_{1},\vec{1}_{0})(\mathbb{C})$ be the image, by the completion map of Definition \ref{Betti definition}, of the homotopy class of the straight path $\gamma : [0,1] \rightarrow [0,1]$. Let $\Phi_{\KZ}$ be the image of $\dch$ by the Betti-De Rham comparison isomorphism of Theorem \ref{Betti De Rham comparison}. $\Phi_{\KZ}$ is called the KZ Drinfeld associator (the origin of this terminology is in \cite{Dr}).
\end{Definition}

\begin{Proposition} Equation (\ref{eq:integral}) is equivalent to : for all $(n_{d},\ldots,n_{1})$ such that $n_{d} \geq 2$, \footnote{The other coefficients of $\Phi_{\KZ}$ are some $\mathbb{Q}$-linear combinations of multiple zeta values, as it can be seen either by computing regularized iterated integrals, or by the shuffle equation of Proposition 2.4 applied to $\Phi_{\KZ}$. They are called regularized multiple zeta values. We have $\Phi_{\KZ}[e_{0}] = \Phi_{\KZ}[e_{1}]=0$.}
$$ \zeta(n_{d},\ldots,n_{1}) = (-1)^{d} \Phi_{\KZ}[e_{0}^{n_{d}-1}e_{1}\ldots e_{0}^{n_{1}-1}e_{1}] $$
\end{Proposition}

\subsubsection{Conjecture of periods for multiple zeta values}

The conjecture of periods in the case of multiple zeta values amounts to say that all polynomial equations of multiple zeta values are preserved by the motivic Galois action evoked in \S\ref{formula for motivic Galois action}, i.e. are "motivic relations".

\subsubsection{Explicit algebraic relations of multiple zeta values and product $\circ^{\smallint_{0}^{1}}$ ; focus on double shuffle relations \label{double shuffle relations definition}}

One has three standard families of explicit algebraic relations over $\mathbb{Q}$ between multiple zeta values, each of them is conjectured to generate all their algebraic relations : they are called respectively the associator relations, the regularized double shuffle relations, and the Kashiwara-Vergne relations. (See \cite{F3} for a review of these three notions and \cite{Dr} for the definition of associators). They are known to be motivic (see \cite{So} \footnote{Although this paper concerns the other definition of motivic multiple zeta values arising from \cite{GoMan} evoked in \S1.3.2, which is conjecturally equivalent to the notion of motivic multiple zeta values arising from $\pi_{1}^{\un}(\mathbb{P}^{1} - \{0,1,\infty\})$}, \cite{F3}). There are several known inclusions between the ideals respectively generated by these relations \cite{AET}, \cite{AT}, \cite{F4}, \cite{Sch}.
\newline In this text we will focus on the regularized double shuffle relations, namely
\begin{itemize} \item two notions of regularized multiple zeta values, extending the Definition \ref{Euler Zagier} of multiple zeta values : $\zeta_{\sh}(w)=\Phi_{\KZ}[w]$ for all words $w$ on $e_{0},e_{1}$, and $\zeta_{\ast}(n_{d},\ldots,n_{1})$ for all words $(n_{d},\ldots,n_{1})$ including $n_{d}=1$, and  a formula relating $\zeta_{\sh}$ and $\zeta_{\ast}$ (see \cite{Ra} for details 
	\footnote{$\zeta_{\sh}$ is defined by replacing an iterated integral from $0$ to $1$ by an iterated integral from $\epsilon$ to $1-\epsilon'$, writing an asymptotic expansion when $\epsilon,\epsilon'\rightarrow 0$ (this depends on a choice of tangential base-points and one takes $\vec{1}_{0}$ and $-\vec{1}_{1}$) and considering its constant term ; similarly, $\zeta_{\ast}$ is obtained by replacing the infinite sum over $0<n_{1}<\ldots<n_{d}$ by a truncated sum over $0<n_{1}<\ldots<n_{d}<n$, writing an asymptotic expansion when $n \rightarrow \infty$ and considering its constant term, in which, in addition, we suppress the terms where the Euler-Mascheroni constant appears ; this gives for example that $\zeta_{\ast}(1)$ is $0$. The notions $\zeta$, $\zeta_{\sh}$, $\zeta_{\ast}$ agree on words $(n_{d},\ldots,n_{1})$ with $n_{d} \geq 2$})
\item the shuffle relation, satisfied by all points of $\Spec(\mathcal{O}^{\sh})$ (Proposition \ref{points of spec o sh}) : for all words $w,w'$, $\zeta_{\sh}(w)\zeta_{\sh}(w')=\zeta_{\sh}(w\text{ }\sh\text{ }w')$. This also follows from equation (\ref{eq:integral}) combined to the canonical way to write a product
$\{0<t_{1}<\ldots<t_{n}<1\} \times \{0<t'_{1}<\ldots<t'_{n'}<1\} \subset \mathbb{R}^{n+n'}$ as a disjoint union of simplices of $\mathbb{R}^{n+n'}$ up to sets of measure $0$. (first example : $\int_{0<t_{1}<1} \times \int_{0<t'_{1}<1} = \int_{0<t_{1}<t'_{1}<1} + \int_{0<t'_{1}<t_{1}<1}$)
\item the quasi-shuffle (or stuffle) relation, already described in \S1.4.3 : it is expressed by the following formula : for all words $w,w'$ whose furthest to the right letter is $e_{1}$ $\zeta_{\ast}(w)\zeta_{\ast}(w') = \zeta_{\ast}(w \ast w')$ where $\ast$ is a bilinear map $\mathcal{O}^{\sh} \times \mathcal{O}^{\sh} \rightarrow \mathcal{O}^{\sh}$ which lifts the canonical expression of a product $\{0<m_{1}<\ldots<m_{d}\} \times \{0<m'_{1}<\ldots<m'_{d'}\} \subset \mathbb{N}^{d} \times \mathbb{N}^{d'}$ as a disjoint union of "simplices" of $\mathbb{N}^{r}$, $r \in \{\max(d,d'),\ldots,d+d'\}$ (first example : $\sum_{0<m_{1}} \times \sum_{0<m'_{1}} = \sum_{0<m_{1}=m'_{1}} + \sum_{0<m_{1}<m'_{1}} + \sum_{0<m'_{1}<m_{1}}$)
\end{itemize}

\noindent For the three standard families of algebraic relations, the product $\circ^{\smallint_{0}^{1}}$ on $\Pi_{1,0}$ reviewed in \S2.2.2 plays the role of an elementary substitute to a motivic Galois action. More precisely, for each field $K$ of characteristic zero and $\mu \in K$, one has \cite{Ra} a pro-affine scheme $\DMR_{\mu}$ of solutions $\varphi$ to the regularized double shuffle equations and the condition $\varphi[e_{0}e_{1}] = \mu$, such that :

\begin{Theorem} (Racinet, \cite{Ra}) The product $\circ^{\smallint_{0}^{1}}$ makes $\DMR_{0}$ into a group scheme and the left multiplication by $\circ^{\smallint_{0}^{1}}$ makes $\DMR_{\mu}$ into a torsor under this group scheme.
\end{Theorem}

\noindent One has similar results for the two other standard families of relations (see \cite{Dr} for associators, and \cite{F3} for a general review including a comparison with the motivic framework.)

\subsection{The crystalline realization of $\pi_{1}^{\un}(\mathbb{P}^{1} - \{0,1,\infty\})$ \label{crystalline aspects}}

\noindent Let $p$ be a prime number. The definition of the crystalline Frobenius $\phi$ of the De Rham $\pi_{1}^{\un}$ reviewed in \S\ref{paragraph crystalline Frobenius} applies to $X_{\mathbb{F}_{p}} = \mathbb{P}^{1} - \{0,1,\infty\} / \mathbb{F}_{p}$ and its lift $X_{\mathbb{Z}_{p}} = \mathbb{P}^{1} - \{0,1,\infty\} / \mathbb{Z}_{p}$ ; we will also denote by $X_{\mathbb{Q}_{p}} = \mathbb{P}^{1} - \{0,1,\infty\} / \mathbb{Q}_{p}$. When $k=\mathbb{F}_{p}$, the Frobenius of $k$ is $(x \mapsto x^{p})=\id_{\mathbb{F}_{p}}$, thus the Frobenius $\sigma$ of $W(\mathbb{F}_{p})=\mathbb{Z}_{p}$, which is the unique automorphism lifting it, is $\id_{\mathbb{Z}_{p}}$. Let $\alpha \in \mathbb{N}^{\ast}$, which will represent the number of iterations of $\phi$.

\subsubsection{The Frobenius of $\pi_{1}^{\un,\DR}(\mathbb{P}^{1} - \{0,1,\infty\} /\mathbb{Q}_{p})$ at $(-\vec{1}_{1},\vec{1}_{0})$ and $p$-adic multiple zeta values}

There are two different, but conjecturally arithmetically equivalent, notions of $p$-adic multiple zeta values. We have extended them in \cite{J1} \cite{J3} by considering the Frobenius iterated a certain number of times instead of the Frobenius itself. Introducing the number of iterations of the Frobenius also enables to formulate the two definitions in a unified way (see \S4.3). The first definition is the following (where $\tau$ is defined in \S2.2.2) :

\begin{Definition} (For $\alpha=1$, Deligne-Goncharov, \cite{DeGo}, \S5.28 ; for any $\alpha$, \cite{J1})
\newline Let, for $\alpha \in \mathbb{N}^{\ast}$, $\Phi_{p,\alpha}= \tau(p^{\alpha})\phi^{\alpha}({}_1 1_{0}) \in \pi_{1}^{\un,\DR}(\mathbb{P}^{1} -\{0,1,\infty\},-\vec{1}_{1},\vec{1}_{0})(\mathbb{Q}_{p})$. The following numbers are called $p$-adic multiple zeta values ($p$MZV's) \footnote{One sometimes encounters in the literature the variant of this definition without the factor $\tau(p^{\alpha})$}
$$ \zeta_{p,\alpha}(n_{d},\ldots,n_{1}) = (-1)^{d}\Phi_{p,\alpha}[e_{0}^{n_{d}-1}e_{1}\ldots e_{0}^{n_{1}-1}e_{1}] $$
\end{Definition}

\noindent These numbers provide a natural way of expressing the Frobenius at $(-\vec{1}_{1},\vec{1}_{0})$ :

\begin{Proposition} \label{description of Frobenius tangential} (Deligne-Goncharov, \cite{DeGo}, \S5.28)
\newline $\tau(p)\phi$ restricted to base points in the set $\{-\vec{1}_{1},\vec{1}_{0}\}$ is a $\mathbb{Q}_{p}$-point of $V_{\omega}$ from Definition \ref{def de V omega}. Thus, by Proposition \ref{prop Ihara action}, for all $\alpha \in \mathbb{N}^{\ast}$, $\tau(p^{\alpha})\phi^{\alpha}$ restricted to  $\pi_{1}^{\un}(\mathbb{P}^{1} - \{0,1,\infty\},-\vec{1}_{1},\vec{1}_{0})(\mathbb{Q}_{p})$ is the map $f \mapsto \Phi_{p,\alpha} \circ^{\smallint_{0}^{1}} f$.
\end{Proposition}

\noindent The second definition is based on the theory of Coleman integration, initiated in \cite{Co}, and generalized and formulated in a Tannakian way by Besser \cite{Bes} and Vologodsky \cite{V} ; one of its main result is the existence and uniqueness of Frobenius-invariant paths.

\begin{Definition} \label{second multiple zeta values}(Furusho, \cite{F1}, \cite{F2})
\newline Let $\Phi_{p}^{\KZ} \in \pi_{1}^{\un,\DR}(\mathbb{P}^{1} -\{0,1,\infty\},-\vec{1}_{1},\vec{1}_{0})(\mathbb{Q}_{p})$ be the unique path invariant by the Frobenius.
The following numbers are called $p$-adic multiple zeta values ($p$MZV's) 
$$ \zeta_{p}^{\KZ}(n_{d},\ldots,n_{1}) = (-1)^{d}\Phi_{p}^{\KZ}[e_{0}^{n_{d}-1}e_{1}\ldots e_{0}^{n_{1}-1}e_{1}] $$ 
\end{Definition}

\noindent Conjecturally, all versions of $p$-adic multiple zeta values satisfy the same algebraic relations, which should be the ones of multiple zeta values "modulo $\zeta(2)$", as reviewed in \S\ref{nature}.

\subsubsection{The Frobenius of $\pi_{1}^{\un,\DR}(\mathbb{P}^{1} - \{0,1,\infty\} / \mathbb{Q}_{p})$ on the fundamental torsor of paths starting at $\vec{1}_{0}$ and $p$-adic multiple polylogarithms \label{Frobenius torsor of paths}}

\noindent In general, to write the Frobenius of $\pi_{1}^{\un,\DR}(X_{K})$, one must consider an open affine covering of the rigid analytic space $\overline{X_{K}}^{\an}$ over $K$ and write the differential equation of Frobenius in coordinates on each open affine of the covering. Here, it is sufficient to consider the open affine $U^{\an} = \mathbb{P}^{1,\an} - \{|z-1|_{p}<1\}$, as suggested first by \cite{De2}, \S19.6.
\footnote{$U^{\an}$ contains $0$ and $\infty$ ; it is chosen instead of $\mathbb{P}^{1} - \{|z-\infty|_{p}<1\}= \mathbb{Z}_{p}^{\an}$ which contains $0$ and $1$ because the good lift of Frobenius on it, $z \mapsto z^{p}$ is simpler and more adapted to our purposes than the good lift of Frobenius on $\mathbb{Z}_{p}^{\an}$, which is $z \mapsto z^{p}$ conjugated by $z \mapsto \frac{z}{z-1}$, the unique homography sending $(0,1,\infty) \mapsto (0,\infty,1)$} Let the four following $\mathbb{Q}_{p}$-algebras of functions on $U^{\an}$ : $A^{\dagger}(U^{\an})$, the algebra of overconvergent rigid analytic functions ; $A^{\rig}(U^{\an})$, the algebra of rigid analytic functions ; $A^{\Col}(U^{\an})$, the algebra of Coleman functions ; $A^{\loc}(U^{\an})$, the algebra of locally analytic functions for the $p$-adic topology. We have $A^{\dagger}(U^{\an}) \subset A^{\rig}(U^{\an}) \subset A^{\Col}(U^{\an}) \subset A^{\loc}(U^{\an})$ \footnote{In particular, the elements of these four algebras have convergent power series expansions on a neighbourhood of $0$. Elements of $A^{\Col}(U^{\an})$ are characterized uniquely by their series expansion at $0$. The values of elements of $A^{\rig}(U^{\an})$ at any point are described explicitly in terms of the series expansion at $0$, see \S4.1.3.}.

\begin{Definition} (for $\alpha=1$, \cite{F2}, Definition 2.13 ; for any $\alpha$, \cite{J1}) Let  $\Li_{p,\alpha}^{\dagger}$ be the map $z \mapsto \tau(p^{\alpha})\phi^{\alpha}({}_z 1_{0})$ on $U^{\an}$ (the fundamental torsors of paths starting at $\vec{1}_{0}$ being trivialized at $\vec{1}_{0}$, $\Li_{p,\alpha}^{\dagger}$ can be seen as an element of $\pi_{1}^{\un,\DR}(X_{\mathbb{Q}_{p}},\vec{1}_{0})(A^{\dagger}(U^{\an})))$. The overconvergent $p$-adic multiple polylogarithms ($p$MPL$^{\dagger}$'s) are the functions $\Li_{p,\alpha}^{\dagger}[w] \in A^{\dagger}(U^{\an})$, for words $w$ on $\{e_{0},e_{1}\}$.
\end{Definition}

\noindent We have :

\begin{Proposition} The Frobenius $\tau(p^{\alpha})\phi^{\alpha}$ on $U^{an}$ (after trivialization at $\vec{1}_{0}$, we can say the Frobenius of $\pi_{1}^{\un,\DR}(X_{\mathbb{Q}_{p}},\vec{1}_{0})(A^{\Col}(U^{\an}))$) is
$$ f(z)(e_{0},e_{1}) \mapsto \Li_{p,\alpha}^{\dagger}(z)(e_{0},e_{1}). f(z^{p^{\alpha}})(e_{0},\Phi_{p,\alpha}^{-1}e_{1}\Phi_{p,\alpha}) $$
\end{Proposition}

\noindent In particular, the Frobenius iterated $\alpha$ times is entirely characterized by the couple $(\zeta_{p,\alpha},\Li_{p,\alpha}^{\dagger})$. 
\newline An essential theorem of Coleman integration is the existence and uniqueness of a primitive in algebras of Coleman functions up to an additive constant\footnote{which of course does not hold for algebras of locally analytic functions because the $p$-adic topology is totally disconnected}.

\begin{Definition} (Furusho, \cite{F1}, Definition 2.9) Let $\Li_{p}^{\KZ}$ be the unique Coleman function 
on $\mathbb{P}^{1} - \{0,1,\infty\}$ characterized by $\nabla_{\KZ}(\Li_{p}^{\KZ})=0$, $\Li_{p}^{\KZ}[w](0)=0$ for $w$ a word of the form $ue_{1}$, and $\Li_{p}^{\KZ}(z)[e_{0}]= \log_{p}(z)$. 	\footnote{Thus $\Li_{p}^{\KZ}$ is relative to the choice of a determination of $\log_{p}$ (we also have $\Li_{p}^{\KZ}(z) = - \log_{p}(1-z)$). However, $p$-adic multiple zeta values and overconvergent $p$-adic multiple polylogarithms do not depend on such a choice. The statements reviewed in this text are true for any branch of $\log_{p}$.} The $p$-adic multiple polylogarithms ($p$MPL's) are the Coleman functions $\Li_{p}^{\KZ}[w]$ where $w$ is a word on $\{e_{0},e_{1}\}$.
\end{Definition}

\noindent This definition can be rewritten in terms of Frobenius-invariant paths on the $\pi_{1}^{\un}$ (\cite{F2}, Theorem 2.3). By this definition, the power series expansion of $\Li_{p}^{\KZ}$ at $0$ is also expressed in terms of multiple harmonic sums, as in Proposition \ref{mhs and mpl} :

\begin{Proposition} \label{p-adic mhs and mpl} (Furusho, \cite{F1}, Lemma 2.7 and Remark 2.10 (1)) For $|z|_{p}<1$,  for $n_{d},\ldots,n_{1} \in \mathbb{N}^{\ast}$, the right-hand side below is absolutely convergent and we have :
\begin{equation} \label{eq:above}\Li_{p}^{\KZ}(z)[e_{0}^{n_{d}-1}e_{1}\ldots e_{0}^{n_{1}-1}e_{1}] = \sum_{0<m_{1}<\ldots<m_{d}} \frac{z^{m_{d}}}{m_{1}^{n_{1}} \ldots m_{d}^{n_{d}}}
\end{equation}
\end{Proposition}

\noindent The topological field $\mathbb{C}_{p}$ is totally disconnected and the disks $\{|z|_{p}<1\}$ and $\{|z-1|_{p}<1\}$ are disjoint. It does not make sense to take the limit $z \rightarrow 1$ in equation (\ref{eq:above}), and we see why an immediate computation of $p$-adic multiple zeta values from (\ref{eq:above}) is prevented. This contrasts with the complex setting where the limit $z \rightarrow 1$ of the analogous equation (\ref{eq: series expansion polylog}) gives the formula (\ref{eq:1}) for multiple zeta values as sums of series.
\newline The fact that the Frobenius commutes to the connection (\S\ref{paragraph crystalline Frobenius}) is equivalent to the following differential equation, which is sometimes referred to as the equation of horizontality of the Frobenius.

\begin{Proposition} \label{prop differential equation of Frobenius}(for $\alpha=1$, \cite{De2}, \S19.6 ; \cite{F2}, Theorem 2.14 ; for any $\alpha$, \cite{J1}) We have the following equivalent equations :
\begin{equation} 
\phi^{\alpha}(\Li_{p}^{\KZ}) = \Li_{p}^{\KZ} \end{equation}
\begin{equation} \label{eq:horizontality1}
\Li_{p,\alpha}^{\dagger}(z)(e_{0},e_{1})\text{ }.\text{ }
\Li_{p,X_{K}}^{\KZ}(z^{p^{\alpha}})\big(e_{0},\Phi_{p,\alpha}^{-1}e_{1}\Phi_{p,\alpha} \big) = \Li_{p,X_{K}}^{\KZ}(z)(p^{\alpha}e_{0},p^{\alpha}e_{1})
\end{equation} 
\begin{equation} \label{eq:horizontality equation}
d\Li_{p,\alpha}^{\dagger} = \bigg( p^{\alpha}\omega_{0}(z) e_{0} + \omega_{1}(z) e_{1} \bigg) \Li_{p,\alpha}^{\dagger} - \Li_{p,\alpha}^{\dagger} \bigg(  \omega_{0}(z^{p^{\alpha}}) e_{0} + \omega_{1}(z^{p^{\alpha}}) \Phi^{-1}_{p,\alpha}e_{1}\Phi_{p,\alpha} \bigg) 
\end{equation}
\end{Proposition}

\noindent The formal properties of Frobenius imply finally :
\begin{Proposition} (\cite{U2}, \S5.2) We have, with $e_{0}+e_{1}+e_{\infty}=0$,
\begin{equation} \label{eq: 1 infinity} e_{0} + \Phi_{p,\alpha}^{-1}e_{1}\Phi_{p,\alpha} + 
	\Li_{p,\alpha}^{\dagger}(\infty)^{-1}e_{\infty}\Li_{p,\alpha}^{\dagger}(\infty) = 0 
\end{equation}
\end{Proposition}

\noindent Since $\Li_{p,\alpha}^{\dagger}[w] \in A^{\dagger}(U^{\an})$, it is possible to express $\Li_{p,\alpha}^{\dagger}[w](\infty)$ in terms of the power series expansion of $\Li_{p,\alpha}^{\dagger}[w]$ at $0$. Through equation (\ref{eq: 1 infinity}), this can lead to a computation of $\zeta_{p,\alpha}$. This is the strategy initiated very briefly in \cite{De2}, \S19.6 in depth one, and followed in depth one and two in \cite{U2}.

\subsubsection{Known results on the algebraic theory of $p$-adic multiple zeta values \label{p-adic double shuffle relations}}

The application of the Tannakian motivic framework to $p$-adic multiple zeta values, evoked in \S1.4.1, follows from work of Yamashita (unpublished, but reviewed and used in \cite{Yam}) who constructed a crystalline realization functor of the categories of mixed Tate motives (which was asked in \cite{DeGo}, \S5.28). A byproduct of this construction is a morphism from the $\mathbb{Q}$-algebra of Goncharov's motivic multiple zeta values defined in \cite{Go1} to the one of $p$-adic multiple zeta values $\zeta_{p}^{\KZ}(w)$, which factors through the ring of $p$-adic periods $B$ evoked in \S1.4.1.
\newline It is known that $p$-adic multiple zeta values satisfy the regularized double shuffle relations \cite{BF}, \cite{FJ}, and that they satisfy the associator relations \cite{U1}. Thus, (by \cite{AET}, \cite{AT}, \cite{Sch}) they also satisfy Kashiwara-Vergne relations.

\section{Principles of the study and notations}

We explain a few principles (\S3.1) and facts, definitions and notations (\S3.2, \S3.3), about multiple harmonic sums and $\pi_{1}^{\un}(\mathbb{P}^{1} - \{0,1,\infty\})$, on which our work will rely.

\subsection{Three frameworks of computations : $\int_{0}^{1}$, $\int_{0}^{z<<1}$ and $\Sigma$ ; and the comparison between them}

\subsubsection{Two types of computations associated to relations between multiple polylogarithms}

The Conjecture \ref{conjecture Kaneko-Zagier} on finite multiple zeta values leads to ask why, and how, certain sequences of multiple harmonic sums - perhaps not only the ones defining finite multiple zeta values - should be equivalent to certain periods. According to the spirit of Grothendieck's period conjecture, we must ask ourselves how multiple harmonic sums are "related to geometry". A basic answer is given by \S\ref{deux point trois point deux} : multiple harmonic sums are essentially the coefficients of the power series expansions of multiple polylogarithms at $0$, both in the complex case (Proposition \ref{mhs and mpl}) and in the $p$-adic case (Proposition \ref{p-adic mhs and mpl}), the connection $\nabla_{\KZ}$ being of course the same in the two frameworks.
\newline Most of the algebraic relations of multiple zeta values, which are conjecturally reflected in finite multiple zeta values by Conjecture \ref{conjecture Kaneko-Zagier}, are actually consequences of equations satisfied by multiple polylogarithms : this is for example the case for double shuffle relations, and also for associator relations. Thus, translating these equations on the coefficients of the power series expansions of multiple polylogarithms at $0$ should give certain equations on multiple harmonic sums, and it would not be surprising that they explain at least partially the analogy between real and finite multiple zeta values.

\begin{Principle} \label{main principle}With any equation of multiple polylogarithms, one can associate two kinds of computations, and they may not necessarily be trivially equivalent to each other
\newline i) one by taking limits at tangential base-points (in our case, $z \rightarrow 1$).
\newline ii) one by considering coefficients of the series expansions at the origin of paths of integration (in our case, at $0$).
\end{Principle}

\noindent This simple principle combined to ideas of $p$-adic analysis will actually guide us both in the parts I and II.

\subsubsection{Computations on multiple harmonic sums without any reference to $\pi_{1}^{\un}$ and comparison}

Multiple harmonic sums are a computational bridge between complex and $p$-adic multiple zeta values. Any result of an explicit computation of $p$-adic multiple zeta values based on the differential equation satisfied by the Frobenius (Proposition \ref{prop differential equation of Frobenius}), should be necessarily expressed in terms of multiple harmonic sums, as it follows from the discussion in \S\ref{Frobenius torsor of paths}.  This motivates to study multiple harmonic sums intrinsically, and not only with the perspective of finding algebraic equations evoked above. They make sense by themselves, independently from the $\pi_{1}^{\un}$ ; given a property of multiple harmonic sums arising from the $\pi_{1}^{\un}$, is there a natural way to retrieve it without making any reference to the $\pi_{1}^{\un}$ ? And vice-versa ?
\newline We will keep these questions in mind, since it is regularly fruitful to compute an object in two different ways and to compare the results of the two computations, as soon as the comparison is not immediate. This general principle, of comparing two different incarnations of a mathematical object, is also at the center of the notion of periods.

\subsubsection{$\int_{0}^{1}$, $\smallint_{0}^{z<<1}$ and $\Sigma$}

The computations will lead us to define several objects, and we will add to the notations symbols reflecting the type of computation from which they are originated. We will add the superscript $\Sigma$ for refering to multiple harmonic sums viewed as iterated sums and $\int$ for refering to iterated integrals. It follows from the discussions of \S3.1.1 and \S3.1.2 that we actually have three different frameworks of computations :

\begin{Convention} In the next parts, we will add to the notations for the new objects :
\begin{itemize} \item the symbol $\int_{0}^{1}$ for objects arising from computations on iterated integrals at the tangential base-points $(\vec{1}_{0},-\vec{1}_{1})$
\item the symbol $\smallint_{0}^{z<<1}$ for objects arising from the computations involving power series expansion of multiple polylogarithms at $0$, and thus multiple harmonic sums viewed in this way via equation (\ref{eq:above})
\item the symbol $\Sigma$ for objects arising from the computations on multiple harmonic sums viewed as elementary finite iterated sums made without any reference to $\pi_{1}^{\un}(\mathbb{P}^{1} - \{0,1,\infty\})$.
\end{itemize}
\end{Convention}

\noindent We will compare the respective outcomes of these three types of computations to each other when it is possible.

\subsection{Framework for weighted multiple harmonic sums}

\begin{Convention} We will add the subscript $_{\har}$ to the notations for objects refering to multiple harmonic sums.
\end{Convention}

\subsubsection{The variant $\mathbb{Q}_{p}\langle\langle e_{0},e_{1} \rangle\rangle_{\har}$ of $\mathbb{Q}_{p} \langle\langle e_{0},e_{1} \rangle\rangle$  \label{paragraph generating series mhs} \label{paragraph notations for harmonic sums} ; non-commutative generating series of multiple harmonic sums}

Equation (\ref{eq:M mhs and mpl}) indicates that the correspondence $(n_{d},\ldots,n_{1}) \leftrightarrow e_{0}^{n_{d}-1}e_{1} \ldots e_{0}^{n_{1}-1}e_{1}$, valid for the indices of multiple zeta values and multiple polylogarithms, must be replaced by $(n_{d},\ldots,n_{1}) \leftrightarrow (e_{0}^{l}e_{1}e_{0}^{n_{d}-1}e_{1}\ldots e_{0}^{n_{1}-1}e_{1})_{l \in \mathbb{N}}$, for multiple harmonic sums. It indicates further
\footnote{by combining it with the fact that we have $\Li_{p}^{\KZ}(z)[e_{0}^{l-1}e_{1}e_{0}^{n_{d}-1}e_{1}\ldots e_{0}^{n_{1}-1}e_{1}e_{0}^{r}] \equiv \sum_{\substack{r_{d+1},\ldots,r_{1}\geq 0 \\ r_{d+1}+\ldots+r_{1}=r}}  {-l_{f} \choose r_{d+1}}\prod_{i=1}^{d} {-s_{i} \choose r_{i}} 
\Li_{p}^{\KZ}(z)[e_{0}^{l-1+r_{d+1}}e_{1}e_{0}^{n_{d}+r_{d}-1}e_{1}\ldots e_{0}^{n_{1}+r_{1}-1}e_{1}]$ modulo $\log_{p}(z)= \Li_{p}^{\KZ}[e_{0}](z)$, which follows from the shuffle equation for $\Li_{p}^{\KZ}(z)$}
an extension of the definition of multiple harmonic sums to words whose furthest to the right letter is $e_{0}$ as follows : let $l_{f}$ be a formal variable ; for $r \in \mathbb{N}$, we define
$$ \har_{m}(n_{d},\ldots,n_{1};r)
= \sum_{\substack{r_{d+1},\ldots,r_{1}\geq 0 \\ r_{d+1}+\ldots+r_{1}=r}}  {-l_{f} \choose r_{d+1}}\prod_{i=1}^{d} {-n_{i} \choose r_{i}} 
\har_{m} \big( n_{d}+r_{d},\ldots,n_{1}+r_{1} \big) $$
\noindent In the next definition, the i) will be motivated by \S4.3 (part I-3) :

\begin{Definition} For any ring $A$, 	
\newline i) Let $A\langle\langle e_{0},e_{1} \rangle\rangle_{\har}^{\smallint_{0}^{1}}
	= \{ f \in A\langle\langle e_{0},e_{1} \rangle\rangle
	\text{ | } \text{ for all words w}, \text{ for all } l\geq 1, f[e_{0}^{l}w]= 0   \}$
	\newline ii) Let $A \langle \langle e_{0},e_{1} \rangle \rangle_{\har}^{\smallint_{0}^{z<<1}} = 
	\{ f \in A\langle\langle e_{0},e_{1} \rangle\rangle
	\text{ | for all words w}, f[e_{0}^{l}w]\text{ is independent of l} \}$
	\newline iii) Let $A \langle \langle e_{0},e_{1} \rangle \rangle_{\har}^{\Sigma} = \prod_{\substack{d \in \mathbb{N} \\ n_{d},\ldots,n_{1} \in \mathbb{N}^{\ast} \\ r \in \mathbb{N}}} A.(n_{d},\ldots,n_{1};r)$ \footnote{since multiple harmonic sums $\har_{m}(n_{d},\ldots,n_{1})$ are $0$ as soon as $d>m$, a variant would be to replace $\prod$ by $\oplus$ in iii), and to make analogous modifications in i) and ii)}
\end{Definition}

\noindent These three modules are obviously isomorphic to each other, and we will use the notation $A\langle\langle e_{0},e_{1} \rangle\rangle_{\har}$ to refer to all of them at the same time.
\newline 
\noindent A function $w \mapsto \har_{n}(w)$ can  thus be seen as an element of $\mathbb{Q}_{p}[l_{f}]\langle\langle e_{0},e_{1} \rangle\rangle_{\har}$ ; for most purposes, it is however sufficient replace it by its image, which is in $\mathbb{Q}_{p}\langle\langle e_{0},e_{1} \rangle\rangle_{\har}$, by the weight-adically continuous linear map
$\mathbb{Q}_{p}[l_{f}]\langle\langle e_{0},e_{1} \rangle\rangle \rightarrow \mathbb{Q}_{p}[l_{f}]\langle\langle e_{0},e_{1} \rangle\rangle$ which sends all words of the form $u e_{1}$ to themselves and all other words to $0$. We adopt this point of view in the rest of this text.

\begin{Convention-Notation} \label{sequence of har n}\noindent
Sequences of weighted multiple harmonic sums are viewed and denoted in the following way :
\newline i) For any word $w$, and for any $I \subset \mathbb{N}$, we denote by $\har_{I}(w) = (\har_{n}(w))_{n\in I}$.  
\newline ii) For any $n \in \mathbb{N}^{\ast}$, 
$\har_{n}$ is the sequence $(\har_{n}(w))_{w\text{ word}}$ viewed in the sense of the discussion above as an element of $\mathbb{Q}_{p}\langle\langle e_{0},e_{1} \rangle\rangle_{\har}$.
\newline iii) For any $I \subset \mathbb{N}$, we denote by $\har_{I}$ the sequence $(\har_{n}(w))_{n \in I,\text{ }w\text{ word}}$, and view it as an element of $\Map(I,\mathbb{Q}_{p}\langle\langle e_{0},e_{1} \rangle\rangle_{\har})$.
\end{Convention-Notation}

\subsubsection{The first term of the $p$-adic expansion of prime weighted multiple harmonic sums \label{convergence of series}}

The prime weighted multiple harmonic sums (Definition \ref{definition weighted multiple harmonic sums}) will play a special role. The lower bound on their valuations below will be used constantly.

\begin{Fact} \label{reduction to finite multiple zeta values}Let $p$ a prime number, $\alpha \in \mathbb{N}^{\ast}$, and a word $w=(n_{d},\ldots,n_{1}$).
\newline We have $v_{p}(\har_{p^{\alpha}}(w)) \geq \weight(w)$ and the reduction of 
$p^{-\weight(w)}\har_{p^{\alpha}}(w)$ modulo $p$ is independent of $\alpha$. In particular, the image of  $\big(p^{-\weight(w)}\har_{p^{\alpha}}(w)\big)_{p \in \mathcal{P}} \in \prod_{p}\mathbb{Z}_{p}$ in $\mathbb{Z}/_{p \rightarrow \infty}$ by reduction modulo infinitely large primes is the finite multiple zeta value $\zeta_{\mathbb{Z}/_{p \rightarrow \infty}}(w)$. 
\footnote{Indeed, for any $\alpha$, the iterated sum 
$\har_{p^{\alpha}}(n_{d},\ldots,n_{1}) =(p^{\alpha})^{n_{d}+\ldots+n_{1}} \sum_{0<m_{1}<\ldots<m_{d}<p^{\alpha}}\frac{1}{m_{1}^{n_{1}}\ldots m_{d}^{n_{d}}}$ is the sum of the terms $(m_{1},\ldots,m_{d})\in (p^{\alpha-1}\mathbb{N})^{d}$, whose contribution is equal to $\har_{p}(n_{d},\ldots,n_{1}) \in p^{n_{d}+\ldots+n_{1}}\mathbb{Z}_{p}$, and the other terms, whose contribution is in  $p^{n_{d}+\ldots+n_{1}+1}\mathbb{Z}_{p}$.}
\end{Fact}

\noindent By this remark, all weight-adically convergent sums of $\mathbb{Q}$-linear combinations of sequences $\har_{I}(w)$ of prime weighted multiple harmonic sums $(I \subset \mathcal{P}^{\mathbb{N}^{\ast}})$ with reasonable rational coefficients are also $p$-adically convergent. Such objects will be omnipresent in this work.

\subsection{Topologies on $\pi_{1}^{\un,\DR}(\mathbb{P}^{1} - \{0,1,\infty\},-\vec{1}_{1},\vec{1}_{0})(\mathbb{Q}_{p})$ \label{topologies}}

For any $n,d \in \mathbb{N}$, let $\mathcal{W}(e_{0},e_{1})$, resp. $\mathcal{W}_{n}(e_{0},e_{1})$, $\mathcal{W}_{\ast,d}(e_{0},e_{1})$, $\mathcal{W}_{n,d}(e_{0},e_{1})$ be the set of words over $\{e_{0},e_{1}\}$, resp. words of depth $d$, of weight $n$, of weight $n$ and depth $d$. We can view $\mathbb{Q}_{p}\langle \langle e_{0},e_{1}\rangle\rangle$ as the set $\Map(\mathcal{W}(e_{0},e_{1}),\mathbb{Q}_{p})$. Let $\Lambda$ and $D$ be formal variables.

\begin{Definition} 
\noindent \newline i) For $f \in \mathbb{Q}_{p}\langle \langle e_{0},e_{1}\rangle\rangle$, let $\mathcal{N}_{\Lambda,D}(f) = \sum_{s,d\geq 0} \sup_{w \in \mathcal{W}_{s,d}(e_{0},e_{1})} |f[w]|_{p} \Lambda^{s} D^{d} \in \mathbb{R}_{+}[[\Lambda,D]]$.
\newline ii)  For $f \in \mathbb{Q}_{p}\langle \langle e_{0},e_{1}\rangle\rangle$ whose restriction to each $\mathcal{W}_{\ast,d}(e_{0},e_{1})$ is bounded, let $\mathcal{N}_{D}(f) = \sum_{d\geq 0} \sup_{w \in \mathcal{W}_{\ast,d}(e_{0},e_{1})} |f[w]|_{p}D^{d} \in \mathbb{R}_{+}[[D]]$.
\end{Definition}

\noindent One can associate with $\mathcal{N}_{\Lambda,D}$, resp. $\mathcal{N}_{D}$ the topology on $\Map(\mathcal{W}(e_{0},e_{1}),\mathbb{Q}_{p})$defined by the simple convergence, resp. uniform convergence on each  $\mathcal{W}_{\ast,d}(e_{0},e_{1})$. This also induces topologies on $\Pi_{1,0}(\mathbb{Q}_{p})$. We will see that these definitions are compatible with most usual operations on $\Pi_{1,0}(\mathbb{Q}_{p})$.

\section{(Part I) Computation of the Frobenius \cite{J1}, \cite{J2}, \cite{J3}}

Let $p$ be a prime number. In this part, we compute the Frobenius of $\pi_{1}^{\un,\DR}((\mathbb{P}^{1} - \{0,1,\infty\})/\mathbb{Q}_{p})$ reviewed in \S2.4. The three sections I-1, I-2, I-3 give three different computations, but depend on each other ; I-2 relies on I-1 and I-3 relies on I-2.

\subsection{(I-1) Direct solution to the equation of horizontality of the Frobenius \cite{J1}}

This first step is inspired from \"{U}nver's work in depth one and two \cite{U2}, who has also generalized it to $\mathbb{P}^{1} - \{0,\mu_{N},\infty\}$, $p\nmid N$ in \cite{U3}
\footnote{with a difference of convention ; in those papers, the Frobenius is $\phi^{-1}$ whereas here the Frobenius is $\tau(p^{\alpha}) \circ \phi^{\alpha}$ for all $\alpha \in \mathbb{N}^{\ast}$} : we will solve the system of equations on $(\Li_{p,\alpha}^{\dagger},\zeta_{p,\alpha})$ formed by (\ref{eq:horizontality equation}) and (\ref{eq: 1 infinity}), initially by induction on the weight (\emph{in fine} by induction on the depth). Because of the motivation coming from finite multiple zeta values (\S1.4.2), we want moreover to express the result in terms of prime multiple harmonic sums.

\subsubsection{A space of convergent series defined via prime weighted multiple harmonic sums}

The main result will be expressed in terms of the following object.

\begin{Definition} i) Let 
\begin{equation} \label{eq: algebra Har}\widehat{\text{Har}}_{\mathcal{P}^{\mathbb{N}^{\ast}}} = \text{Vect}_{\mathbb{Q}} \big\{ \big(\sum_{L\in \mathbb{N}} F_{L} \prod_{\eta=1}^{\eta_{0}} \har_{p^{\alpha}}(w_{L,\eta})\big)_{(p,\alpha)}\text{ | } (\ast)\}\subset 
\prod_{(p,\alpha)\in \mathcal{P} \times \mathbb{N}^{\ast}} \mathbb{Q}_{p}
\end{equation}
\noindent where $(\ast)$ means that $(w_{L,1})_{L\in\mathbb{N}},\ldots,(w_{L,\eta_{0}})_{L\in\mathbb{N}}$ $(\eta_{0} \in \mathbb{N}^{\ast})$ are sequences of words on $\{e_{0},e_{1}\}$ satisfying $\sum_{\eta=1}^{\eta_{0}}\text{weight}(w_{L,\eta}) \rightarrow_{l \rightarrow \infty} \infty$ and $ \limsup_{L \rightarrow \infty}\sum_{\eta=1}^{\eta_{0}}\text{depth}(w_{L,\eta}) < \infty$, and $(F_{L})_{L\in \mathbb{N}}$ is a sequence of rational numbers (independent of $(p,\alpha)$), such that the series is convergent (see \S\ref{convergence of series}).
\newline ii) For all $(n,d) \in (\mathbb{N}^{\ast})^{2}$, let $$ \widehat{\text{Har}}_{\mathcal{P}^{\mathbb{N}^{\ast}}}^{n,d} \subset \widehat{\text{Har}}_{\mathcal{P}^{\mathbb{N}^{\ast}}} $$
\noindent be the subspace generated by the elements as above such that, for all $L\in \mathbb{N}$,
\newline $\sum_{\eta=1}^{\eta_{0}}\text{weight}(w_{L,\eta})\geq n$ and $\sum_{\eta=1}^{\eta_{0}}\text{depth}(w_{L,\eta})\leq d$. We have $\widehat{\text{Har}}_{\mathcal{P}^{\mathbb{N}^{\ast}}} = \cup_{(n,d) \in (\mathbb{N}^{\ast})^{2}} \widehat{\text{Har}}_{\mathcal{P}^{\mathbb{N}^{\ast}}}^{n,d}$.
\end{Definition}

\subsubsection{Main result}

The main result of this part is the following :
\newline
\newline \textbf{Theorem 1.} $p$-adic multiple zeta values at a word $w$ of weight $n$ and depth $d$, viewed as sequences $(\zeta_{p,\alpha}(w))_{(p,\alpha)\in \mathcal{P} \times \mathbb{N}^{\ast}}$, are elements of $\widehat{\text{Har}}_{\mathcal{P}^{\mathbb{N}^{\ast}}}^{n-d,d}$.
\newline Furthermore (see \S6.2 of \cite{J1} for a precise statement) :
\newline i) There is an explicit formula, inductive with respect to the depth, for $(\Li_{p,\alpha}^{\dagger},\zeta_{p,\alpha})$ where the maps $\Li_{p,\alpha}^{\dagger}(w)$ are viewed as the maps $m \in \mathbb{N}^{\ast} \mapsto \Li_{p,\alpha}^{\dagger}[w][z^{m}]\in \mathbb{Q}_{p}$, \footnote{where $[z^{m}]$ means the coefficient of degree $m$ in the power series expansion at $0$} themselves viewed as elements of a specific subalgebra (defined in terms of multiple harmonic sums) of the one of locally analytic functions $\mathbb{Z}_{p} \rightarrow \mathbb{Q}_{p}$.
\newline ii) There are explicit lower bounds for the $p$-adic valuations of the numbers $\Li_{p,\alpha}^{\dagger}[w][z^{m}]$, and of $p$-adic multiple zeta values $\zeta_{p,\alpha}(w)$, which depend only on ($p$,weight($w$), depth($w$)).

\subsubsection{Strategy of proof and example of explicit formula}

1) In order to solve the system of equations (\ref{eq:horizontality equation}), (\ref{eq: 1 infinity}), we must first express $\Li_{p,\alpha}^{\dagger}(\infty)$ in terms of the coefficients $\Li_{p,\alpha}^{\dagger}[z^{m}]$. A way to do it is to remark that $U^{\an}=\mathbb{P}^{1,\an} - \{|z-1|_{p}<1\}$ is the image by $z \mapsto \frac{z}{z-1}$ (the homography of $\mathbb{P}^{1}$ sending $(0,1,\infty)$ to $(0,\infty,1)$) of $\mathbb{Z}_{p}^{\an}=\mathbb{P}^{1,\an} - \{|z-\infty|_{p}<1\}$ ; by definition, the algebra of rigid analytic functions on $\mathbb{Z}_{p}^{\an}$ is $A^{\rig}(\mathbb{Z}_{p}^{\an}) = \{ \sum_{m=0}^{\infty} a_{m}z^{m} \text{ | }|a_{m}|_{p} \rightarrow_{m \rightarrow \infty} 0\}$. Then, the characterization of the functions $\mathbb{N} \rightarrow \mathbb{Q}_{p}$ which extend to continuous functions $\mathbb{Z}_{p} \rightarrow \mathbb{Q}_{p}$ in terms of their Mahler coefficients \cite{Mah} can be reformulated as an explicit description of $A^{\rig}(\mathbb{P}^{1,\an} - ]1[)$, which implies that the maps $m \in \mathbb{N}^{\ast} \subset \mathbb{Z}_{p} \mapsto \Li_{p,\alpha}^{\dagger}[w][z^{m}]$ are continuously interpolable on $\mathbb{Z}_{p}$, and the formula $\Li_{p,\alpha}^{\dagger}(\infty)= - \lim_{|m|_{p} \rightarrow 0} \Li_{p,\alpha}^{\dagger}(\infty)[z^{m}]$ (see Proposition 2 of \cite{U2} for another version of this result, proved differently.)
\newline 
\newline 2) We must also check that equation (\ref{eq: 1 infinity}) gives an expression of $\Phi_{p,\alpha}$ and $\Li_{p,\alpha}^{\dagger}(\infty)$ in terms of each other which is compatible with the depth filtration. This extends some results of \cite{U2}, \S5 in low depth. We prove that we can express $\Phi_{p,\alpha}$ in depth $\leq d$ and $\Li_{p,\alpha}^{\dagger}(\infty)$ in depth $\leq d+1$ in terms of each other by polynomial expressions. These results do not alterate lower bounds of valuations since these polynomials have coefficients in $\mathbb{Z}$.
\newline 
\newline 3) Combining the two previous steps gives a solution to (\ref{eq:horizontality equation}), (\ref{eq: 1 infinity}) ; however, the recursive expression of $(\Li_{p,\alpha}^{\dagger},\zeta_{p,\alpha})$ obtained in this way is expressed in terms of certain limits of $\mathbb{Q}_{p}$-linear combinations of finite iterated sums, and this is not enlightening enough. In order to make any progress, we must understand via the formulas why the maps $m \in \mathbb{N}^{\ast} \mapsto \Li_{p,\alpha}^{\dagger}[w][z^{m}]$ are continuously interpolable on $\mathbb{Z}_{p}$. For $w$ of depth $1$ \footnote{in depth one we have : $\Li_{p,\alpha}^{\dagger}[e_{0}^{n-1}e_{1}](z) = (p^{\alpha})^{n}\sum_{\substack{0<m \\ p^{\alpha} \nmid m}} \frac{z^{m}}{m^{n}}$} and for $w = e_{1}e_{0}^{n-1}e_{1}$, this is easy, and it gives back the following essentially known result,
\footnote{see \cite{Co}, equation (4), and \cite{W}}
\begin{Example} \label{example 4.2}For all $n \in \mathbb{N}^{\ast}$, we have
	$\zeta_{p,\alpha}(n) = \frac{1}{n-1} \sum_{l \in \mathbb{N}} {-n \choose l} B_{l} \text{ } \har_{p^{\alpha}}(n+l-1)$
\end{Example}
\noindent which also leads to a proof for other words $w$ of depth $2$, using a power series expansion of $m \mapsto \Li_{p,\alpha}^{\dagger}[w][z^{m}]$ when $|m|_{p} \rightarrow 0$. In \cite{U2}, \S5.15, it is explained (for $\alpha=1$) that each integration operator $f \mapsto \int f \omega$, with $\omega\in \{\frac{p^{\alpha}dz}{z}, \frac{p^{\alpha}dz}{z-1}, \frac{d(z^{p^{\alpha}})}{z^{p^{\alpha}}-1}\}$, appearing in the resolution of (\ref{eq:horizontality equation}) in low depth, can be replaced by a regularized variant ; i.e the linear combinations of iterated integrals of these forms can be replaced by the same linear combinations of regularized variants, in such a way that the evaluation at $\infty$ makes sense for each term of the linear combination. Then, in \cite{U2}, \S5.16-\S5.19, for certain $w$ of depth $2$ and $3$ (which leads to compute certain $\zeta_{p,\alpha}$ in depth $2$), it is proven that certain terms of the maps $m \in p\mathbb{N}^{\ast} \subset \mathbb{Z}_{p} \mapsto \Li_{p,\alpha}^{\dagger}[w][z^{m}]$, are analytic
\footnote{by an indirect method using an equation satisfied by the logarithm of the $p$-adic Gamma function ; the formulas of the type of Example \ref{example 4.2}, infinite sums of prime weighted multiple harmonic sums, do not appear in \cite{U2}, \cite{U3}}.
\newline What we prove, by induction, is that the maps $m \in \mathbb{N}^{\ast} \mapsto \Li_{p,\alpha}^{\dagger}[w][z^{m}]$ extend to locally analytic functions on $\mathbb{Z}_{p}$, with coefficients of a very specific type, made out of elements of $\widehat{\text{Har}}_{\mathcal{P}^{\mathbb{N}^{\ast}}}$ and of multiple harmonic sums of upper bound $m \in \{1,\ldots,p^{\alpha}-1\}$, and subject to certain bounds of valuations of a specific type.
\newline 
\newline 4) It remains to lift each step of the inductive process above to an explicit map : these steps are the dualizations (expressed in terms of shuffle Hopf algebras) of equations (\ref{eq:horizontality equation}) and (\ref{eq: 1 infinity}), and the computation of regularized $p$-adic iterated integrals associated with a word in the differential forms $\frac{p^{\alpha}dz}{z}$, $\frac{p^{\alpha}dz}{z-1}$, and $\frac{d(z^{p^{\alpha}})}{z^{p^{\alpha}}-1}$. A fourth step is the use of the formula for $\Li_{p,\alpha}^{\dagger}(\infty)$ in terms of $\Li_{p,\alpha}^{\dagger}[z^{m}]$, but it adds no complexity. We can also upgrade the induction on the weight into an induction on the depth, which is much more efficient algorithmically.

\subsection{(I-2) Indirect solution to the equation of horizontality of the Frobenius \cite{J2}}

We are going to see that there exists a way to compute the Frobenius which is entirely different from I-1 but relies on I-1. Let us try to follow the ideas of \S3, and to find an application of the Principle \ref{main principle} to the differential equation of Frobenius.

\subsubsection{The Frobenius viewed as an operation on multiple harmonic sums \label{Frobenius harmonically}}

We are going to view the Frobenius $\tau(p^{\alpha})\phi^{\alpha}$, as well as the map $z \mapsto z^{p^{\alpha}}$, lift to $U^{\an}$ of the Frobenius of $X_{\mathbb{F}_{p}}$, as transformations of multiple harmonic sums, via their action on $\Li_{p}^{\KZ}$ and via the expression of multiple harmonic sums in terms of the power series expansions of $\Li_{p}^{\KZ}$ of equation (\ref{eq:above}).
\newline The map $\mathcal{F} : z \mapsto z^{p^{\alpha}}$ induces the map on power series $\sum_{m\in\mathbb{N}} c_{m} z^{m} \mapsto \sum_{m\in\mathbb{N}} c_{m} z^{p^{\alpha}m}$, i.e. the map on sequences of coefficients $(c_{m})_{m \in \mathbb{N}} \mapsto (1_{p^{\alpha}\mathbb{N}}c_{\frac{m}{p^{\alpha}}})_{m \in \mathbb{N}}$. Let us consider its right inverse $(c_{m})_{m\in\mathbb{N}} \mapsto (c_{p^{\alpha}m})_{m\in\mathbb{N}}$ ; let us say that $\mathcal{F}$ "sends $\har_{m}$ to $\har_{p^{\alpha}m}$".
\newline The invariance of $\Li_{p}^{\KZ}$ by the Frobenius $\tau(p^{\alpha})\phi^{\alpha}$ is expressed by equation (\ref{eq:horizontality1}) ; let us take, in that equation, the coefficient of the series expansion at $0$ of degree $p^{\alpha}m$, with $m \in \mathbb{N}^{\ast}$, then, the coefficient a word $e_{0}^{l-1}e_{1}e_{0}^{n_{d}-1}e_{1} \ldots e_{0}^{n_{1}-1}e_{1}$ ; and, finally, let us apply $\tau(p^{\alpha}m)$. We get an equation of the form
\begin{multline} \label{eq:harmonic version of Frobenius} \har_{p^{\alpha}m}(n_{d},\ldots,n_{1}) = 
\text{ a function indexed by }(l,n_{d},\ldots,n_{1})\text{ applied to }
\\ (\Li_{p,\alpha}^{\dagger}, \zeta_{p,\alpha} ,m, (\har_{m'})_{1 \leq m' \leq m})
\end{multline}
\noindent let us say that $\tau(p^{\alpha})\phi^{\alpha}$ "sends $(\har_{m'})_{1 \leq m' \leq m}$ to $\har_{p^{\alpha}m}$".
\newline In particular, the right-hand side of (\ref{eq:harmonic version of Frobenius}), which depends a priori on $l$ in a complicated way, is actually a constant function of $l$. The idea leading to the next theorem is then to show that equation (\ref{eq:harmonic version of Frobenius}) can be simplified a lot if we consider the limit of the right-hand side when $l \rightarrow \infty$.

\subsubsection{Main result \label{main result of I-2}}

We define a condition of summability in $\mathbb{Q}_{p} \langle\langle e_{0},e_{1}\rangle\rangle$ :

\begin{Definition} 
\noindent\newline For any $A \subset \mathbb{Q}_{p}\langle \langle e_{0},e_{1} \rangle\rangle$, let 
$A_{S} = \{f \in A \text{ | for all }d \in \mathbb{N}^{\ast},\text{ } 
\sup_{w \in \mathcal{W}_{n,d}} |f[w]|_{p} \displaystyle \rightarrow_{n \rightarrow \infty} 0 \}$.
\newline Let also $\tilde{\Pi}_{1,0}(\mathbb{Q}_{p})_{S} \subset \Pi_{1,0}(\mathbb{Q}_{p})_{S}$ defined by the elements $f$ such that $f[e_{0}] = f[e_{1}] = 0$.
\end{Definition}

\noindent We can show that $\tilde{\Pi}_{1,0}(\mathbb{Q}_{p})_{S}$ is a subgroup of $\Pi_{1,0}(\mathbb{Q}_{p})$ both for the usual multiplication and $\circ^{\smallint_{0}^{1}}$, and that it is closed for both topologies of \S\ref{topologies}. We equip it, as well as $\mathbb{Q}_{p} \langle \langle e_{0},e_{1}\rangle\rangle_{\har}$, with the $\mathcal{N}_{D}$-topology. The topological group $(\tilde{\Pi}_{1,0}(\mathbb{Q}_{p})_{S},\circ^{\smallint_{0}^{1}})$ is isomorphic via the map $\Ad(e_{1})$
\footnote{since we read the groupoid multiplication of $\pi_{1}^{\un,\DR}(\mathbb{P}^{1} - \{0,1,\infty\})$ from the right to the left, we take the convention that $\Ad_{g}(h) = g^{-1}hg$} to the topological group $(\Ad_{\tilde{\Pi}_{1,0}(\mathbb{Q}_{p})_{S}}(e_{1}),\circ_{\Ad}^{\smallint_{0}^{1}})$ where $\circ_{\Ad}^{\smallint_{0}^{1}}$ is defined as $g \circ_{\Ad}^{\smallint_{0}^{1}} f = f(e_{0},g)$.
\newline 
\newline \textbf{Theorem 2}
\newline i) (framework $\smallint_{0}^{z<<1}$) There exists an explicit free continuous group action of $(\Ad_{\tilde{\Pi}_{1,0}(\mathbb{Q}_{p})_{S}}(e_{1}),\circ_{\Ad}^{\smallint_{0}^{1}})$, the $p$-adic pro-unipotent $\int_{0}^{z<<1}$-harmonic action of $\mathbb{P}^{1} - \{0,1,\infty\}$
$$ \circ_{\har}^{\smallint_{0}^{z<<1}} : \Ad_{\tilde{\Pi}_{1,0}(\mathbb{Q}_{p})_{S}}(e_{1}) \times \Map(\mathbb{N},\mathbb{Q}_{p}\langle \langle e_{0},e_{1}\rangle\rangle_{\har}^{\smallint_{0}^{z<<1}}) \rightarrow  \Map(\mathbb{N},\mathbb{Q}_{p}\langle \langle e_{0},e_{1}\rangle\rangle_{\har}^{\smallint_{0}^{z<<1}}) $$ 
\noindent such that we have 
\begin{equation} \label{eq:DR RT har} \har_{p^{\alpha}\mathbb{N}} = \Phi_{p,\alpha}^{-1}e_{1}\Phi_{p,\alpha} \circ_{\har}^{\smallint_{0}^{z<<1}} \har_{\mathbb{N}}
\end{equation}
\newline ii) (framework $\Sigma$) There exists an explicit continuous map, the $p$-adic pro-unipotent $\Sigma$-harmonic action of $\mathbb{P}^{1} - \{0,1,\infty\}$
$$ \circ_{\har}^{\Sigma} : (\mathbb{Q}_{p}\langle \langle e_{0},e_{1} \rangle\rangle_{\har}^{\Sigma})_{S} \times \Map(\mathbb{N},\mathbb{Q}_{p}\langle \langle e_{0},e_{1} \rangle\rangle_{\har}^{\Sigma}) \rightarrow \Map(\mathbb{N},\mathbb{Q}_{p}\langle \langle e_{0},e_{1} \rangle\rangle_{\har}^{\Sigma}) $$
\noindent such that we have 
\begin{equation} \label{eq:RT RT har}\har_{p^{\alpha}\mathbb{N}} = \har_{p^{\alpha}} \circ_{\har}^{\Sigma} \har_{\mathbb{N}} 
\end{equation}
\noindent iii) (comparison) By relating to each other the two formulas above, we get the definition of an explicit map $\comp^{\smallint \leftarrow \Sigma}$ satisfying
\begin{equation} \label{eq: RT sigma}
\Phi_{p,\alpha}^{-1}e_{1}\Phi_{p,\alpha} = \comp^{\smallint \leftarrow \Sigma} \har_{p^{\alpha}}
\end{equation}

\noindent 
\newline We call $\int_{0}^{z<<1}$-harmonic Frobenius the map $ (\phi^{\alpha})_{\har}^{\smallint_{0}^{z<<1}} : h \mapsto \Phi_{p,\alpha}^{-1}e_{1}\Phi_{p,\alpha} \circ_{\har}^{\smallint_{0}^{z<<1}} h$ and $\Sigma$-harmonic Frobenius the map  $(\phi^{\alpha})_{\har}^{\Sigma} : h \mapsto \har_{p^{\alpha}} \circ_{\har}^{\Sigma} h$. 
\newline We note that, by contrast with the formulas of I-1 which involves the five objects $\Phi_{p,\alpha}$, $\Phi_{p,\alpha}^{-1}e_{1}\Phi_{p,\alpha}$, $\Li_{p,\alpha}^{\dagger}(\infty)$, 
$\Li_{p,\alpha}^{\dagger}(\infty)^{-1}(e_{0}+e_{1})\Li_{p,\alpha}^{\dagger}(\infty)$, $m \mapsto \Li_{p,\alpha}^{\dagger}[z^{m}]$, and relations between them, the formula (\ref{eq:DR RT har}) above involves only $\Phi_{p,\alpha}^{-1}e_{1}\Phi_{p,\alpha}$, and similarly, the formula (\ref{eq:RT RT har}) above involves only $\har_{p^{\alpha}}$. Moreover, by contrast with the formulas of I-1 which are recursive, the three formulas (\ref{eq:DR RT har}), (\ref{eq:RT RT har}), (\ref{eq: RT sigma}) can be written without recursion and in a few lines (see \cite{J2}). Instead of solving a system of equations, as we did in I-1, here, we only write two different incarnations of the Frobenius and compare them to each other. The proof will however rely on I-1.

\subsubsection{The $m=1$ case of the main result}

Let us look at the $m=1$ term of (\ref{eq:DR RT har}) ; since $\har_{1}(w)=0$ for all non-empty words $w$, almost all the terms of the equation vanish and the explicit formula for $\circ_{\har}^{\smallint_{0}^{z<<1}}$ gives directly :
\newline 
\newline 
\textbf{Corollary 2} 
\newline For all $(n_{d},\ldots,n_{1})$, we have :

\begin{multline} \label{eq:corollary I 2 a} \har_{p^{\alpha}}(n_{d},\ldots,n_{1}) = (\Phi_{p,\alpha}^{-1}e_{1}\Phi_{p,\alpha})\big[\frac{1}{1-e_{0}}e_{1}e_{0}^{n_{d}-1}e_{1}\ldots e_{0}^{n_{1}-1}e_{1} \big]
\\ = \sum_{d'=0}^{d} \sum_{l_{d'+1},\ldots,l_{d} \in \mathbb{N}} \prod_{i=d'}^{d} {-n_{i} \choose l_{i}}  \zeta_{p,\alpha}(n_{d'+1}+l_{d'+1},\ldots,n_{d}+l_{d})\zeta_{p,\alpha}(n_{d'},\ldots,n_{1}) \end{multline}

\noindent This is an inversion of the  expansion of $p$-adic multiple zeta values in terms of prime weighted multiple harmonic sums arising from Theorems 1 and 2. The $\alpha=1$ case of (\ref{eq:corollary I 2 a}) solves a conjecture of Akagi, Hirose and Yasuda
\footnote{In depth one and $\alpha=1$, equation (\ref{eq:corollary I 2 a}) is equivalent to a result on Washington on the values of the Kubota-Leopoldt zeta function \cite{W}, via a result of Coleman \cite{Co}. In depth two and for $\alpha=1$, equation (\ref{eq:corollary I 2 a}) had been proved by M. Hirose, a priori by a different method.}. One can lift the formula of (\ref{eq:corollary I 2 a}) to the definition of a map $\comp^{\Sigma \leftarrow \smallint}$, which reformulates (\ref{eq:corollary I 2 a}) as :
\begin{equation} \har_{p^{\alpha}} = \comp^{\Sigma \leftarrow \smallint}
(\Phi_{p,\alpha}^{-1}e_{1}\Phi_{p,\alpha}) \label{eq: inv DR sigma}
\end{equation}
\noindent This map satisfies the following property, which re-implies equation (\ref{eq:corollary I 2 a}) through equation (\ref{eq: RT sigma}) :
\begin{equation} \label{eq: inv de sigma} \comp^{\Sigma \leftarrow \smallint} \circ \comp^{\smallint \leftarrow \Sigma} = \id 
\end{equation}
\noindent This last statement (\ref{eq: inv de sigma}) is used in the proof of iii) of the Theorem 2, and is obtained as a "lift" of the proof of (\ref{eq:RT RT har}) restricted to $m=1$.

\subsubsection{Application to finite multiple zeta values according to Akagi, Hirose and Yasuda \label{paragraph applications to finite}}

We now review, following \cite{AHY}, some consequences on finite multiple zeta values of Corollary 2 brought together with other results ; they answer to the question of \S1.4.3. The following integrality property of Furusho's version of $p$-adic multiple zeta values holds : for all words $w$, we have $\displaystyle\zeta_{p}^{\KZ}(w) \in \sum_{n \geq \weight(w)} \frac{p^{n}}{n!}\mathbb{Z}_{p}$ ; this is proved in \cite{AHY}, by the logarithmic generalization \cite{OFcrystals} of Mazur's theorem on crystalline cohomology \cite{Maz}, and in \cite{Cha}. It implies easily the same result for $\zeta_{p}=\zeta_{p,1}$. \footnote{and more generally for all $\zeta_{p,\alpha}$ ; there is a difference of convention : in \cite{AHY} the authors do not consider $\zeta_{p,1}(w)$ but $\zeta_{p}^{\text{De}}(w) = p^{-\weight(w)}\zeta_{p,1}(w)$ in accordance with the notations of \cite{F2}, and they do not consider $\har_{p}(w)$ but $p^{-\weight(w)}\har_{p}(w)$.} The integrality of $\zeta_{p}$ and our equation (\ref{eq:corollary I 2 a}) brought together imply, for all $p>n_{d}+\ldots+n_{1}$, that below the right-hand side is in $p^{n_{d}+\ldots+n_{1}}\mathbb{Z}_{p}$ and :
\begin{equation} \label{eq: reduction modulo primes}\har_{p}(n_{d},\ldots,n_{1}) \equiv \sum_{d'=0}^{d}(-1)^{n_{d'+1}+\ldots+n_{d}}\zeta_{p}(n_{d'+1},\ldots,n_{d})\zeta_{p}(n_{d'},\ldots,n_{1})
 \mod p^{n_{d}+\ldots+n_{1}+1}\mathbb{Z}_{p} 
\end{equation}
\noindent whence a formula for $\zeta_{\mathbb{Z}/_{p \rightarrow \infty}}(n_{d},\ldots,n_{1})$ in terms of $p$-adic multiple zeta values ; this formula coincides with the formula in Kaneko-Zagier's Conjecture \ref{conjecture Kaneko-Zagier}, and thus explains it.
\newline This also enables to recast Kaneko-Zagier's conjecture as the injectivity of the map of reduction of $p$-adic multiple zeta values "modulo infinitely large primes", which it is possible to define by the integrality property, and which is a surjection onto the $\mathbb{Q}$-algebra of finite multiple zeta values by (\ref{eq: reduction modulo primes}).
\newline Finally, the dimensions of the spaces of $p$-adic multiple zeta values being upper bounded by the dimensions of the motivic Galois group \cite{Yam} - this is a particular case of a central fact of the Galois theory of periods, see \S1.1.2 - the same bounds of dimensions for finite multiple zeta values are deduced by reduction modulo infinitely large primes : $\sum_{n=0}^{\infty} \dim(\mathcal{Z}_{n}^{\text{fin}})\Lambda^{n} \leq \frac{1-\Lambda^{2}}{1 - \Lambda^{2} - \Lambda^{3}}$ where $\mathcal{Z}_{n}^{\text{fin}}$ is the $\mathbb{Q}$-vector space generated by finite multiple zeta values of weight $n$. This is conjecturally an equality by Kaneko-Zagier's conjecture. In the end, the Galois theory of periods applies to finite multiple zeta values as if they were as periods, despite that they live in the unusual ring $\mathbb{Z}/_{p \rightarrow \infty} = \big( \prod_{p \in \mathcal{P}} \mathbb{Z}/p\mathbb{Z} \big) / \big( \bigoplus_{p \in \mathcal{P}} \mathbb{Z}/p\mathbb{Z} \big)$.

\subsubsection{Strategy of proof of i) and definition of $\circ^{\smallint_{0}^{z<<1}}_{\har}$ and $\comp^{\Sigma \leftarrow \smallint}$}

\noindent 1) Let us go back to equation (\ref{eq:harmonic version of Frobenius}). Having made the observation of \S4.2.1, we prove
\footnote{an intuition behind the next lemma is that the numbers $m^{\weight(w)}\Li_{p,\alpha}^{\dagger}[w][z^{m}]$ should be smaller $p$-adically than the numbers $\har_{m}(w)$ because of the difference between the respective radii of convergence of $\Li_{p,\alpha}^{\dagger}$ and $\Li_{p}^{\KZ}$.}

\begin{Key Lemma} (\cite{J2}, \S3)
	The terms of equation  (\ref{eq:harmonic version of Frobenius}) containing a factor $\Li_{p,\alpha}^{\dagger}$ tend to zero when $l \rightarrow \infty$.
\end{Key Lemma}

\noindent This is a consequence of the bounds of valuations of $\Li_{p,\alpha}^{\dagger}[w]$ in Theorem 1. By this lemma, taking the limit $l \rightarrow \infty$ in (\ref{eq:harmonic version of Frobenius}) is a big simplification of (\ref{eq:harmonic version of Frobenius}). By expressing the limit of the other terms, we have a formula for $\har_{p^{\alpha}m}$ in terms of $\har_{m}$ (instead of $(\har_{m'})_{1 \leq m' \leq m}$) and certain power series of $n$ whose coefficients are made out of $\zeta_{p,\alpha}$ (instead of $(\zeta_{p,\alpha},\Li_{p,\alpha}^{\dagger})$) : this is a preliminary version of equation (\ref{eq:DR RT har}) :

\begin{Example} \label{example DR RT depth two}In depth one, equation (\ref{eq:DR RT har}) is for all $s \in \mathbb{N}^{\ast}$,
	\begin{equation}\har_{p^{\alpha}m}(n) = \har_{m} (n) + \sum_{b \in \mathbb{N}} m^{b+n} (\Phi_{p,\alpha}^{-1}e_{1}\Phi_{p,\alpha})[e_{0}^{b}e_{1}e_{0}^{n-1}e_{1}]
	\end{equation}
\noindent In depth two, equation (\ref{eq:DR RT har}) is for all $n_{1},n_{2} \in \mathbb{N}^{\ast}$,
	\begin{multline}\har_{p^{\alpha}m}(n_{2},n_{1}) = \har_{m}(n_{2},n_{1}) + 
	\sum_{b\in \mathbb{N}} m^{b+n_{2}+n_{1}} (\Phi_{p,\alpha}^{-1}e_{1}\Phi_{p,\alpha})[e_{0}^{b}e_{1}e_{0}^{n_{2}-1}e_{1} e_{0}^{n_{1}-1}e_{1}] 
	\\ + \sum_{r_{2}=0}^{n_{2}-1} \har_{m}(n_{2}-r_{2}) m^{r_{2}+n_{1}} (\Phi_{p,\alpha}^{-1}e_{1}\Phi_{p,\alpha})[e_{0}^{r_{2}}e_{1} e_{0}^{n_{1}-1}e_{1}] 
	\\ + \sum_{r_{1}=0}^{n_{1}-1} \har_{m}(n_{1}-r_{1})  \sum_{b \in \mathbb{N}} m^{b+n_{2}+r_{1}}(\Phi_{p,\alpha}^{-1}e_{1}\Phi_{p,\alpha})[e_{0}^{b}e_{1}e_{0}^{n_{2}-1}e_{1} e_{0}^{r_{1}}]
	\end{multline}
\end{Example}
\noindent We also observe that specializing this equation to $n=1$ gives equation (\ref{eq:corollary I 2 a}).
\newline 
\newline 2) Let $\mathbb{Q}_{p} \langle\langle e_{0},e_{1} \rangle\rangle^{\lim} = \{ f \in \mathbb{Q}_{p}\langle\langle e_{0},e_{1} \rangle\rangle \text{ |  for all words }w,\text{ }(f[e_{0}^{l}w])_{l\in \mathbb{N}}\text{ has a limit when }l \rightarrow \infty\}$. We have a map
$\lim : \mathbb{Q}_{p} \langle \langle e_{0},e_{1} \rangle\rangle^{\lim} \rightarrow \mathbb{Q}_{p} \langle \langle e_{0},e_{1} \rangle\rangle_{\har}^{\smallint_{0}^{z<<1}}$ 
\noindent defined by, for all words $w$, $\displaystyle(\lim f)[w] = \lim_{l\rightarrow \infty} f[e_{0}^{l}w]$. We define $\circ_{\har}^{\smallint_{0}^{z<<1}}$ by the equation
$$ g \circ_{\har}^{\smallint_{0}^{z<<1}} (m \mapsto h_{m}) = \big(m \mapsto \lim \big( \Ad_{\tau(m)(g)}(e_{1}) \circ^{\smallint_{0}^{1}}_{\Ad} h_{m} \big)\big) $$
\noindent We prove that it is well-defined, that it is a group action, that it is continuous, and that it is free. Finally, equation (\ref{eq:corollary I 2 a}) lifts to the definition of a map $\comp^{\Sigma \leftarrow \smallint}$.

\subsubsection{Strategy of proof of ii) and iii) and definition of $\circ_{\har}^{\Sigma}$ and $\comp^{\Sigma \leftarrow \smallint}$\label{the proof of ii in I-2}}

1) We consider a multiple harmonic sum $\har_{p^{\alpha}m}(n_{d},\ldots,n_{1})$, and we partitionate the set its indices $\{0<m_{1}<\ldots<m_{d}<p^{\alpha}m\}$ into the $2^{d}$ subsets defined by $d$ conditions of the type $p^{\alpha}|m_{i}$ or $p^{\alpha}\nmid m_{i}$, where $i$ runs through $\{1,\ldots,d\}$. For each of these subsets, we write the Euclidean division of each $m_{i}$ by $p^{\alpha}$ : $m_{i}=p^{\alpha}u_{i}$ or $m_{i} = p^{\alpha}u_{i}+r_{i}$ with $r_{i} \in \{1,\ldots,p^{\alpha}-1\}$ ; in this second case, we write $m_{i}^{-n_{i}} = r_{i}^{-s_{i}} (1 + \frac{p^{\alpha}u_{i}}{r_{i}})^{n_{i}} = \sum_{l_{i}\geq 0} {-n_{i} \choose l_{i} } (\frac{p^{\alpha}u_{i}}{r_{i}})^{l_{i}}$.
\newline We then sum over the values of $r_{i} \in \{1,\ldots,p^{\alpha}-1\}$ and $u_{i}$. The sums over $r_{i}$ give directly prime weighted multiple harmonic sums $\har_{p^{\alpha}}$. The sums over $u_{i}$ are of the type $\sum_{0 \leq u_{i_{1}} \leq \ldots \leq u_{i_{r}} \leq m-1} u_{i_{1}}^{\pm l_{i_{1}}} \ldots u_{i_{r}}^{\pm l_{i_{r}}}$, $l_{i_{j}} \in \mathbb{N}$. A subsum indexed by an $u_{i_{j}}$ with a positive exponent, i.e. $+l_{i_{j}}$, can be expressed as a polynomial of the upper and lower bounds of its domain of summation, namely, $u_{i_{j-1}}$ and $u_{i_{j+1}}$ ; by iterating this result, we obtain that the sums $\sum_{0 \leq u_{i_{1}} \leq \ldots \leq u_{i_{r}} \leq m-1} u_{i_{1}}^{\pm l_{i_{1}}} \ldots u_{i_{r}}^{\pm l_{i_{r}}}$ are linear combinations of weighted multiple harmonic sums $\har_{m}$ with coefficients in the ring of $\mathbb{Q}$-polynomial functions on $m$. After inverting an absolutely convergent double sum, we obtain an expression of $\har_{p^{\alpha}m}$ of the type of equation (\ref{eq:RT RT har}) :
\begin{Example} \label{example of RT RT har}In depth one and two 
\footnote{With $\mathcal{B}_{b}^{l_{r},\ldots,l_{1}}$, defined by for all $n \in \mathbb{N}^{\ast}$, $\sum_{0\leq n_{1}<\ldots<n_{r}<n} n_{1}^{l} = \sum_{b=1}^{l_{1}+\ldots+l_{r}+1} \mathcal{B}_{b}^{l_{r},\ldots,l_{1}} n^{b}$ ($r \geq 1$, $l_{1},\ldots,l_{r},b \geq 0$ and $1 \leq b \leq l_{1}+\ldots+l_{r}+r$) ; example : $\mathcal{B}_{b}^{l} = \frac{1}{l+1} {l+1 \choose b} B_{b+1-l}$ where $B$ denotes Bernoulli numbers.
\newline The generalization $\mathcal{B}_{b}^{\pm l_{r},\ldots,\pm l_{1}}$ is defined similarly as certain coefficients arising from the map of localization above}, equation (\ref{eq:RT RT har}) is, respectively :
\begin{equation} \har_{p^{\alpha}m}(n) = \har_{m}(n) + \sum_{b \geq 1} m^{b+n} \sum_{l \geq b-1} {-n \choose l} \mathcal{B}_{b}^{l} \har_{p^{\alpha}}(n+l)
\end{equation}
\begin{multline} \label{eq: har RT RT in depth 2}
\har_{p^{\alpha}m}(n_{2},n_{1}) = 
\har_{m}(n_{2},n_{1}) +
\\ \sum_{t \geq 1} m^{n_{2}+n_{1}+t} 
	\sum_{l \geq t-1} \bigg[ {-n_{1} \choose l+n_{2}} \mathcal{B}_{t}^{l+n_{2},-n_{2}} - {-n_{2} \choose l+n_{1}} \mathcal{B}_{t}^{l+n_{1},-n_{1}} \bigg] \har_{p^{\alpha}}(n_{1}+n_{2}+t)  + \sum_{t \geq 1} m^{n_{1}+n_{2}+t}
	\bigg[
	\\ \sum_{\substack{l_{1},l_{2} \geq 0 \\ l_{1}+l_{2} \geq t-2}} 
	\mathcal{B}_{t}^{l_{2},l_{1}} 
	\prod_{i=1}^{2} {-n_{i} \choose l_{i}} \har_{p^{\alpha}}(n_{i}+l_{i}) +  \sum_{\substack{l_{1},l_{2} \geq 0 \\ l_{1}+l_{2} \geq t-1}}
	\mathcal{B}_{t}^{l_{1}+l_{2}} 
	\bigg( \prod_{i=1}^{2} {-n_{i} \choose l_{i}} \bigg) \har_{p^{\alpha}}(n_{2}+l_{2},n_{1}+l_{1})
	\bigg] 
\\ - m^{n_{2}+n_{1}} \bigg[ \sum_{l_{1} \geq n_{2}-1} \mathcal{B}_{n_{2}}^{l_{1}} {-n_{1} \choose l_{1}} \har_{p^{\alpha}}(n_{1}+l_{1}) 
	- \sum_{l_{2} \geq n_{1}-1} \mathcal{B}_{n_{1}}^{l_{2}} {-s_{2} \choose l_{2}} \har_{p^{\alpha}}(n_{2}+l_{2})  \bigg] +
	\\ \sum_{\substack{ 1 \leq t < n_{2} \\ l \geq t-1}}
	m^{n_{1}+t} \har_{m}(n_{2}-t) \mathcal{B}_{t}^{l} {-n_{1} \choose l} \har_{p^{\alpha}}(n_{1}+l) - \sum_{\substack{1 \leq t < n_{1} \\ l' \geq t-1}} 
	m^{n_{2}+t} \har_{m}(n_{1}-t) \mathcal{B}_{t}^{l'} {-n_{2} \choose l'} \har_{p^{\alpha}}(n_{2}+l')
\end{multline}
\end{Example}
\noindent 2) Let $\mathbb{Q}_{p}\langle\langle e_{0}^{\pm 1},e_{1}\rangle\rangle$ be the set of linear maps $\mathbb{Q}\langle e_{0}^{\pm 1},e_{1}\rangle \rightarrow \mathbb{Q}_{p}$ where $\mathbb{Q}\langle e_{0}^{\pm 1},e_{1}\rangle$ is the localization of the non-commutative ring $\mathbb{Q}\langle e_{0},e_{1}\rangle$ equipped with the concatenation product at the multiplicative part generated by $e_{0}$. We define through it a variant $\mathbb{Q}_{p}\langle\langle e_{0}^{\pm 1},e_{1}\rangle\rangle_{\har}$ in the way of \S\ref{paragraph generating series mhs} which contains, for each $m \in \mathbb{N}$, the generating sequence $\har_{m}^{\loc}$ of multiple harmonic sums $\har_{m}(\pm n_{d},\ldots,\pm n_{1})$ in the sense above.
\newline The first step of 1), of writing a $p$-adic expansion of $\har_{p^{\alpha}m}$ and summing over $r_{i}$'s and $u_{i}$'s, lifts in a natural way into a map, for which one can write an explicit formula, which we call a $p$-adic pro-unipotent $\Sigma$-harmonic action for $\mathbb{P}^{1} - \{0,1,\infty\}$ localized at the source :
$$ (\circ_{\har}^{\Sigma})_{\loc} : (\mathbb{Q}_{p}\langle \langle e_{0},e_{1} \rangle\rangle_{\har}^{\Sigma})_{S} \times \Map(\mathbb{N},\mathbb{Q}_{p}\langle \langle e_{0}^{\pm 1},e_{1} \rangle\rangle_{\har}) \rightarrow 
\Map(\mathbb{N},\mathbb{Q}_{p}\langle\langle e_{0},e_{1}\rangle\rangle_{\har}) $$
\noindent The second step of 1), of expressing the maps $n \mapsto \sum_{0 \leq u_{i_{1}} \leq \ldots \leq u_{i_{r}} \leq m-1} u_{i_{1}}^{\pm l_{i_{1}}} \ldots u_{i_{r}}^{\pm l_{i_{r}}}$ as $\mathbb{Q}[m]$-linear combinations of $m \mapsto \har_{m}$ lifts to a map for which it is also possible to write an explicit formula
$$ \loc^{\vee} : \Map(\mathbb{N},\mathbb{Q}_{p}\langle\langle e_{0},e_{1}\rangle\rangle_{\har}) \rightarrow \Map(\mathbb{N},\mathbb{Q}_{p}\langle\langle e_{0}^{\pm 1},e_{1}\rangle\rangle_{\har}) $$
\noindent Then we define 
$$\circ_{\har}^{\Sigma}= (\circ_{\har}^{\Sigma})_{\loc} \circ (\id \times \loc^{\vee})$$
\noindent The explicit formulas are shown in \cite{J2}, \S5. The formula for $\loc^{\vee}$ is indexed by the paths from the root to the leaves of a certain tree.
\newline 3) We define the map $\comp^{\Sigma\rightarrow \smallint}$ as follows : $(\comp^{\Sigma\rightarrow \smallint}\har_{m})[e_{0}^{b}e_{1}e_{0}^{n_{d}-1}e_{1}\ldots e_{0}^{n_{1}-1}e_{1}]$ is the coefficient of $m^{b + n_{d}+\ldots+n_{1}}\har_{m}(\emptyset)$ in equation (\ref{eq:RT RT har}), and similarly when we replace $(\har_{m})_{m \in \mathbb{N}}$ by any sequence $(h_{m})_{m \in \mathbb{N}}$.
\newline A first observation is that we have $\comp^{\Sigma\leftarrow \smallint} \circ \comp^{\smallint \leftarrow \Sigma} = \id$ ; this follows from that, by specializing the equation $\har_{\mathbb{N}}^{\loc} = \loc^{\vee} \har_{\mathbb{N}}$ to $m=1$, we get relations between the coefficients $\mathcal{B}$ from the fact that $\har_{m=1}(w)=0$ for all non-empty words $w$ ; this is hidden behind the fact that specializing equation (\ref{eq:RT RT har}) to $m=1$ yields only the tautological equality $\har_{p^{\alpha}}=\har_{p^{\alpha}}$.
\newline A second observation is that, under certain conditions on $g$ we have $g\text{ } \circ_{\har}^{\Sigma}\text{ }h =$
\newline $(\comp^{\smallint \leftarrow \Sigma}g)\text{ } \circ_{\har}^{\smallint_{0}^{z<<1}}\text{ } h$ : this follows essentially from the fact that $\circ_{\har}^{\smallint_{0}^{z<<1}}$ and $\circ_{\har}^{\Sigma}$ are built out of similar processes of taking subwords and quotient words.

\subsection{(I-3) The number of iterations of the Frobenius viewed as a variable \cite{J3}}

In I-1 and I-2, the number $\alpha\in \mathbb{N}^{\ast}$ of iterations of Frobenius was fixed. We now study the Frobenius iterated $\alpha$ times as a function of $\alpha$ viewed as a $p$-adic variable. More generally, for any $(\alpha_{0},\alpha) \in (\mathbb{N}^{\ast})^{2}$ such that $\alpha_{0}|\alpha$, we are going to consider  $\phi$ iterated $\alpha$ times as $\phi^{\alpha_{0}}$ iterated $\frac{\alpha}{\alpha_{0}}$ times.
\newline For prime weighted multiple harmonic sums $\har_{p^{\alpha}}$, the variable $p^{\alpha}$ has two roles : it reflects the number of iterations of the Frobenius $\alpha$ and it is the upper bound of the domain of summation of the underlying iterated sum. This relates these considerations to the study of the map $n \mapsto \har_{n}$ viewed as a function of the $p$-adic variable $n$, whose certain aspects (the regularization) played a central role in I-1. By this relation, it will be useful in \cite{J7} (part III-1) to reason in terms of $(\alpha_{0},\alpha)$ and not only in terms of $\alpha$.

\subsubsection{The Frobenius as a contraction mapping}

We are going to use again, but this time in a more essential way, the topological framework of \S3.3. One can prove that $\mathcal{N}_{\Lambda,D}$, and thus $\mathcal{N}_{D}$ as well, is sub-multiplicative with respect to the product $\circ^{\smallint_{0}^{1}}$ : for all $f,f' \in \Pi_{1,0}(\mathbb{Q}_{p})$, we have
\begin{equation} \label{eq:submultiplicativity} \mathcal{N}_{\Lambda,D}( f' \circ^{\smallint_{0}^{1}} f ) \leq \mathcal{N}_{\Lambda,D}(f') \times \mathcal{N}_{\Lambda,D}(f)
\end{equation}
\noindent This inequation, together with the completeness of $\Pi_{1,0}(\mathbb{Q}_{p})$ authorizes the following definition : below, the exponent $-1(\circ^{\smallint_{0}^{1}})$ means "inverse for the product $\circ^{\smallint_{0}^{1}}$" :

\begin{Definition} \label{def contracting}Let $\kappa \in \mathbb{Q}_{p}^{\ast}$ with $|\kappa|_{p}<1$. A map $\psi : \Pi_{1,0}(\mathbb{Q}_{p})\rightarrow \Pi_{1,0}(\mathbb{Q}_{p})$ is said to be a $\kappa$-contraction  if, for all $f,f' \in \Pi_{1,0}(\mathbb{Q}_{p})$, we have :
\begin{equation}
\label{eq:inequality contractance} \mathcal{N}_{\Lambda,D} \big( \psi(f')^{-1(\circ^{\smallint_{0}^{1}})} \circ^{\smallint_{0}^{1}} \psi(f) \big) (\Lambda,D) \leq \mathcal{N}_{\Lambda,D} \big( f'^{-1(\circ^{\smallint_{0}^{1}})} \circ^{\smallint_{0}^{1}} f \big)(\kappa \Lambda,D) 
\end{equation}
\end{Definition}

\noindent Indeed, by (\ref{eq:submultiplicativity}) and the fact that $\Pi_{1,0}(\mathbb{Q}_{p})$ is complete, the contractions in this sense satisfy the usual properties of contractions of complete metric spaces regarding fixed points : they have a unique fixed point and any sequence $(\phi^{a}(f))_{a \in \mathbb{N}^{\ast}}$, with $f \in \Pi_{1,0}(\mathbb{Q}_{p})$ tends to the fixed point when $a \rightarrow \infty$.
\newline One can see that, for all 
$(\lambda,f') \in \mathbb{G}_{m}(\mathbb{Q}_{p}) \times \Pi_{1,0}(\mathbb{Q}_{p})$, if $|\lambda|_{p}<1$, the map $f \mapsto f' \circ^{\smallint_{0}^{1}} \tau(\lambda)(f)$ is a contraction in this sense
\footnote{and the inequality (\ref{eq:inequality contractance}) is in this case an equality with $\kappa=\lambda$}. This implies, by the formula for the Frobenius of $\Pi_{1,0}(\mathbb{Q}_{p})$ (Proposition \ref{description of Frobenius tangential}) that the inverse of the Frobenius iterated $\alpha_{0}$ times, i.e. $\phi^{-\alpha_{0}}$, is a contraction in this sense with $\lambda = p^{\alpha_{0}}$ 
\footnote{$\phi^{-\alpha_{0}}$ is of the form above with $g =\Phi_{p,\alpha_{0}}^{-1(\circ^{\smallint_{0}^{1}})}$  and $\lambda = p^{\alpha_{0}}$.}. This was hidden behind the existence and uniqueness of a Frobenius-invariant path in $\Pi_{1,0}(\mathbb{Q}_{p})$ subjacent to Definition \ref{second multiple zeta values} (see also the proof of Proposition 3.1 in \cite{F2} for an implicit evocation of this fact). Denoting by $\Phi_{p,-\alpha} = \Phi_{p,\alpha}^{-1(\circ^{\smallint_{0}^{1}})}$ for $\alpha \in \mathbb{N}^{\ast}$, we thus have $\displaystyle \Phi_{p}^{\KZ} = \lim_{\alpha \rightarrow \infty} \Phi_{p,-\alpha}$ and $\displaystyle (\Phi_{p}^{\KZ})^{-1(\circ^{\smallint_{0}^{1}})} = \lim_{\alpha \rightarrow \infty} \Phi_{p,\alpha}$. We thus have for each $\alpha \in \mathbb{Z} \cup \{\pm \infty\} - \{0\}$ a notion of $p$-adic multiple zeta values associated with the Frobenius iterated $\alpha$ times, in a sense given by the discussion above when $\alpha= \pm \infty$, and this formulation unifies all the definitions of $p$-adic multiple zeta values. We can use the notations $\Phi_{p}^{\KZ}= \Phi_{p,-\infty}$ and $(\Phi_{p}^{\KZ})^{-1(\circ^{\smallint_{0}^{1}})}=\Phi_{p,\infty}$.

\subsubsection{Main result}

Below, $\Lambda$ and $\textbf{a}$ are formal variables ; we use the notations of \S\ref{main result of I-2} ; and, for $n \in \mathbb{N}^{\ast}$, $\pr_{n} : \mathbb{Q}_{p} \langle \langle e_{0},e_{1}\rangle\rangle \rightarrow \mathbb{Q}_{p} \langle\langle e_{0},e_{1} \rangle\rangle$ is the map of "projection onto the terms of weight $n$" i.e. the sequence $(\pr_{n})_{n \in \mathbb{N}}$ is characterized by : for all $f \in \mathbb{Q}_{p} \langle \langle e_{0},e_{1}\rangle\rangle$, and $\lambda \in \mathbb{Q}_{p}^{\ast}$, $\tau(\lambda)(f) = \sum_{n \in \mathbb{N}} \pr_{n}(f)\lambda^{n}$. The parts o) and i) below rely on Corollary 2.
\newline 
\newline 
\textbf{Theorem 3}
\newline Let $(\alpha_{0},\alpha) \in (\mathbb{N}^{\ast})^{2}$ such that $\alpha_{0} | \alpha$.
\newline o) (framework $\int_{0}^{1}$) There exists a continuous free group action of $(\Ad_{\tilde{\Pi}_{1,0}(\mathbb{Q}_{p})_{S}}(e_{1}),\circ_{\Ad}^{\smallint_{0}^{1}})$, the $p$-adic pro-unipotent $\int_{0}^{1}$-harmonic action of $\mathbb{P}^{1} - \{0,1,\infty\}$
$$ \circ_{\har}^{\smallint_{0}^{1}} :
\Ad_{\tilde{\Pi}_{1,0}(\mathbb{Q}_{p})_{S}}(e_{1}) \times \mathbb{Q}_{p} \langle\langle e_{0},e_{1} \rangle\rangle^{\smallint_{0}^{1}}_{\har} \rightarrow \mathbb{Q}_{p} \langle\langle e_{0},e_{1} \rangle\rangle^{\smallint_{0}^{1}}_{\har} $$
\noindent which extends to a map $\mathbb{Q}_{p} \langle\langle e_{0},e_{1} \rangle\rangle_{S}
\times
\mathbb{Q}_{p} \langle\langle e_{0},e_{1} \rangle\rangle^{\smallint_{0}^{1}}_{\har} \rightarrow \mathbb{Q}_{p} \langle\langle e_{0},e_{1} \rangle\rangle^{\smallint_{0}^{1}}_{\har}$ which we will denote in the same way, such that we have the analytic expansion (equality of functions of $\alpha \in \mathbb{N}^{\ast}$) :
\begin{equation}
\label{eq:first of I-3}\har_{p^{\mathbb{N}}} = \sum_{s \geq 0} \big( \pr_{s+1}\big(
\Phi_{p,\infty}^{-1}e_{1}\Phi_{p,\infty}\big) \circ_{\har}^{\smallint_{0}^{1}} \comp^{\smallint \rightarrow \Sigma}(\Phi_{p,-\infty}^{-1}e_{1}\Phi_{p,-\infty})\big) (p^{\mathbb{N}})^{s}
\end{equation}
\noindent i) (framework $\int_{0}^{1}$) There exists an explicit map, the $\int_{0}^{1}$-harmonic iteration of the Frobenius $$(\widetilde{\text{iter}}_{\har}^{\smallint_{0}^{1}})^{\textbf{a},\Lambda} : \Ad_{\tilde{\Pi}_{1,0}(\mathbb{Q}_{p})_{S}}(e_{1}) \rightarrow \mathbb{Q}_{p}[[\Lambda^{\textbf{a}}]][\textbf{a}](\Lambda)\langle\langle e_{0},e_{1}\rangle\rangle^{\smallint_{0}^{1}}_{\har}$$\noindent such that, the map $(\text{iter}_{\har}^{\smallint_{0}^{1}})^{\frac{\alpha}{\alpha_{0}},p^{\alpha_{0}}} : \tilde{\Pi}_{1,0}(\mathbb{Q}_{p})_{S} \rightarrow \mathbb{Q}_{p}\langle\langle e_{0},e_{1}\rangle\rangle^{\smallint_{0}^{1}}_{\har}$ defined as the post-composition of $(\widetilde{\text{iter}}_{\har}^{\smallint_{0}^{1}})^{\textbf{a},\Lambda}$ by the reduction modulo $(\textbf{a}-\frac{\alpha}{\alpha_{0}},\Lambda-p^{\alpha_{0}})$ satisfies, at words $w$ such that $\frac{\alpha}{\alpha_{0}} > \depth(w)$ :
\begin{equation} \label{eq:second of I-3} \har_{p^{\alpha}}(w) = (\iter_{\har}^{\smallint_{0}^{1}})^{\frac{\alpha}{\alpha_{0}},p^{\alpha_{0}}}(\Phi^{-1}_{p,\alpha_{0}}e_{1}\Phi_{p,\alpha_{0}})\text{ }(w)
\end{equation}
\noindent ii) (framework $\Sigma$) There exists an explicit map, the $\Sigma$-harmonic iteration of the Frobenius $$(\widetilde{\text{iter}}_{\har}^{\Sigma})^{\textbf{a},\Lambda} : (\mathbb{Q}_{p}\langle\langle e_{0},e_{1}\rangle \rangle_{\har}^{\Sigma})_{S} \rightarrow \mathbb{Q}_{p}[[\Lambda^{\textbf{a}}]][\textbf{a}](\Lambda)\langle\langle e_{0},e_{1}\rangle\rangle^{\smallint_{0}^{1}}_{\har} $$
\noindent such that, the map $(\text{iter}_{\har}^{\Sigma})^{\frac{\alpha}{\alpha_{0}},p^{\alpha_{0}}} : (\mathbb{Q}_{p}\langle\langle e_{0},e_{1}\rangle \rangle_{\har}^{\Sigma})_{S} \rightarrow \mathbb{Q}_{p}\langle\langle e_{0},e_{1}\rangle\rangle^{\smallint_{0}^{1}}_{\har}$ defined as the composition of $\widetilde{\text{iter}}_{\har,\Sigma}^{\textbf{a},\Lambda}$ by the reduction modulo  $(\textbf{a}-\frac{\alpha}{\alpha_{0}},\Lambda-p^{\alpha_{0}})$, satisfies,
\begin{equation} \label{eq:third of I-3}\har_{p^{\alpha}} = \iter_{\har,\Sigma}^{\frac{\alpha}{\alpha_{0}},p^{\alpha_{0}}} (\har_{p^{\alpha_{0}}})
\end{equation}
\noindent iii) (comparison) Let us fix $\alpha_{0}$, and view (\ref{eq:first of I-3}), (\ref{eq:second of I-3}), (\ref{eq:third of I-3}) as three expansions of a function $\mathbb{N}^{\ast} \rightarrow \mathbb{Q}_{p} \langle \langle e_{0},e_{1}\rangle\rangle_{\har}$ of $\frac{\alpha}{\alpha_{0}}$, in the ring $\mathbb{Q}_{p}[[(p^{\alpha_{0}})^\frac{\alpha}{\alpha_{0}} ]][\frac{\alpha}{\alpha_{0}}]\langle \langle e_{0},e_{1} \rangle\rangle_{\har}$. Then, the coefficients of these expansions are identical.
\newline 
\newline We thus have a natural way to compute the Frobenius-invariant path $\Phi_{p}^{\KZ}$ indirectly, i.e. by identifying the computations in the frameworks $\int_{0}^{1}$ and $\Sigma$ of a same object.
\newline We note that it is easy to check that we have an isomorphism of topological spaces with continuous group action between $\Map(\mathbb{N}, \mathbb{Q}_{p} \langle \langle e_{0},e_{1} \rangle\rangle^{\smallint_{0}^{1}}_{\har})$ equipped with $m \mapsto (\circ_{\har}^{\smallint_{0}^{1}}) \circ (\tau(m) \times \id)$ and $\Map(\mathbb{N}, \mathbb{Q}_{p} \langle \langle e_{0},e_{1} \rangle\rangle^{\smallint_{0}^{z<<1}}_{\har})$ equipped with $\circ_{\har}^{\smallint_{0}^{z<<1}}$.

\subsubsection{Strategy of proof and examples of explicit formulas}

o) The formula for the Frobenius on $\Pi_{1,0}(\mathbb{Q}_{p})$ (Proposition \ref{description of Frobenius tangential}) combined to the definition of $\Phi_{p}^{\KZ}$ as the Frobenius-invariant yields a formula expressing $\Phi_{p,\alpha}$ in terms of $\Phi_{p}^{\KZ}$ (see also \cite{F2}, Theorem 2.8). We apply to this formula, first, $\Ad(e_{1})$, which fits into a commutative diagram involving $\circ^{\smallint_{0}^{1}}$ and $\circ^{\smallint_{0}^{1}}_{\Ad}$, and then $\comp^{\Sigma\leftarrow\smallint} : \Ad_{\tilde{\Pi}_{1,0}(\mathbb{Q}_{p})_{S}}(e_{1}) \rightarrow \mathbb{Q}_{p}\langle \langle e_{0},e_{1}\rangle\rangle_{\har}^{\smallint_{0}^{1}}$, and we apply Corollary 2 (equation \ref{eq:corollary I 2 a}). We get the result, where we view $\comp^{\Sigma\leftarrow \smallint}$ as a surjective map onto 
$\mathbb{Q}_{p}\langle \langle e_{0},e_{1}\rangle\rangle_{\har}^{\smallint_{0}^{1}}$, where $\circ_{\Ad}^{\smallint_{0}^{1}}$ is extended as a map $\mathbb{Q}_{p}\langle\langle e_{0},e_{1}\rangle\rangle \times \mathbb{Q}_{p}\langle\langle e_{0},e_{1}\rangle\rangle \rightarrow \mathbb{Q}_{p}\langle\langle e_{0},e_{1}\rangle\rangle$ by the same formula with the one of \S4.2.2 and where $\circ_{\har}^{\smallint_{0}^{1}}$ is defined by
 $$ g \circ_{\har}^{\smallint_{0}^{1}} (\comp^{\Sigma \leftarrow \smallint} h) =  \comp^{\Sigma \leftarrow \smallint}( g \circ_{\Ad}^{\smallint_{0}^{1}} h) $$
\noindent We check that $\circ_{\har}^{\smallint_{0}^{1}}$ is well-defined, and that is a continuous free group action.
\begin{Example} In depth one, this gives : for all $n \in \mathbb{N}^{\ast}$,
	$$ \har_{p^{\alpha}}(n) = (\Phi_{p,-\infty}^{-1}e_{1}\Phi_{p,-\infty}) \big[ \frac{1}{1-e_{0}}e_{1}e_{0}^{n-1}e_{1} \big] + 
	\sum_{b \in \mathbb{N}} (p^{\alpha})^{b+n} (\Phi_{p,\infty}^{-1}e_{1}\Phi_{p,\infty}) \big[ e_{0}^{l}e_{1}e_{0}^{n-1}e_{1} \big]$$
\end{Example}
\noindent i) For $a \in \mathbb{N}^{\ast}$, and $\lambda \in \mathbb{Q}_{p}^{\ast}$ which is not a root of unity, let us call "adjoint iteration $a$ times weighted by $\lambda$" the map $g \in \Ad_{\tilde{\Pi}_{1,0}(\mathbb{Q}_{p})_{S}}(e_{1}) \mapsto g \circ_{\Ad}^{\smallint_{0}^{1}} \tau(\lambda)(g) \circ_{\Ad}^{\smallint_{0}^{1}} \ldots \circ_{\Ad}^{\smallint_{0}^{1}} \tau(\lambda^{a-1})(g)$. It can be written in terms of the standard multiplication of $\mathbb{Q}_{p}\langle\langle e_{0},e_{1}\rangle\rangle$, which enables to write its dual. Then, we can analyze how it depends on $a$, provided we restrict to $a$ greater than the depth of the word under consideration. The map of the statement is then obtained by applying $\comp^{\Sigma \leftarrow \smallint}$ and Corollary 2 to the result.
\newline 
\newline ii) We introduce in $(p^{\alpha})^{n_{d}+\ldots+n_{1}}\sum_{0<m_{1}<\ldots<m_{d}<p^{\alpha}} \frac{1}{m_{1}^{n_{1}}\ldots m_{d}^{n_{d}}}$, the following additional parameters : we write each $m_{i}$ as $p^{\alpha_{0} v_{i}}(p^{\alpha_{0}}u_{i}+r_{i})$, where $v_{i},u_{i} \in \mathbb{N}$ and $r_{i} \in \{1,\ldots,p^{\alpha_{0}}-1\}$. Then, we write
$\frac{1}{m_{i}^{n_{i}}} = \frac{1}{p^{\alpha_{0} v_{i}n_{i}}} \sum_{l_{i} \geq 0} {n_{i} \choose - l_{i}}  \frac{(p^{\alpha_{0}}u_{i})^{l_{i}}}{r_{i}^{n_{i}+l_{i}}}$. We first have to translate the condition $0<m_{1}<\ldots<m_{d}<p^{\alpha}$ in terms of the new parameters, $u_{i}, v_{i}, r_{i}$. Summing over $u_{i}$ and $r_{i}$ has some similarities with a step of the proof of Theorem 2 in \ref{the proof of ii in I-2}; summing over $v_{i}$ is different since they appear as exponents. The result is obtained after having summed over all additional parameters.
\begin{Example} In depth one, this gives : for all $n \in \mathbb{N}^{\ast}$,
	$$ \har_{p^{\alpha}}(n) = \sum_{b \geq 1} \frac{p^{\alpha(n+b)}-1}{p^{\alpha_{0}(n+b)}-1} \sum_{L\geq -1} \mathcal{B}^{L+b}_{b} \har_{p^{\alpha_{0}}}(n+b+L) $$
\end{Example}
\noindent iii) the uniqueness of the expansion follows from a standard elementary argument.

\section{(Part II) Algebraic relations \cite{J4}, \cite{J5}, \cite{J6}}

In this part, we want to understand explicitly the algebraic relations of $p$-adic multiple zeta values, using the computations of part I. We will apply again the principles of \S3.
\newline We keep in mind the concepts of the Galois theory of periods and, at the same time, the fact that the formulas exchanging $\zeta_{p,\alpha}$ and $\har_{p^{\alpha}}$ involve (weight-adically and $p$-adically convergent) infinite summations and the group $\tilde{\Pi}_{1,0}(\mathbb{Q}_{p})_{S}$ defined by bounds of valuations ; our framework is thus not strictly speaking algebraic, although it comes from an algebraic framework.

\subsection{(II-1) Standard algebraic relations of adjoint multiple zeta values and multiple harmonic values \cite{J4}}

We are led to make certain infinite sequences of numbers $\har_{p^{\alpha}}(w)$ into a notion of their own, studied in the three frameworks : $\Sigma$ by their definition, $\int_{0}^{z<<1}$ by equation (\ref{eq:above}) and $\int_{0}^{1}$ by equation (\ref{eq: inv DR sigma}) ; this will shed light on the notion of finite multiple zeta values, by \S\ref{paragraph def of finite}, \S\ref{paragraph applications to finite} and Fact \ref{reduction to finite multiple zeta values}.
\footnote{In \cite{Ro}, Rosen has considered the sequences $(\har_{p}(w))_{p \in \mathcal{P}}$, and made them into a notion in the complete topological ring $\varprojlim (\prod_{p\in \mathcal{P}} \mathbb{Z}/p^{n}\mathbb{Z}) / (\oplus_{p\in \mathcal{P}} \mathbb{Z}/p^{n}\mathbb{Z}) = (\prod_{p\in \mathcal{P}} \mathbb{Z}_{p}) / \overline{\oplus_{p\in \mathcal{P}}\mathbb{Z}_{p}}$, where the bar refers to the adherence with respect to the uniform topology on $\prod_{p\in \mathcal{P}} \mathbb{Z}_{p}$ relative to the $p$-adic topologies. We denote this ring by $\mathbb{Z}_{p \rightarrow \infty} = \varprojlim \mathbb{Z} / _{(p \rightarrow \infty)^{n}}$. The proofs in \cite{Ro} do not use $\pi_{1}^{\un}(\mathbb{P}^{1} - \{0,1,\infty\})$ (they are in what we called the $\Sigma$ setting in \S3.1.3). We prove in \cite{J4} that the "asymptotic reflexion theorem" and the "asymptotic duality theorem" of \cite{Ro} are equivalent certain particular cases of our main theorems of \cite{J4}, which concern the standard families of algebraic relations : double shuffle relations, associator relations and Kashiwara-Vergne relations.} We are also led to consider instrinsically the coefficients of  $\Phi_{p,\alpha}^{-1}e_{1}\Phi_{p,\alpha}$ rather than those of $\Phi_{p,\alpha}$, because they are the ones directly involved in the explicit formulas of part I such as equation (\ref{eq: RT sigma}) and (\ref{eq: inv DR sigma}) ; we will call them adjoint $p$-adic multiple zeta values. Both these notions are intermediates for the explicit study of algebraic relations of $p$-adic multiple zeta values, but we also develop them for their intrinsic interest.

\subsubsection{Multiple harmonic values}

We define several types of infinite sequences of $\har_{p^{\alpha}}(w)$ which we intend to view as "periods" comparable to multiple zeta values ; the notations are the one of \S\ref{paragraph notations for harmonic sums} :

\begin{Definition} \label{first harmonic}For any word $w=(n_{d},\ldots,n_{1})$ :
\newline i) for each $p \in \mathcal{P}$, we call $p$-adic multiple harmonic value the sequence $\har_{p^{\mathbb{N}^{\ast}}}(w)$ 
\newline ii) for each $\alpha \in \mathbb{N}^{\ast}$, we call adelic multiple harmonic value the sequence  $\har_{\mathcal{P}^{\alpha}}(w)$
\newline iii) we call $\mathcal{P}$-adic adelic multiple harmonic value the sequence  $\har_{\mathcal{P}^{\mathbb{N}^{\ast}}}(w)$  
\end{Definition}

\noindent We will view i) as an explicit substitute to $p$-adic multiple zeta values, ii) as the canonical lift of finite multiple zeta values (via Fact \ref{reduction to finite multiple zeta values}), and iii) as the "product" of i) by ii), which provides the natural way to state the properties of i) and ii) in a unified way.

\subsubsection{Adjoint multiple zeta values}

In the framework $\int_{0}^{1}$, the study of multiple harmonic values via equation (\ref{eq:corollary I 2 a}) will be equivalent to the study of a slightly different object. Indeed, the known algebraic relations of ($p$-adic) multiple zeta values are homogeneous for the weight, and it is conjectured that it is true for all algebraic relations (Conjecture \ref{conjecture Zagier for multiple zeta values}, Conjecture \ref{conjecture p-adic multiple zeta values}) ; separating the terms of different weights in equation (\ref{eq:corollary I 2 a}) leads to :

\begin{Definition} Let $(p,\alpha) \in \mathcal{P} \times \mathbb{N}^{\ast}$,  $(n_{d},\ldots,n_{1})$ and $b \in \mathbb{N}$.
	\newline i) \label{def adjoint}We call adjoint $p$-adic multiple zeta value (Ad$p$MZV) the number 
	$$ \zeta^{\Ad}_{p,\alpha}
	\big(b; n_{d},\ldots,n_{1} \big) = (-1)^{d}  ({\Phi_{p,\alpha}}^{-1}e_{1}\Phi_{p,\alpha})
	\big[ e_{0}^{b}e_{1}e_{0}^{n_{d}-1}e_{1}\ldots e_{0}^{n_{1}-1}e_{1} \big] \in \mathbb{Q}_{p} $$
\noindent ii) Let $\Lambda$ be a formal variable. We call $\Lambda$-adic adjoint $p$-adic multiple zeta values ($\Lambda$Ad$p$MZV) the power series
$$ \zeta^{\Lambda,\Ad}_{p^{\alpha}}\big( n_{d},\ldots,n_{1} \big) = (-1)^{d} (\Phi_{p,\alpha}^{-1}e_{1}\Phi_{p,\alpha})
\big[ \frac{\Lambda^{n_{d}+\ldots+n_{1}}}{1-\Lambda e_{0}} e_{1} e_{0}^{n_{d}-1}e_{1} \ldots e_{0}^{n_{1}-1}e_{1} \big] \in \mathbb{Q}_{p}[[\Lambda]] $$
\end{Definition}

\noindent We can define in the same way complex and motivic analogues of these notions, with notations $\zeta^{\Ad}$,
$\zeta_{\text{mot}}^{\Ad}$, $\zeta^{\Lambda,\Ad}$,   $\zeta_{\text{mot}}^{\Lambda,\Ad}$, respectively where $\zeta_{\mot}$ refers to Goncharov's motivic multiple zeta values \cite{Go1} \footnote{They satisfy $\zeta_{\mot}(2)=0$ which is convenient for our purposes, contrarily to the other notions of motivic multiple zeta values defined later}. The numbers  $\zeta_{p,\alpha}^{\Ad}(0;n_{d},\ldots,n_{1})$,  which we will denote more simply by $\zeta_{p,\alpha}^{\Ad}(n_{d},\ldots,n_{1})$, are conjecturally equivalent to finite multiple zeta values by Kaneko-Zagier's Conjecture \ref{conjecture Kaneko-Zagier}, and have additional properties of their own ; their reduction modulo infinitely large primes is equal to finite multiple zeta values by Akagi-Hirose-Yasuda's equation (\ref{eq: reduction modulo primes}). The equation 
\begin{equation} \label{eq:reduction to adjoint multiple zeta values}
\zeta_{p,\alpha}^{\Lambda,\Ad}(n_{d},\ldots,n_{1}) \equiv \zeta^{\Ad}(n_{d},\ldots,n_{1}) \mod \Lambda^{n_{d}+\ldots+n_{1}+1}
\end{equation}
\noindent is the counterpart for adjoint multiple zeta values of the Fact \ref{reduction to finite multiple zeta values} which can be rewritten in terms of multiple harmonic values.

\subsubsection{From pro-affine schemes to affine formal schemes over $\mathbb{Q}[[\Lambda]]$}

Whereas the polynomial equations of periods considered earlier are homogeneous for the weight, their analogues for multiple harmonic values which we will find will be in general convergent infinite sums of algebraic relations of $p$-adic multiple zeta values, involving sequences $(\har_{p^{\alpha}}(w_{n}))_{n \in \mathbb{N}^{\ast}}$ or $(\zeta^{\Ad}_{p,\alpha}(w_{n}))_{n \in \mathbb{N}^{\ast}}$ with $\displaystyle \weight(w_{n}) \rightarrow_{n \rightarrow \infty} \infty$ which we will call completed algebraic relations. A family of completed algebraic relations defines an affine formal scheme over $\mathbb{Q}[[\Lambda]]$ equipped with the $\Lambda$-adic topology, $\Lambda$ being a formal variable ; more precisely, a sequence $(T_{n})_{n \in \mathbb{N}}$, such that, for all $n \in \mathbb{N}$, $T_{n}$ is an affine scheme over $\mathbb{Q}[[\Lambda]]/(\Lambda^{n})$, with $T_{n} \equiv T_{n+1} \mod \Lambda^{n}$.
\newline There is a faithful functor $\widehat{\tau(\Lambda)}$ sending projective limits of affine schemes over $\mathbb{Q}$ such as the fibers of $\pi_{1}^{\un}(\mathbb{P}^{1} - \{0,1,\infty\})$ or Racinet's scheme $\DMR_{0}$, to formal schemes over $\mathbb{Q}[[\Lambda]]$, defined by multiplying each equation of weight $n$ by $\Lambda^{n}$. We will view all the pro-affine schemes over $\mathbb{Q}$ and the morphisms between them as their images by $\widehat{\tau(\Lambda)}$. A functor in the converse direction can be defined by considering the lowest weight term of the equations of affine formal schemes over $\mathbb{Q}[[\Lambda]]$ ; let us call it $\gr$ ; we have then $\gr \text{ }\circ\text{ }\widehat{\tau(\Lambda)} = \id$.

\subsubsection{Main result}

One of the main results is a description of an analogue the double shuffle relations for adjoint multiple zeta values and multiple harmonic values.
We denote by $\DMR_{\smallint_{0}^{z<<1}}^{p-sym} = \DMR_{\Sigma}^{p-sym}$ the $\mathbb{Z}[[\Lambda]]$-formal scheme of solutions to the double shuffle equations satisfied by multiple polylogarithms, or equivalently by the functions $\har_{\mathbb{N}}$, and of a certain equation encoding a basic $p$-adic property of prime weighted multiple harmonic sums, which we call the "symmetry equation" (see \cite{J4}).
\newline
\newline \textbf{Theorem 4}
\newline i) (framework $\int_{0}^{1}$) There exists an explicit pro-affine scheme
$\DMR_{0,\Ad}$ over $\mathbb{Z}$, which is a group scheme for $\circ_{\Ad}^{\smallint_{0}^{1}}$, such that the map $\Ad(e_{1})$ defines a morphism of group schemes
\newline  $(\DMR_{0},\circ^{\smallint_{0}^{1}}) \rightarrow (\DMR_{0,\Ad},\circ_{\Ad}^{\smallint_{0}^{1}})$.
\newline There exists an explicit affine $\mathbb{Z}[[\Lambda]]$-formal scheme $\DMR_{\har_{\mathcal{P}^{\mathbb{N}}}}^{\smallint_{0}^{1}}$, such that $\widehat{\tau(\Lambda)}$ sends $\DMR_{0,\Ad}$ to $\DMR^{\smallint_{0}^{1}}_{\har_{\mathcal{P}^{\mathbb{N}}}}$, and which is stable by the action of  $(\DMR_{0,\Ad},\circ_{\Ad}^{\smallint_{0}^{1}})$ by $\circ_{\har}^{\smallint_{0}^{1}}$.
\newline ii) (frameworks $\int_{0}^{z<<1}$ and $\Sigma$) There exists an explicit sequence $\DMR_{\har_{\mathcal{P}^{\mathbb{N}}}}^{\smallint_{0}^{z<<1}} = \DMR_{\har_{\mathcal{P}^{\mathbb{N}}}}^{\Sigma}$, indexed by $\mathbb{N}$, of affine $\mathbb{Z}[[\Lambda]]$-formal schemes, such that the map of "restriction to $\mathcal{P}^{\mathbb{N}}$" defines a morphism : $(\DMR^{p-sym}_{\smallint_{0}^{z<<1}} = \DMR^{p-sym}_{\Sigma}) \rightarrow \DMR_{\har_{\mathcal{P}^{\mathbb{N}}}}^{\smallint_{0}^{z<<1}} = \DMR_{\har_{\mathcal{P}^{\mathbb{N}}}}^{\Sigma}$.
\newline iii) (comparison) We have $\DMR_{\har_{\mathcal{P}^{\mathbb{N}}}}^{\smallint_{0}^{1}} = \DMR_{\har_{\mathcal{P}^{\mathbb{N}}}}^{\smallint_{0}^{z<<1}}=\DMR_{\har_{\mathcal{P}^{\mathbb{N}}}}^{\Sigma}$.
\newline
\newline By i), $\Ad_{\Phi_{p,\alpha}}(e_{1})$ which is the non-commutative generating series of adjoint multiple zeta values (\S5.1.2) is a point of $\DMR_{0,\Ad}$ and $\Lambda\Ad_{\Phi_{p,\alpha}}(e_{1})$, the non-commutative generating series of $\Lambda$-adic adjoint multiple zeta values, is a point of  $\DMR_{\har_{\mathcal{P}^{\mathbb{N}}}}^{\smallint_{0}^{1}}$. By ii), $\har_{\mathcal{P}^{\mathbb{N}}}$, the non-commutative generating series of multiple harmonic values, is a point of $\DMR_{\har_{\mathcal{P}^{\mathbb{N}}}}^{\smallint_{0}^{z<<1}}=\DMR_{\har_{\mathcal{P}^{\mathbb{N}}}}^{\Sigma}$. We note that in this case, the results for $\Sigma$ and $\smallint_{0}^{z<<1}$ are not different from each other, whereas comparing the $\int_{0}^{1}$ and $\smallint_{0}^{z<<1}$ results requires a proof.
\newline 
\newline We can also define a variant $\DMR_{\har_{\mathcal{P}^{\mathbb{N}}}^{\dagger}}$ of $\DMR_{\har_{\mathcal{P}^{\mathbb{N}}}}$ which admits as a point the non-commutative generating series of the numbers $\har_{p^{\alpha}}^{\dagger_{p,\alpha}}(w) =(p^{\alpha})^{\weight(w)}\Li_{p,\alpha}^{\dagger}[w][z^{p^{\alpha}}]$, which are equal to the remainders in the series (\ref{eq:corollary I 2 a}).
\newline 
\newline \textbf{Corollary 4}
\newline The pro-affine schemes $\DMR_{\zeta_{\mathbb{Z}/_{p \rightarrow \infty}}}^{\smallint_{0}^{1}} = \gr(\DMR_{\har_{\mathcal{P}^{\mathbb{N}}}}^{\smallint_{0}^{1}})$ and $\DMR_{\zeta_{\mathbb{Z}/_{p \rightarrow \infty}}}^{\smallint_{0}^{z<<1}}= \gr(\DMR_{\har_{\mathcal{P}^{\mathbb{N}}}}^{\smallint_{0}^{z<<1}})$ are equal to each other and admit as points the non-commutative generating series of the numbers "$\zeta^{\Ad}(n_{d},\ldots,n_{1}) \mod \zeta(2)$", their $p$-adic and motivic analogues (\S5.1.2), and the non-commutative generating series of finite multiple zeta values.
\newline 
\newline This follows from Remark \ref{reduction to finite multiple zeta values}, equation (\ref{eq:reduction to adjoint multiple zeta values}) and equation (\ref{eq: reduction modulo primes}).

\subsubsection{Strategy of proof and explicit formulas}

Let, for all $d \in \mathbb{N}^{\ast}$ and $n_{1},\ldots,n_{d} \in \mathbb{N}^{\ast}$, $\shft_{\ast} (e_{0}^{n_{d}-1}e_{1} \ldots e_{0}^{n_{1}-1}e_{1}) = \frac{e_{0}^{n_{d}-1}}{(1+e_{0})^{n_{d}}}e_{1} \ldots \frac{e_{0}^{n_{1}-1}}{(1+e_{0})^{n_{1}}}e_{1}$ and 
$S_{Y}( e_{0}^{n_{d}-1}e_{1} \ldots e_{0}^{n_{1}-1}e_{1}) = (-1)^{n_{d}+\ldots+n_{1}}e_{0}^{n_{1}-1}e_{1}\ldots e_{0}^{n_{d}-1}e_{1}$
\newline 
\newline 
i) The quasi-shuffle part of the result is a purely technical checking of a commutation between the quasi-shuffle product $\ast$ and $\widehat{\tau(\Lambda)}$ resp. $\Ad(e_{1})$. The shuffle part of the result requires to guess the equation what we must prove : what we find is a formula involving $\har_{p^{\alpha}}^{\Lambda}$ equivalent to saying that $\Phi^{-1}_{p,\alpha}e_{1}\Phi_{p,\alpha}$ satisfies the shuffle equation modulo products : for all words $w,w'$, $(\Phi^{-1}_{p,\alpha}e_{1}\Phi_{p,\alpha})[w\text{ } \sh \text{ }w']=0$.
\footnote{Let $\Delta_{\sh}$ the "shuffle coproduct", defined as the continuous multiplicative map
	$\mathbb{Q}_{p}\langle \langle e_{0},e_{1}\rangle\rangle \rightarrow \mathbb{Q}_{p}\langle \langle e_{0},e_{1}\rangle\rangle \otimes \mathbb{Q}_{p}\langle \langle e_{0},e_{1}\rangle\rangle$ satisfying 
	$\Delta_{\sh}(e_{i})= 1 \otimes e_{i} + e_{i} \otimes 1$ for $i=0,1$ ; for $f \in \mathbb{Q}_{p}\langle \langle e_{0},e_{1}\rangle\rangle$, the shuffle equation amounts to $\Delta_{\sh}(f) = f \otimes f$, and the shuffle equation modulo products amounts to $\Delta_{\sh}(f) = 1 \otimes f + f \otimes 1$ ; whence $\Phi_{p,\alpha}^{-1}e_{1}\Phi_{p,\alpha}$ satisfies the shuffle equation modulo products} $\DMR_{\har_{\mathcal{P}^{\mathbb{N}}}}^{\smallint_{0}^{1}}$ is defined by the quasi-shuffle equation (for all words $w,w'$, $h(w)h(w')=h(w \ast w')$) and for all words $w,w'$ and $n \in \mathbb{N}^{\ast}$,
\begin{equation} \label{eq:first shuffle} h\big( (e_{0}^{n-1}e_{1}w) \text{ }\sh\text{ } w') \big) = h \big( w \text{ }\sh\text{ } (-1)^{n}\shft_{\ast}(e_{0}^{n-1}e_{1}w') \big) 
\end{equation}
\noindent The formula comparing the two regularizations (integral and series) of $\Phi_{p,\alpha}$ \cite{Ra} and the vanishing of the coefficients $\Phi_{p,\alpha}[e_{0}^{2n-1}e_{1}]$ imply that the two regularizations of $\Phi^{-1}_{p,\alpha}e_{1}\Phi_{p,\alpha}$ are equal to each other. One can read the equations of $\DMR_{0,\Ad}$ through the ones of $\DMR_{\har_{\mathcal{P}^{\mathbb{N}}}}^{\smallint_{0}^{1}}$ and the definition of $\zeta_{p^{\alpha}}^{\Lambda \Ad}$.
\newline The compatibilities involving $\circ^{\smallint_{0}^{1}}$, $\circ_{\Ad}^{\smallint_{0}^{1}}$ and $\circ_{\har}$ are obtained by showing that all the maps and equations necessary for Racinet's proof in \cite{Ra} can appear as lines of commutative diagrams involving the maps $\Ad(e_{1})$ and $\widehat{\tau(\Lambda)}$.
\newline ii) is proved by the usual strategy of the framework $\int_{0}^{z<<1}$ (explained in \S3.1.1) applied to the case of the double shuffle relations of multiple polylogarithms, combined to the fact that, for $m \in \{1,\ldots,p^{\alpha}\}$, we have $v_{p}(m)< v_{p}(p^{\alpha})$, thus $v_{p}(\frac{p^{\alpha}}{m}) > 0$, thus, for all $n \in \mathbb{N}^{\ast}$, $(1-\frac{p^{\alpha}}{m})^{-n} = \sum_{l\in \mathbb{N}} {-n \choose l} \big( \frac{p^{\alpha}}{m}\big)^{n}$.
$\DMR_{\har_{\mathcal{P}^{\mathbb{N}}}}^{\Sigma}$ is defined by the quasi-shuffle equation and, for all words $w,w'$,
\begin{equation} \label{eq:second shuffle} h(w\text{ }\sh\text{ }w') = h ( \shft_{\ast}(S_{Y}(w)w'))
\end{equation}
\noindent iii) We prove by that (\ref{eq:first shuffle}) and (\ref{eq:second shuffle}) are equivalent by analyzing the parameters they depend on : we prove a triple equivalence involving these equations and the following one : for all words, $u,w,w'$,
\begin{equation} h (u w\text{ }\sh\text{ }w') = h( w\text{ }\sh\text{ }\shft_{\ast}(S_{Y}(u)w')) 
\end{equation}

\subsection{(II-2) From properties of multiple harmonic sums to those of adjoint $p$-adic multiple zeta values \cite{J5}}

We now want to show that algebraic properties of adjoint $p$-adic multiple zeta values can be understood via the formulas of part I, and find converse to some arrows of \S5.1.

\subsubsection{The implicit and explicit points of view}

The question can be tackled using two different parts of the formulas, corresponding to two different points of view :
\newline 1) (By the implicit formula) Adjoint $p$-adic multiple zeta values are characterized, by equation (\ref{eq:DR RT har}), as the coefficients of the unique element of $\Ad_{\tilde{\Pi}_{1,0}(\mathbb{Q}_{p})_{S}}(e_{1})$ whose action by $\circ_{\har}^{\smallint_{0}^{z<<1}}$ sends $\har_{\mathbb{N}}$ to $\har_{p^{\alpha}\mathbb{N}}$ ; i.e. they are characterized as the coefficients of the $\int_{0}^{z<<1}$-harmonic Frobenius $(\phi^{\alpha})^{\smallint_{0}^{z<<1}}$. 
\newline 2) (By the explicit formula) Adjoint $p$-adic multiple zeta values are given explicitly in terms of multiple harmonic values by the formula (\ref{eq: RT sigma}).
\newline The first point is natural because it keeps track of the fact that the Frobenius is an isomorphism $\pi_{1}^{\un,\DR}(X_{\mathbb{Q}_{p}}) \simlra \pi_{1}^{\un,\DR}(X_{\mathbb{Q}_{p}}^{(p)})$, whereas the second point of view is natural because it is explicit.

\subsubsection{Main results}

Let $\text{O}_{\har_{\mathbb{N}}}^{\smallint_{0}^{z<<1}}$ be the orbit of $\har_{\mathbb{N}}$ for the action of $(\Ad_{\tilde{\Pi}_{1,0}(\mathbb{Q}_{p})_{S}}(e_{1}),\circ^{\smallint_{0}^{1}}_{\Ad})$ by $\circ_{\har}^{\smallint_{0}^{z<<1}}$ ; since $\circ_{\har}^{\smallint_{0}^{z<<1}}$ is free (Theorem 2), $\text{O}_{\har_{\mathbb{N}}}^{\smallint_{0}^{z<<1}}$ is a torsor under this group action. We call adjoint quasi-shuffle relation the variant of the quasi-shuffle relation which is part of the equations defining the scheme $\DMR_{0,\Ad}$ of II-1.
\newline 
\newline \textbf{Theorem 5}
\newline i) (point of view of the implicit formula) Let $h \in \text{O}_{\har_{\mathbb{N}}}^{\smallint_{0}^{z<<1}}$, and $g \in \Ad_{\tilde{\Pi}_{1,0}(\mathbb{Q}_{p})_{S}}(e_{1})$ such that $h$ and $\tilde{h} = g \circ_{\har}^{\smallint_{0}^{z<<1}} h$ satisfy the quasi-shuffle relation. Then $g$ satisfies the adjoint quasi-shuffle relation.
\newline ii) (point of view of the explicit formula) The map $\comp^{\smallint \leftarrow \Sigma}$ sends a solution to the quasi-shuffle relation to a solution to the adjoint quasi-shuffle relation.
\newline Moreover, we have somewhat analogous statements for the integral shuffle relations.
\newline 
\newline Bringing together i) resp. ii) with equation (\ref{eq:DR RT har}), resp. equation (\ref{eq: RT sigma}), both from Theorem 2, and with the fact that multiple harmonic sums satisfy the quasi-shuffle relation, reproves in two different ways that $\Phi_{p,\alpha}^{-1}e_{1}\Phi_{p,\alpha}$ satisfies the adjoint quasi-shuffle relation. This gives an answer to an "adjoint variant" of the question of Deligne and Goncharov of \S1.4.2.

\subsubsection{Strategy of proof}

We write the strategy of proof for the quasi-shuffle part.
\newline 
i) We write the quasi-shuffle relation for $\tilde{h}$ : $\tilde{h}(w)\tilde{h}(w') = \tilde{h}(w \ast w')$, and we translate it in terms of $g$ and $h$, via the equation $\tilde{h} = g \circ_{\har}^{\smallint_{0}^{z<<1}} h$. 
\newline When we encounter a product $h(z')h(z'')$, we linearize it by the quasi-shuffle relation of $h$, writing it as $h(z' \ast z'')$. It remains a linear relation between the maps $h(w) : m \in \mathbb{N} \mapsto h_{m}(w) \in \mathbb{Q}_{p}$, with coefficients in a ring of power series of $m$, which are actually overconvergent over $\mathbb{Z}_{p}$, by the bounds of valuations of $g$.
\newline By a lemma of \"{U}nver (\cite{U4}, Proposition 2.0.3) the maps $\har_{\mathbb{N}}(w)$ are linearly independent over the ring of overconvergent analytic functions of $m \in \mathbb{Z}_{p}$. For any $g \in \Ad_{\tilde{\Pi}_{1,0}(\mathbb{Q}_{p})_{S}}(e_{1})$, the map $f \mapsto g \circ_{\har}^{\smallint_{0}^{z<<1}} f$ is a pro-unipotent linear invertible operator, and it thus preserves the property of linear independence above.
\newline Thus for any $h$ in the orbit of $\har_{\mathbb{N}}$ for $\circ_{\har}^{\smallint_{0}^{z<<1}}$, all coefficients of the relation above are zero. The vanishing of the coefficients of $\har_{m}(\emptyset)=1$ is the adjoint quasi-shuffle relation for $\Phi_{p,\alpha}^{-1}e_{1}\Phi_{p,\alpha}$.
\newline ii) The quasi shuffle relation, which is true for multiple harmonic sums $\har_{m}(n_{d},\ldots,n_{1})$, remains true for the localized multiple harmonic sums $\har_{m}(\pm n_{d},\ldots,\pm n_{1})$ defined in \S\ref{the proof of ii in I-2} ; this combined with the formula for the localization map (\S\ref{the proof of ii in I-2}) yields a version of the quasi-shuffle relation for the coefficients $\mathcal{B}$ of the map of elimination of positive powers.
\newline On the other hand, $\circ_{\har}^{\Sigma}$ satisfies certain properties of symmetry with respect to the variables $(n_{d},\ldots,n_{1})$, and more generally, to certain combinatorial indices subjacent to the explicit formula for $\circ_{\har}^{\Sigma}$ written in \cite{J2}.
\newline Bringing together these two facts give the result.

\begin{Example} In depth (1,1), we have the following incarnation of the quasi-shuffle relation for the coefficients of the elimination of positive powers : $ \mathcal{B}_{b}^{l_{2}+l_{1}} + \mathcal{B}_{b}^{l_{2},l_{1}} + \mathcal{B}_{b}^{l_{1},l_{2}} = \sum_{\substack{b',b'' \geq 0 \\ b'+b''=b}} \mathcal{B}_{b'}^{l_{1}}\mathcal{B}_{b''}^{l_{2}}$, for $b,l_{1},l_{2} \in \mathbb{N}$ such that $1 \leq b \leq l_{1}+l_{2}+2$. 
\newline On the other hand, one can see on equation (\ref{eq: har RT RT in depth 2}) of Example \ref{example of RT RT har}, which expresses the formula for $\circ_{\har}^{\Sigma}$ in depth $2$, that certain lines of the formula are unchanged, resp. multiplied by $-1$ when we exchange $n_{2}$ and $n_{1}$.
\newline These two facts brought together, and combined with the other formulas of Example \ref{example DR RT depth two}, Example \ref{example of RT RT har} for $\circ_{\har}^{\smallint_{0}^{z<<1}}$ and $\circ_{\har}^{\Sigma}$ in depth one and two imply the adjoint quasi-shuffle relation in depth $(1,1)$ for $\Phi_{p,\alpha}^{-1}e_{1}\Phi_{p,\alpha}$ :
$$ (\Phi_{p,\alpha}^{-1}e_{1}\Phi_{p,\alpha})[e_{0}^{b}e_{1}e_{0}^{n_{2}-1}e_{1}e_{0}^{n_{1}-1}e_{1}\text{ }+\text{  }e_{0}^{b}e_{1}e_{0}^{n_{1}-1}e_{1}e_{0}^{n_{2}-1}e_{1}\text{ }+\text{ }e_{0}^{b}e_{1}e_{0}^{n_{2}+n_{1}-1}e_{1}] =$$
$$ \sum_{\substack{b',b''\geq 0\\ b'+b''=b}}(\Phi_{p,\alpha}^{-1}e_{1}\Phi_{p,\alpha})[e_{0}^{b'}e_{1}e_{0}^{n_{1}-1}e_{1}] \times
	(\Phi_{p,\alpha}^{-1}e_{1}\Phi_{p,\alpha})[e_{0}^{b''}e_{1}e_{0}^{n_{1}-1}e_{1}] $$
\end{Example}

\subsection{(II-3) Multiple harmonic values viewed as periods \cite{J6}}

We now formalize the idea, implicit in \S5.1 and \S5.2, that the Galois theory of periods can be applied to multiple harmonic values and $\Lambda$-adic adjoint $p$-adic multiple zeta values as if they were periods. We distinguish again the three frameworks $\int_{0}^{1}$, $\int_{0}^{z<<1}$ and $\Sigma$.

\subsubsection{In the framework $\int_{0}^{1}$ : motivic multiple harmonic values and a period conjecture}

The Galois theory of multiple zeta values is formulated by means of the notions of motivic multiple zeta values and period map and the conjecture of injectivity of the period map. We are going to see that these have analogues for multiple harmonic values.
\newline Let $\mathcal{Z}^{\mot}$ be the $\mathbb{Q}$-vector space generated by motivic multiple zeta values of \cite{Go1}. For all $n \in \mathbb{N}^{\ast}$, let $\mathcal{Z}_{n}^{\mot}$ be the subspace generated by elements of weight $n$. The weight is a grading on $\mathcal{Z}^{\mot}$, namely we have  $\mathcal{Z}^{\mot} = \oplus_{n\in \mathbb{N}} \mathcal{Z}^{\mot}_{n}$. 
\newline Let $\widehat{\mathcal{Z}^{\mot}}$ be the completion of $\mathcal{Z}^{\mot}$ with respect to the decreasing weight filtration whose $n$-th term is $\oplus_{n'\geq n}\mathcal{Z}^{\mot}_{n'}$ for all $n \in \mathbb{N}$. Since the weight on $\mathcal{Z}^{\mot}$ is a grading, we have a canonical isomorphism $\widehat{\mathcal{Z}^{\mot}} \simeq \prod_{n\in\mathbb{N}} \mathcal{Z}_{n}^{\mot}$. This enables to represent an element of $\widehat{\mathcal{Z}^{\mot}}$ as a formal sum $\sum_{n\in \mathbb{N}} \zeta^{\mot}(w_{n})$, where $w_{n}$ is a $\mathbb{Q}$-linear combination of words of weight $n$, in such a way that we have the equivalence $\sum_{n\in \mathbb{N}} \zeta^{\mot}(w_{n})=0$ $\Leftrightarrow$ for all $n \in \mathbb{N}$, $\zeta^{\mot}(w_{n})=0$. We will use implicitly this representation below. Let $\Phi_{\mot} = \sum_{w\text{ word}} \Phi_{mot}[w]w$.

\begin{Definition}
For any $n_{d},\ldots,n_{1} \in \mathbb{N}^{\ast}$, we call motivic multiple harmonic value and, at the same time, motivic $\Lambda$-adic adjoint multiple zeta value the following infinite sequence of motivic multiple zeta values (to be compared with equation (\ref{eq:corollary I 2 a}))
$$ \har_{\mathcal{P}^{\mathbb{N}}}^{\mot}(n_{d},\ldots,n_{1}) = 
(\Phi_{\mot}^{-1}e_{1}\Phi_{\mot})\big[\frac{1}{1-e_{0}}e_{1}e_{0}^{n_{d}-1}e_{1}\ldots e_{0}^{n_{1}-1}e_{1} \big] \in \widehat{\mathcal{Z}^{\mot}} $$
\end{Definition}

\noindent In the next statement, we assume that Yamashita's (currently unavailable) work on the crystalline realization of mixed Tate motives implies, for any $(p,\alpha) \in \mathcal{P} \times (\mathbb{Z}\cup \{+\infty\}-\{0\})$, the existence of a period map sending $\zeta^{\mot}(w) \mapsto \zeta_{p,\alpha}(w)$ for all words $w$ - to our knowledge, this is known at least for certain values of $\alpha$.

\begin{Definition}
\noindent 
\newline i) Let $\widehat{\mathcal{Z}_{\conv}^{\mot}} \subset \widehat{\mathcal{Z}^{\mot}}$ be the set of elements $\sum_{n \in \mathbb{N}} \zeta^{\mot}(w_{n})$, such that, for all $(p,\alpha) \in \mathcal{P} \times \mathbb{N}^{\ast}$, $\sum_{n\in\mathbb{N}} \zeta_{p,\alpha}(w_{n})$ converges in $\mathbb{Q}_{p}$.
\newline ii) For each $I \subset\mathcal{P}\times\mathbb{N}$, let
$\widehat{\text{per}}_{\zeta_{I}} : \sum_{n \in \mathbb{N}} \zeta^{\mot}(w_{n}) \in \widehat{\mathcal{Z}_{\conv}^{\mot}} \mapsto (\sum_{n \in \mathbb{N}} \zeta_{p,\alpha}(w_{n}))_{(p,\alpha)\in I}$.
\newline iii) For each $I \subset\mathcal{P}\times\mathbb{N}$, let $\widehat{\text{per}}_{\har_{I}}$ be the restriction of $\widehat{\text{per}}_{\zeta_{I}}$ to the subset of weight-adically convergent sums of motivic multiple harmonic values.
\end{Definition}

\noindent The maps of ii) are well-defined by the fact that the weight is a grading on $\mathcal{Z}^{\mot}$ and from the ensuing description of $\widehat{\mathcal{Z}^{\mot}}$ explained above. The maps of iii) are well-defined since $\widehat{\mathcal{Z}_{\conv}^{\mot}}$ contains the motivic multiple harmonic values

\begin{Conjecture} If $I$ contains a set of the form $p^{\mathbb{N}^{\ast}}$ with $p$ a prime number, or  $\mathcal{P}^{\alpha}$, with $\alpha \in \mathbb{N}^{\ast}$, then $\widehat{\text{per}}_{\har_{I}}$ is injective. In other terms, any relation expressing the vanishing of a weight-adically and $p$-adically convergent sum of multiple harmonic values lifts to an infinite collection of equalities among motivic multiple zeta values.
\end{Conjecture}

\noindent This setting can be essentialized by defining a notion of "summable subgroupoid of $\pi_{1}^{\un,\DR}(\mathbb{P}^{1} - \{0,1,\infty\})$" : the details are written in \cite{J6}.

\subsubsection{In the framework $\smallint_{0}^{z<<1}$ : the "Taylor periods" of $\pi_{1}^{\un}(\mathbb{P}^{1} - \{0,1,\infty\})$}

The framework of computation which involves power series expansion of multiple polylogarithms can also be essentialized by a conjecture ; however, we have to involve not only $\mathbb{P}^{1} - \{0,1,\infty\}=\mathcal{M}_{0,4}$ but also at least the moduli space $\mathcal{M}_{0,5}$ since, for instance, the quasi-shuffle relation is a property of multiple polylogarithms on $\pi_{1}^{\un}(\mathcal{M}_{0,5})$ ; it is actually convenient to involve all the spaces $\mathcal{M}_{0,n}$.

\begin{Conjecture} Any relation expressing the vanishing of an $p$-adically and weight-adically convergent sum of multiple harmonic values is a consequence of equations obtained from equations on multiple polylogarithms on $\mathcal{M}_{0,n}$, $n \geq 4$, by taking power series expansions, and from the "$p$-adic symmetry equation" of prime weighted multiple harmonic sums evoked in \S5.2 and defined in II-1.
\end{Conjecture}

\subsubsection{In the framework $\Sigma$ : a conjecture à la Kontsevich-Zagier for multiple harmonic values viewed as elementary iterated sums}

\noindent Let us define a counterpart for multiple harmonic values viewed as finite iterated sums of the elementary framework for periods stated by Kontsevich and Zagier (Conjecture \ref{conjecture of periods Kontsevich Zagier}).

\begin{Conjecture}
Let us consider
\newline i) the subsets of $\amalg_{d\in \mathbb{N}}(\mathbb{N}^{\ast})^{d+1}$ obtained by taking finite disjoint unions and finite products of subsets of the form $\{(m_{1},\ldots,m_{d},m=m_{d+1}) \in (\mathbb{N}^{\ast})^{d+1} \text{ }|\text{ }0 \ast_{0} m_{1} \ast _{1}\ldots \ast_{d-1} m_{d} \ast_{d} m\}$, with $\ast_{i} \in \{=,<\}$ for all i, 
\newline ii) the $\mathbb{Q}$-algebras of functions $(\mathbb{N}^{\ast})^{d+1} \rightarrow \mathbb{Q}$ generated by the maps sending $(m_{1},\ldots,m_{d+1})$ to $m_{i}^{t_{i}}$ with $t_{i} \in \mathbb{Z}$, binomial coefficients of affine functions of the $m_{i}'s$ (including of the form ${u \choose u'}$ with $u<0$), factorials of the $m_{i}$'s. 
\newline iii) the $\mathbb{Q}$-algebra generated by numbers $\sum_{x \in D} f(x)$ with $f$ as in ii) and $D$ as in i), viewed as elements of $\mathbb{Q}_{p}$.
\newline Then, any equation expressing the vanishing of an $p$-adically and weight-adically convergent sum of $\mathbb{Q}$-linear combinations of multiple harmonic values, can be obtained using i) ii) iii) and the properties of the following type :
\newline iv) equalities among domains of summation, equations satisfied by summands, equations of changes of variables
\newline v) operations arising from the structure of complete topological field of $\mathbb{Q}_{p}$ (such as expanding $\frac{1}{1-x}= \sum_{l\in \mathbb{N}} x^{l}$ in $\mathbb{Q}_{p}$ for $|x|_{p}<1$).
\end{Conjecture}

\subsubsection{Conclusion}

In the end, what we can call "explicit $p$-adic periods" is the couple formed by ($\Lambda$-adic) adjoint $p$-adic multiple zeta values and multiple harmonic values, with the usual conjecture of periods for $p$-adic multiple zeta values, the conjectures of this \S5.3, and the equations relating adjoint $p$-adic multiple zeta values and multiple harmonic sums from part I.

\section{Conclusion}

Let us summarize what seem to us to be the two main outcomes about $p$-adic multiple zeta values of this work. The initial question was to compute $p$-adic multiple zeta values and apply it to understand explicitly algebraic relations (\S\ref{initial question}).

\subsection*{a - About the explicitness}

The formulas for $p$-adic multiple zeta values provided by our indirect method for solving the differential equation of Frobenius are fully explicit, i.e. they can be written concisely provided the introduction of certain combinatorial tools (the full formulas, which we did not reproduced here, are written in the papers I-2 and I-3) ; they also enable to visualize explicitly algebraic relations (\S5.1-\S5.2).
\newline 
\newline What we have proved is that the Frobenius, because of bounds of the valuations of its values at the canonical paths, becomes much simpler computationally when one considers a certain limit of its parameters, and that the corresponding limit of the Frobenius, which we call harmonic Frobenius, is equivalent to a particular explicit structure on multiple harmonic sums. Moreover, this limit of the Frobenius is enough to reconstruct the whole of the Frobenius.
\newline 
\newline The explicitness is obtained by introducing and using the notion of adjoint $p$-adic multiple zeta values ; namely, for each question on $p$MZVs which we want to tackle via explicit formulas, we find the appropriate variant of this question for Ad$p$MZVs, we solve the variant, and the step of going back to $p$MZVs from adjoint $p$MZVs, if we do it, is non $p$-adic, it amounts to a problem of algebra over $\mathbb{Q}$. $p$MZVs and Ad$p$MZVs are equally natural families of values which characterize the Frobenius, and it seems to us that Ad$p$MZVs are actually the most natural from certain points of view.

\subsection*{b - About the formulation related to the motivic Galois theory of periods}

We have formulated this work in a way which keeps in mind the concepts of the motivic Galois theory of periods, despite the presence of infinite summations in the formulas.
\newline Aside from the regularization of $p$-adic iterated integrals, the central objects are the $p$-adic pro-unipotent harmonic actions $\circ_{\har}^{\smallint_{0}^{1}}$, $\circ_{\har}^{\smallint_{0}^{z<<1}}$, $\circ_{\har}^{\Sigma}$, the maps $\comp^{\Sigma \rightarrow \smallint}$, $\comp^{\smallint \rightarrow \Sigma}$ expressing the relations between them, and the maps of iterations of the harmonic Frobenius $\iter_{\har}^{\smallint_{0}^{1}}$, $\iter_{\har}^{\Sigma}$ which we have defined in part I. The $p$-adic pro-unipotent harmonic actions keep in mind the motivic Galois theory : they are related to the product $\circ^{\smallint_{0}^{1}}$ (the Ihara product or twisted Magnus product) on $\Pi_{1,0}=\pi_{1}^{\un,\DR}(\mathbb{P}^{1} -\{0,1,\infty\},-\vec{1}_{1},\vec{1}_{0})$, which is itself related to the motivic Galois action on $\Pi_{1,0}$. In part I, these maps were the way to express the formulas for the Frobenius. In part II, the central idea was that these maps "reflect algebraic relations" as motivic Galois actions do for periods. More precisely, we have defined \emph{multiple harmonic values} as certain infinite sequences of prime weighted multiple harmonic sums, and viewed them as periods, we have studied their (completed) algebraic relations, and compared them to the algebraic relations among adjoint $p$-adic multiple zeta values.
\newline 
\newline In a summary, what we have showed is that 1) there exists a computation of the Frobenius which keeps a track of the motivic Galois action, and 2) the motivic Galois theory of $p$-adic multiple zeta values can be extended to multiple harmonic values : we have an "explicit elementary version" of the motivic Galois theory of $p$-adic multiple zeta values, formulated as a comparison between adjoint $p$-adic multiple zeta values and multiple harmonic values via the $p$-adic pro-unipotent harmonic actions.

\subsection*{c - Other perspectives}

The ideas summarized above, and some variants of them, also play a central role in the three other parts of this work (parts III, IV, V : \cite{J7} to \cite{J12}), and in other forthcoming works on pro-unipotent fundamental groupoids.

\section*{Acknowledgments}

I warmly thank Hidekazu Furusho for having invited me to the conference \emph{Various aspects of multiple zeta values} held in July of 2016 in Kyoto. This text was written for the proceedings of that conference. I thank Kyoto University for hospitality and financial support during the conference. I thank Pierre Cartier for his support in the end of my PhD at Universit\'{e} Paris Diderot. I thank Benjamin Enriquez for his support during my post-doctorate in Universit\'{e} de Strasbourg. This text has been written at Institut de Recherche Mathématique Avancée of Universit\'{e} de Strasbourg, supported by Labex IRMIA, and revised at Universit\'{e} de Gen\`{e}ve, supported by NCCR SwissMAP.

\noindent Universit\'{e} de Gen\`{e}ve, Section de mathématiques, 2-4 rue du Li`{e}vre,
Case postale 64
1211 Gen\`{e}ve 4, Suisse
\newline e-mail : david.jarossay@unige.ch

\end{document}